\documentclass[11pt,leqno]{article}
\usepackage{amssymb}
\usepackage{hyperref}
\newcommand\qed{\hfill$\sqcap\kern-7.5pt\hbox{$\sqcup$}$}
\newcommand{\CC}{\mathbb{C}}
\newcommand{\NN}{\mathbb{N}}
\newcommand{\ZZ}{\mathbb{Z}}
\newcommand{\RR}{\mathbb{R}}
\newcommand{\Ll}{\mathbb{L}}
\newcommand{\R}{\mathbb{R}}
\newcommand{\N}{\mathbb{N}}
\newcommand{\Sp}{\mathbb{S}}

\newcommand{\Ww}{\mathbb{W}}
\newtheorem{theo}{Theorem}
\newtheorem{prop}[theo]{Proposition}
\newtheorem{lem}[theo]{Lemma}
\newtheorem{cor}[theo]{Corollary}
\newtheorem{rem}[theo]{Remark}
\newtheorem{rems}[theo]{Remarks}

\newcommand{\beqn}{\begin{equation}}
\newcommand{\eeqn}{\end{equation}}
\newcommand{\bear}{\begin{eqnarray}}
\newcommand{\eear}{\end{eqnarray}}
\newcommand{\bean}{\begin{eqnarray*}}
\newcommand{\eean}{\end{eqnarray*}}

\newcommand{\Sph}{\mathbb{S}}
\newcommand{\Aa}{\mathcal{A}}

\newcommand{\Cc}{\mathcal{C}}
\newcommand{\DD}{\mathcal{D}}
\newcommand{\EE}{\mathcal{E}}
\newcommand{\GG}{\mathcal{G}}
\newcommand{\HH}{\mathcal{H}}

\newcommand{\LL}{\mathcal{L}}
\newcommand{\OO}{\mathcal{O}}

\newcommand{\Rr}{\mathcal{R}}

\newcommand{\Ss}{\mathcal{S}}

\newcommand{\e}{{\varepsilon}}

\newcommand{\eps}{\varepsilon}
\newcommand{\wto}{\rightharpoonup}

\begin{document}

\title{Stability, convergence to self-similarity and elastic limit for 
the Boltzmann equation for inelastic hard spheres} 

\date{}

\author{S. {\sc Mischler}$^1$, C. {\sc Mouhot}$^2$}

\footnotetext[1]{CEREMADE, Universit\'e Paris IX-Dauphine,
Place du Mar\'echal de Lattre de Tassigny, 75775 Paris, France. 
E-mail: \texttt{mischler@ceremade.dauphine.fr}}

\footnotetext[2]{CEREMADE, Universit\'e Paris IX-Dauphine,
Place du Mar\'echal de Lattre de Tassigny, 75775 Paris, France. 
E-mail: \texttt{cmouhot@ceremade.dauphine.fr}}

\maketitle

\begin{abstract} 
We consider the spatially homogeneous Boltzmann equation for {\em inelastic hard spheres}, 
in the framework of so-called {\em constant normal restitution coefficients} $\alpha \in [0,1]$. 
In the physical regime of a small inelasticity (that is $\alpha \in [\alpha_*,1)$ for some constructive 
$\alpha_*>0$) we prove uniqueness of the self-similar profile for 
given values of the restitution coefficient $\alpha \in [\alpha_*,1)$, the mass 
and the momentum; therefore we deduce the uniqueness of the self-similar solution (up 
to a time translation). 

Moreover, if the initial datum lies in $L^1_3$, and 
under some smallness condition on $(1-\alpha_*)$ depending on the mass, energy and 
$L^1 _3$ norm of this initial datum, 
%depending on the mass,  energy and $L^2$ norm of the initial datum, 
we prove time asymptotic convergence 
(with polynomial rate) of the solution towards the self-similar solution 
(the so-called {\em homogeneous cooling state}). 

These uniqueness, stability and convergence results are expressed in the self-similar 
variables and then translate into corresponding results for the original 
Boltzmann equation.  The proofs are based on the identification of a suitable 
elastic limit rescaling, and the construction of a smooth path of self-similar profiles 
connecting to a particular Maxwellian equilibrium in the elastic limit, together 
with tools from perturbative theory of linear operators. 
Some universal quantities, such as the ``quasi-elastic self-similar temperature'' 
and the rate of convergence towards self-similarity at first order 
in terms of $(1-\alpha)$, are obtained from our study. 

These results provide a positive answer and a mathematical proof of the  
Ernst-Brito conjecture~\cite{EBJSP} in the case of inelastic hard spheres with small inelasticity. 
\end{abstract}



%\textbf{Ce qu'il reste \`a faire}: 

%- peux-tu revoir la sous-section 5.1 (il faudra \'galement relire la toute fin de la section 5)

%-  A droite dans (\ref{convLetoL1}) c'est 
%$\| f \|_{W^{3,1}_3(m^{-1})} \, \| g \|_{W^{3,1}_3(m^{-1})}$ et il faut donc suivre les modifications?

%- notation: $\alpha_5 = \alpha_{**}$?

%- revoir l'\'economie g\'en\'erale des remarques, et en particulier "1.8 Plan of the paper"
\newpage 
\tableofcontents

\vspace{0.3cm}

\section{Introduction and main results}
\setcounter{equation}{0}
\setcounter{theo}{0}

%{\bf VIRER TOUS LES "EXPLICIT CONSTANT" } 

\subsection{The model}

We consider the spatially homogeneous Boltzmann equation for hard spheres
undergoing inelastic collisions with a constant normal restitution coefficient $\alpha \in [0,1)$ 
(see~\cite{GPV**,BGP**,MMRI,MMII}).
More precisely, the gas is described by the distribution density of particles
$f= f_t = f(t,v) \ge 0$ with velocity $v \in \RR^N$ ($N \ge 2$) at time $t \ge 0$ and 
it satisfies the evolution equation
     \bear  \label{eqBol1}
     {\partial f \over \partial t}& = & Q_\alpha(f,f) \quad\hbox{ in }\quad
     (0,+\infty)\times \RR^N,\\
     \label{eqBol2}
     f(0,\cdot) & = &  f_{\mbox{\scriptsize{in}}} \quad\hbox{ in }\quad  \RR^N.
     \eear
\smallskip

The quadratic collision operator $Q_\alpha(f,f)$ models the interaction of
particles by means of inelastic binary collisions (preserving mass and momentum
but dissipating kinetic energy). We define the collision operator by its action on test functions, or 
{\it observables}. Taking $\psi = \psi(v)$ to be a suitably regular test function, 
we introduce the following weak formulation of the collision operator 
\beqn  \label{Qinel}
\int_{\R^N} 
Q_\alpha(g,f) \, \psi \, dv =  \int \!\! \int \!\!\int_{\R^N \times \R^N \times \Sp^{N-1}} 
b \, |u| \, g_* \, f \, (\psi' - \psi ) \, d\sigma \, dv \, dv_*, 
\eeqn
where we use the shorthand notations $f:= f(v)$,  $g_*:=g(v_*)$, $\psi':=\psi(v')$, etc.
Here and below $u=v-v_*$ denotes the relative velocity and $v',v'_*$ denotes the possible post-collisional velocities (which encapsule the inelasticity of the collision operator 
in terms of $\alpha$). They are defined by
\beqn\label{vprimvprim*}
\quad\quad   v' = {w \over 2} +  {u' \over 2}, \quad 
v'_*= {w \over 2} - {u' \over 2}, 
\eeqn
with
$$
w = v+v_*, \qquad u' = \left( {1- \alpha \over 2} \right) \, u + \left( {1+\alpha \over 2} \right) \, |u| \, \sigma.
$$
We also introduce the notation $\hat{x} = x/|x|$ for any $x \in \R^N$, $x\not= 0$.
The function $b=b(\hat u \cdot \sigma)$ in (\ref{Qinel}) is (up to a multiplicative factor)
the {\em differential collisional cross-section}. We assume that 
  \begin{equation}\label{hypb1}
  b \mbox{  is Lipschitz, non-decreasing and convex on } (-1,1) 
  \end{equation}
and that 
  \begin{equation}\label{hypb2}
  \exists \, b_m,b_M \in (0,\infty) \quad \mbox{ s.t. } \quad 
     \forall \, x \in [-1,1], \quad b_m \le b(x) \le b_M. 
   \end{equation}   
Note that the ``physical'' cross-section for hard spheres is given by
(see \cite{GPV**,Cerci?})
  \beqn\label{HScs} 
  b(x) = b'_0\, (1-x)^{-{N-3 \over 2}}, \quad b'_0 \in (0,\infty),
  \eeqn
so that it fulfills the above hypothesis (\ref{hypb1},\ref{hypb2}) when $N=3$. These hypothesis are 
needed in the proof of moments estimates (see~\cite[Proposition~3.2]{MMRI} 
and~\cite[Proposition~3.1]{MMII}).

We also define the symmetrized (or polar form of the) bilinear collisional operator $\tilde Q_\alpha$ by setting
 \beqn  \label{Qinelsym}
 \left\{
 \begin{array}{l}
 \displaystyle{
 \int_{\R^N} \tilde Q_\alpha(g,h) \, \psi \, dv = {1 \over 2} 
 \int \!\! \int \!\!\int_{\R^N \times \R^N \times \Sp^{N-1}} 
 b \, |u| \, g_* \, h \, \Delta_\psi \, d\sigma\,dv \, dv_*, }\vspace{0.3cm} \\
 \displaystyle{ \hbox{with} \quad \Delta_\psi = \left( \psi' + \psi'_* - \psi - \psi_*\right) }.
 \end{array}
 \right.
 \eeqn
In other words, $ \tilde Q_\alpha(g,h) =  ( Q_\alpha(g,h) +  Q_\alpha(h,g))/2$. 
The formula~(\ref{Qinel}) suggests the natural splitting 
$Q_\alpha = Q^+ _\alpha - Q^- _\alpha$ 
between gain and loss part. The loss part $Q^- _\alpha$ can be 
defined in strong form noticing that 
 $$
 \langle Q^-_\alpha(g,f) , \ \psi \rangle =  \int \!\! \int \!\!\int_{\R^N \times \R^N \times \Sp^{N-1}} 
 b \, |u| \, g_* \, f \, \psi  \, d\sigma \, dv \, dv_* =:  \langle f \, L(g) , \ \psi \rangle,
 $$
where $\langle \cdot, \cdot \rangle$ is the usual scalar product in $L^2$ and 
$L$ is the convolution operator
 \beqn\label{defL}\qquad
 L(g) (v) = (b_0 \,  |\cdot| * g )(v) = b_0 \, \int_{\R^N} g(v_*) \, |v-v_*| \, dv_*,
\,\,\,\hbox{with}\,\,\, b_0 = \int_{S^{N-1}} b(\sigma_1) \, d\sigma. 
 \eeqn
In particular note that $L$ and $Q^- _\alpha=Q^-$ are indeed independent of the 
normal restitution coefficient $\alpha$. 
\smallskip

The Boltzmann equation~(\ref{eqBol1}) is complemented with an initial
datum~(\ref{eqBol2}) which satisfies 
  \beqn  \label{initialcond}
  \left\{
  \begin{array}{l}
  \displaystyle{
  0 \le f_{\mbox{\scriptsize{in}}} \in L^1(\RR^N), \qquad
  \rho(f_{\mbox{\scriptsize{in}}}) :=  \int_{\RR^N} f_{\mbox{\scriptsize{in}}} \, dv  
  = \rho \in (0,\infty)  }\vspace{0.1cm} \\
  \displaystyle{ \int_{\RR^N} f_{\mbox{\scriptsize{in}}} \, v \, dv  = 0, \qquad
  \EE( f_{\mbox{\scriptsize{in}}}) :=    \int_{\RR^N} f_{\mbox{\scriptsize{in}}} \, |v| ^2 \, dv  < \infty     .}
  \end{array}
  \right.
  \eeqn

%{\bf ENLEVER: Notice that, without loss of generality, we can assume the two first moment
%conditions in~(\ref{initialcond}) with $\rho = 1$, since we
%may always reduce to that case by a scaling and translation argument 
%(see for instance \cite[Section 1.5]{GPV**}). However in the sequel 
%we shall indicate the dependency in terms of $\rho$ of the main formulas, 
%since it enlightens their homogeneity and therefore the 
%physical meaning of the mathematical quantities handled. } 
% CM : OK j'enlève 
\smallskip

As explained in~\cite{MMRI,MMII}, the operator~(\ref{Qinel})
preserves mass and momentum, and so does the evolution equation:
\beqn
 {d \over dt} \int_{\RR^N} f_t \,
\left( \begin{array}{ll}
       1 \\
       v
       \end{array} \right)
 \, dv = 0,
 \eeqn
while kinetic energy is dissipated
\beqn \label{eqdiffEE}
{d \over dt} \EE(f_t) = - (1-\alpha^2) \, D_\EE(f_t). 
\eeqn
The {\em energy dissipation functional} is given by
\[
D_\EE (f):= b_1 \int\!\!\int_{\RR^N \times \RR^N} f \, f_* \, |u|^3 \, dv \, dv_*, 
\]
where  $b_1$ is (up to a multiplicative factor) the angular momentum defined by 
\beqn \label{defb1}
b_1 :=  {1 \over 8} \int_{\Sp^{N-1}} \left( 1-(\hat{u}\cdot\sigma) \right) \, 
b(\hat{u}\cdot\sigma) \, d\sigma .
%=  \frac{|\Sp^{N-2}|}8 \, \int_0 ^\pi \left( 1- \cos \theta \right) \,  b(\cos \theta) \, \sin^{N-2} \theta/2 \, d\theta.  
\eeqn
In order to establish (\ref{eqdiffEE}) we have used (\ref{Qinelsym}) and the elementary computation 
$$
\Delta_{|\cdot|^2} (v,v_*,\sigma) = - {1 - \alpha^2 \over 4} \, \, (1 -(\hat{u}\cdot\sigma)) \, |u|^2.
$$

\smallskip
The study of the Cauchy theory and the cooling process
of~(\ref{eqBol1})-(\ref{eqBol2}) was
done in~\cite{MMRI}. The
equation is well-posed for instance in $L^1 _2$: for $0 \le 
f_{\mbox{\scriptsize{in}}} \in L^1 _2$,
there is a unique global solution in 
$C(\RR_+; L^1 _2) \cap L^1(\RR_+;  L^1 _3)$ (see Subsection~\ref{MMR:subsec:not} 
for the notation of functional spaces).
This solution preserves mass, momentum and
has a positive and decreasing kinetic energy. 
Moreover, as time goes to infinity, it satisfies:
\beqn  \label{asymptTc}
   \EE(t) \to 0 \qquad \hbox{and}  \qquad f(t,\cdot) \ \wto \  \delta_{v=0}
   \,\,\hbox{ in }\,\, M^1(\RR^N)\hbox{-weak *},
\eeqn
where $M^1(\RR^N)$ denotes the space of probability measures on $\RR^N$.

\subsection{Introduction of rescaled variables}

Let us introduce some rescaled variables (which  can be found in~\cite{EBphysrev,BGP**,MMII} 
for instance), in order to study more precisely the asymptotic behavior (\ref{asymptTc}) of the solution.
For any solution $f$ to the Boltzmann equation~(\ref{eqBol1}), we may associate for any $\tau \in (0,\infty)$ the  self-similar rescaled solution $g$ by the relation
$$
g(t,v) = e^{-N \, \tau \, t} \, f\left( {e^{\tau \, t}  -1 \over \tau} , e^{-\tau \, t} \, v \right).
$$
Using the homogeneity property 
$Q_\alpha(g(\lambda \cdot),g(\lambda \cdot))(v) = \lambda^{-(N+1)} \, Q_\alpha(g,g)(\lambda v)$, 
it is straightforward that $g$ satisfies the evolution equation
\beqn\label{eqresc}
  \frac{\partial g}{\partial t} = Q_\alpha(g,g) - \tau \, \nabla_v \cdot (v g). 
\eeqn
Any non-negative steady state $0 \le G = G(v)$ of~(\ref{eqresc}), that is $G$ satisfying
\beqn\label{eqprofGtau}
 Q_\alpha(G,G) - \tau \, \nabla_v \cdot (v \, G) = 0, 
\eeqn
is called a {\em self-similar profile}. It translates into a {\em self-similar solution} (or 
{\em homogeneous cooling state}) $F$ of the original equation~(\ref{eqBol1}) by setting
\beqn\label{GtoF}
F(t,v) = (V_0+\tau \, t)^N \, G( (V_0+ \tau \, t)v), 
\eeqn
for a given constant $V_0 \in (0,\infty)$. 
Reciprocally, let us consider a {\em self-similar solution}  $F$ of 
the original equation~(\ref{eqBol1}). That means a solution 
$F$ of (\ref{eqBol1}) with the specific shape 
\begin{equation}\label{shapeselfsim}
F(t,v) = V(t)^N \, G(V(t) \, v)
\end{equation}
for some given non-negative distribution $G = G(v)$ and some $C^1$, positive, increasing time rescaling function $V(t)$. One can easily show (see for instance \cite[section 1.2]{MMII}) that $V(t) = \tau \, t + V_0$ for some constants $\tau, \, V_0 >0$ and $G$ satisfies (\ref{eqprofGtau}) associated to the velocity rescaling parameter $\tau$. For a given self-similar profile $G$, associated to a velocity rescaling parameter $\tau$ and with mass $\rho$ and energy $\EE$, 
we may associate a new  self-similar profile $\tilde G$, associated to a velocity rescaling parameter 
$\tilde \tau$ and with mass $\tilde \rho$ by setting 
$$
\tilde G (v) = K \, G(V \, v), \quad V = {\tilde \rho \over \rho} \, {\tau \over \tilde \tau}, \quad
K = V^N \ {\tilde \rho \over \rho}.
$$
The energy of $\tilde G$ is then $\tilde \EE = {\tilde \rho \over \rho} \, \left( {\tau \over \tilde \tau}\right)^2 \, \EE$. We thus see that there exists a two real parameters family of self-similar profiles which can be either parametrized by $(\rho,\tau)$ or by $(\rho,\EE)$.   For fixed mass, changing the velocity rescaling parameter $\tau$ in (\ref{eqprofGtau}) corresponds to a change of the energy of the profile, or equivalently to an homothetic change of variable of the solution. Therefore it is no restriction to choose arbitrarily  this constant. Also note that modifying $V_0$ just corresponds to a time translation in the self-similar solution $F$ defined by (\ref{GtoF}). 

It was proved in~\cite[Theorem~1.1]{MMII} that for any inelastic 
parameter $\alpha \in (0,1)$, mass $\rho \in (0,\infty)$ and (thanks to the preceding discussion) 
any velocity rescaling parameter $\tau \in (0,\infty)$, there exists at least one positive and smooth 
self-similar profile $G$ with given mass $\rho$ and vanishing 
momentum:
  \begin{equation}  \label{profileeq1}
  \left\{ 
  \begin{array}{l} \displaystyle 
  Q_\alpha(G,G) - \tau \, \nabla_v \cdot (v \, G) = 0  
  \quad \hbox{ in }\quad \RR^N,\vspace{0.3cm} \\ \displaystyle 
  \int_{\R^N} G \, dv = \rho, \quad \int_{\RR^N} G \, v \, dv = 0, 
  \quad  0 <  G  \in \Ss(\R^N),
  \end{array}
  \right.    
  \end{equation}
where $\Ss(\R^N)$ denotes the Schwartz space of $C^\infty$ functions decreasing at infinity 
faster than any polynomial. 

Finally, for any solution $g$ to the Boltzmann equation in self-similar variables (\ref{eqresc}),
we may associate a solution $f$ to the evolution problem~(\ref{eqBol1}),
defining $f$ by the relation
\beqn\label{gtof}
f(t,v) = (V_0+ \tau \, t)^N \, g \left( \frac{\ln(V_0 + \tau \, t)}{\tau},(V_0+\tau \, t)v \right).
\eeqn

\subsection{Rescaled variables and elastic limit $\alpha \to 1$}

We now make the choice 
\beqn\label{taualpha}
\tau = \tau_\alpha = \rho \, (1-\alpha), 
\eeqn
and denote by  $G_\alpha$ a solution to the problem (\ref{profileeq1}). 
At a formal level, it is immediate that with this choice of scaling, 
in the elastic limit $\alpha \to1$, the equation (\ref{profileeq1}) becomes
  \begin{equation} \label{profileeq2}
  \left\{ 
  \begin{array}{l} \displaystyle 
  Q_1(G_1,G_1)  = 0  \,\,\hbox{ in }\,\,  \RR^N, \vspace{0.3cm}  \\ \displaystyle 
  \int_{\R^N} G_1 \, dv = \rho, \quad \int_{\RR^N} G_1 \, v \, dv = 0, 
     \quad 0 \le  G_1  \in \Ss(\R^N).  
  \end{array}
  \right. 
  \end{equation}
Moreover, multiplying the first equation of (\ref{profileeq1}) by $|v|^2$, integrating in 
the velocity variable as in (\ref{eqdiffEE}) and taking into account the additional term 
coming from the additional drift term in (\ref{eqresc}), one gets
  \beqn\label{EnergyDEGe}
  2 \, (1-\alpha) \, \rho \, \EE(G_\alpha) - (1-\alpha^2) D_\EE(G_\alpha)  = 0.
  \eeqn
Dividing the above equation by $(1-\alpha)$ and passing to the limit $\alpha \to 1$, one obtains 
  \beqn\label{profilee3}
  \rho \, \EE(G_1) -  D_\EE(G_1)  = 0.
  \eeqn
It is straightforward (see Proposition~\ref{ConvGexplicit} below) that the only 
function satisfying the constraints (\ref{profileeq2}) and (\ref{profilee3}) 
is the Maxwellian function
  \beqn\label{Maxlim}
  \bar G_1 := M_{\bar\theta_1} = M_{\rho, 0, \bar\theta_1} 
  \eeqn
where, for any $\rho,\theta >0$, $u \in \R^N$, the function $M_{\rho,u,\theta}$ 
denotes the  Maxwellian with mass $\rho$, momentum $u$ and temperature 
$\theta$ given  by  
\beqn\label{defMrut}
M_{\rho,u,\theta} (v) :=  
  {\rho \over (2\pi \theta)^{N/2}} \, e^{- {|v-u|^2 \over 2 \theta}},
\eeqn
and where the temperature $\bar\theta_1 \in (0,\infty)$ is given by 
(we recall that $b_1$  is defined in~(\ref{defb1}))
  \beqn \label{tempSS}
  \bar\theta_1 = \frac{ N^2}{8 \, b_1 ^2} \, 
         \left( \int_{\RR^N} M_{1,0,1}(v) \, |v|^3 \, dv \right)^{-2}.
  \eeqn
For instance in dimension $N=3$ we obtain 
$$
\bar\theta_1 = \frac{9 \pi}{64 b_1 ^2}.
$$
Moreover,  in the particular case of the hard-spheres cross-section (\ref{HScs}) 
in dimension $3$, we find $b_1 = b'_0 (4 \pi)/3$ and therefore 
$$
\bar\theta_1 = \frac{81}{1024  \pi \, (b'_0) ^2}.
$$

\subsection{Physical and mathematical motivation}

%\`A mettre ici~: justification par la physique de regarder les petites in\'elasticit\'e. 
%Cf. exemples dans articles de physique (Ernst-Brito ?) o\`u les in\'elasticit\'e typiques 
%sont de l'ordre de $0,99$. \c Ca peut \^etre reli\'e \'egalement aux sph\`eres dures 
%visco\'elastiques... 

% ici je fais un résumé de qqs arguments choppés dans les premiers chapitre du 
% bouquin de Brilliantov et Pöschel... 

For a detailed physical introduction to granular gases we refer to~\cite{BPlivre}. 
As can be seen from the references included in the latter, granular 
flows have become a subject of physical research on their own in the last decades, 
and for certain regimes of dilute and rapid flows this studies are based on kinetic theory. 
By contrast, the mathematical kinetic theory of granular gas is rather young and 
began in the late 1990 decade. 
We refer to~\cite{MMRI,MMII} for some (short) mathematical introduction to this theory 
and a (non exhaustive) list of references. As explained in these papers, granular gases are 
composed of grains of macroscopic size with contact collisional interactions, when one does 
not consider other additional possible self-interaction mechanisms such as gravitation 
-- for cosmic clouds for instance -- or electromagnetism -- for ``dusty plasmas" for instance --. Therefore the natural assumption about the 
binary interaction  between grains is that of inelastic hard spheres, with no loss of ``tangential relative 
velocity" (according to the impact direction) and a loss in  ``normal relative velocity" quantified in some 
(normal) restitution coefficient. The latter is either assumed to be constant as a first approximation 
(as in this paper) or can be more intricate: for instance it is a function of the modulus $|v'-v|$ of 
the normal relative velocity in the case of ``visco-elastic hard spheres" for instance 
(see~\cite{BPlivre}), which shall be studied in a forthcoming work~\cite{MM4}.
% CM : j'ai pas enlevé complètement la référence au travail en cours mais je l'ai allégé pour la 
% rendre moins précise... 

Simplified Boltzmann models like inelastic Maxwell molecules 
or pseudo inelastic hard spheres have been proposed (see~\cite{BCG00}) for which existence, uniqueness 
and global stability of a self-similar profile has been shown (see~\cite{BoCeTo:03,BCT2}), see also~\cite{BCT1} 
for similar results in the driven case of a thermal bath. However these models 
do not capture some crucial physical features of the cooling process of granular gas, like the tail behavior of 
the velocity distribution of the rate of decay of temperature (the so-called Haff's law). 
For (spatially homogeneous) inelastic hard spheres Boltzmann models, the existing mathematical works are: 
\begin{itemize}
\item the paper~\cite{BGP**} which shows {\em a priori}  polynomial 
and exponential moments bounds on any possible self-similar profile (resp. stationary solutions), 
whose existence is assumed, 
for freely cooling (resp. driven by a thermal bath) inelastic hard spheres with constant restitution coefficient;
\item the paper~\cite{GPV**} which shows existence of stationary solutions for 
inelastic hard spheres driven by a thermal bath, and improves the estimates on their tails of the previous paper into pointwise ones in this case; 
\item the paper~\cite{MMRI} which provides a Cauchy theory for freely cooling 
inelastic hard spheres with a broad family of collision kernels  
(including in particular restitution coefficients possibly depending on the relative 
velocity and/or the temperature), and studies the question of cooling in finite time or not 
for these various interactions; 
\item the paper~\cite{MMII}  which shows, for freely cooling inelastic hard spheres with constant restitution coefficient, 
existence of self-similar profile(s) as well as propagation of regularity and damping with time of singularity. 
\end{itemize}

In this paper we want to study the self-similarity properties of Boltzmann equation for 
inelastic hard spheres. Therefore as a natural 
first step we consider constant restitution coefficient $\alpha$ in order to 
have a self-similar scaling, which translates the study of self-similar solutions 
(often called homogeneous cooling states) to the study of stationary solutions for a rescaled 
equation. We also reduce to the case of restitution coefficients $\alpha$ close to $1$, 
that is, of small inelasticity. 
There are several physical as well as mathematical motivations for such a choice: 
\begin{itemize}
\item the first reason is related to the physical regime of the validity of kinetic theory: 
as explained in~\cite[Chapter~6]{BPlivre} for instance, 
the more inelasticity, the more correlations between grains are created during the binary collisions, and 
therefore the molecular chaos assumption, which is at the basis of the valdidity of Boltzmann's theory, 
suggests weak inelasticity to be the most effective; 
\item second as emphasized in~\cite{BPlivre} again, the case of restitution coefficient $\alpha$ close to $1$ 
has been widely considered in physics or mathematical physics since it allows to use expansions around 
the elastic case, and since conversely it is an interesting question to understand the connection of 
the inelastic case (dissipative at the microscopic level) to the elastic case (``hamiltonian" at the microscopic level); 
\item finally this case of a small inelasticity is reasonable from the viewpoint of applications, since it applies to interstellar dust clouds in astrophysics, or sands and dusts in earth-bound experiments, and more generally 
to visco-elastic hard spheres whose restitution coefficient is not constant but close to $1$ on the average. 
\end{itemize}

In this framework we shall show uniqueness and attractivity of self-similar solutions (in a suitable sense), and 
thus give a complete answer to the Ernst-Brito conjecture~\cite{EBJSP} (stated there for the simplified inelastic Maxwell model), for inelastic hard spheres with a 
small inelasticity. Moreover we give precise results about the elastic limit 
and deduce some quantitative informations about the weakly inelastic case. 

\subsection{Notation}\label{MMR:subsec:not}

%{\bf We shall denote by ``$C$'' various constants which do not depend on the
%collision kernel $b$ and by $K$ various constants which only depend  on the
%collision kernel $b$ (???)}
Throughout the paper we shall use the notation
$\langle \cdot \rangle = \sqrt{1+|\cdot|^2}$.
We denote, for any $p \in [1,+\infty]$, $q \in \R$ and weight function $\omega:\R^N \to \R_+$, the 
weighted Lebesgue space $L^p_q(\omega)$ by 
     \[ 
     L^p_q(\omega) := \left\{f: \RR^N \mapsto \RR \hbox{ measurable }; \; \;
     \| f \|_{L^p_q(\omega)}   < + \infty \right\},
     \]
with, for  $p < +\infty$,
   \[ \| f \|_{L^p _q (\omega)} = \left[ \int_{\R^N} |f (v)|^p \, \langle v
        \rangle^{pq} \, \omega(v) \, dv \right]^{1/p} \]
and, for $p = +\infty$, 
    \[ \| f \|_{L^\infty _q (\R^N)} = \sup_{v \in \R^N} |f (v)| \, 
         \langle v \rangle^{q} \omega(v). \]
We shall in particular use the exponential weight functions 
\beqn\label{defdem}
m  = m_{s,a}(v) := e^{ -a \, |v|^s} \quad\hbox{for}\quad a  \in (0,\infty), \,\, s \in (0,1),
\eeqn    
or a smooth version $m(v):= e^{- \zeta(|v|^2)}$ with  $\zeta \in C^\infty$ is a positive 
function such that $\zeta(r) =  r^{s/2}$ for any $r \ge 1$,  with $s \in (0,1)$.
%{\bf Enlever par la suite les expressions de $m := ...$ et pr\'eciser si $s < 1$ ou $s \le 1$. }

In the same way, the weighted Sobolev space $W^{k,p} _q (\omega)$ ($k \in \N$)
is defined by the norm
\[ 
\| f \|_{W^{k,p} _q (\omega)} =  \left[ \sum_{|s| \le k}
       \left\| \partial^s f (v) \right\|_{L^p _q (\omega)} ^p  \right]^{1/p}, 
\]
and as usual in the case $p=2$ we denote $H^k_q(\omega) =  W^{k,2} _q (\omega)$. The weight $\omega$ 
shall be omitted when it is $1$. 
Finally, for $g \in L^1_{2k}$, with $k \ge 0$, we introduce the following notation for the homogeneous moment of order $2k$
$$
{\bf m}_k(g) := \int_{\R^N} g \, |v|^{2 \, k} \, dv, 
$$
and we also denote by  $\rho(g) = {\bf m}_0(g)$ the mass of $g$, $\EE(g) = {\bf m}_1(g)$ 
the energy of $g$ and by $\theta(g) = \EE(g)/(\rho(g) \, N)$ the temperature associated 
to $g$ (when the distribution $g$ has $0$ mean). For any $\rho,\EE \in (0,\infty)$, $u \in \R^N$ we then introduce the subsets of $L^1$ of functions of given mass, mean velocity and energy
\bean
\Cc_{\rho,u}&:= &\{ h \in L^1_1; \, \int_{\R^N} h \, dv = \rho, \,\,  \int_{\R^N} h \, v \, dv = \rho \, u \}, \\
\Cc_{\rho,u,\EE} &:=& \{ h \in L^1_2; \, \int_{\R^N} h \, dv = \rho, \,\,  \int_{\R^N} h \, v \, dv = \rho \, u, \,\,  \int_{\R^N} h \, |v|^2 \, dv = \EE \}.
\eean
For any (smooth version of) exponential weight function $m$ we introduce the Banach space
$$ %\label{defLL}
\Ll^1 (m^{-1}) = L^1(m^{-1}) \cap \Cc_{0,0}.
$$

%%%%%%%%%%%%%%%%%%%%%%%%%%%%%%%%%%%%%%%%%%%%%%%%%

\subsection{Main results in self-similar variables}

\smallskip
Our main result, that we state now, deals with the evolution equation in self-similar variables
\beqn\label{eqresca}
  \frac{\partial g}{\partial t} = Q_\alpha(g,g) - \tau_\alpha \, \nabla_v \cdot (v g), \quad 
  g(0,.) = g_{\mbox{\scriptsize{in}}} \in \Cc_{\rho,0}
\eeqn
and the associated stationary equation, namely the self-similar profile equation
\beqn\label{profileeq3}
Q_\alpha(G,G) - \tau_\alpha \, \nabla_v \cdot (v \, G) = 0, \quad G \in \Cc_{\rho,0}.
\eeqn

\medskip
\begin{theo} \label{theo:uniq}
There is some constructive $\alpha_* \in (0,1)$  
such that for $\alpha \in [\alpha_*,1]$, and any given mass $\rho \in (0,\infty)$, we have:
\begin{itemize}
\item[(i)] For any $\tau>0$, the equation (\ref{eqresc}) admits a unique non-negative 
stationary solution with  mass $\rho$ and vanishing momentum. 
We denote by $\bar G_\alpha$ the self-similar profile obtained 
by fixing $\tau = \tau_{\alpha}$ (defined by (\ref{taualpha})). 
\item[(ii)] Let define $\bar G_1 = M_{\rho,0,\bar\theta_1}$ the Maxwellian  distribution with mass $\rho$, 
momentum $0$ and ``{\em quasi-elastic self-similar temperature}" $\bar\theta_1$ defined in~(\ref{tempSS}).
The path of self-similar profiles $\alpha \to \bar G_\alpha$  parametrized by the normal restitution coefficient 
is $C^1$ from $[\alpha_*,1]$ into $W^{k,1} \cap L^1(e^{a \, |v|})$ for any $k \in \N$ and some $a \in (0,\infty)$.
\item[(iii)] For any $\alpha \in [\alpha_*,1]$, the linearized collision operator 
\beqn\label{defLa}
h \mapsto \LL_\alpha \, h := 2 \, \tilde Q_\alpha (\bar G_\alpha, h) - \tau_\alpha \, \nabla_v \cdot (v \, h)
\eeqn
is well-defined and closed on $\Ll^1(m^{-1})$ for any exponential weight function $m$ with exponent $s\in(0,1)$ 
(defined in (\ref{defdem})).  
Its spectrum decomposes between a part 
which lies in the half-plane $\{ \mbox{{\em Re}} \, \xi \le \bar \mu \}$ for some constructive 
$\bar \mu < 0$, and some remaining discrete eigenvalue $\mu _\alpha$. This eigenvalue is real negative and satisfies 
\beqn\label{expanmua}
\mu _\alpha = - \rho \, (1-\alpha) + \OO(1-\alpha)^2 \quad\hbox{when}\quad \alpha \to 1.
\eeqn
The associated eigenspace is of dimension $1$ and then denoting by $\phi_\alpha = \phi_\alpha(v)$ the unique associated 
eigenfunction such that $\| \phi_\alpha \|_{L^1_2} = 1$ and $\phi_\alpha(0) < 0$, there holds 
$\phi_\alpha \in \Ss(\R^N)$ (with bounds of regularity independent of $\alpha$) and 
\beqn\label{phiaTOphi1}
\phi_\alpha \rightarrow \phi_1 := c_0 \, \big( |v|^2 - N \, \bar\theta_1 \big) \, \bar G_1
\quad\hbox{ as }\quad  \alpha \to 1, 
\eeqn
where $c_0$ is the positive constant such that $\| \phi_1 \|_{L^1_2} = 1$. 
Finally one has constructive decay estimates on the semigroup 
associated to this spectral decomposition in this Banach space (see the 
key Theorem~\ref{thLalpha}  and the following point).
%$$
%c_0 = \left( \int_{\R^N} \big( |v|^2 - N \, \bar \theta \big)^2 \, M_{\rho, 0, \bar \theta} \, dv \right)^{-1/2}.
%$$
%
%\item[(iv)] For any $\alpha \in [\alpha_*,1]$ and any exponential weight function $m$ of exponent $s\in(0,1)$, 
%the self-similar profile $\bar G_\alpha$ is locally attractive in $L^1(m^{-1}) \cap H^{N+4}$,  with exponential rate. 
%More precisely, there are $\eps_*>0$ and $C_*>0$ such that if the initial datum satisfies
%$$
%\| g_{in} \|_{ L^1(m^{-1}) \cap H^{N+4} } \le \eps_*, 
%$$ 
%the solution $g$ to~(\ref{eqresca}) satisfies 
%\beqn\label{cvgcess}
%g(t,.) = \bar G_\alpha + c_\alpha (t) \, \phi_\alpha + r_\alpha(t,.), 
%\eeqn
%with $c_\alpha (t) \in \R$ and $r_\alpha(t,.) \in L^1_2(\R^N)$ such that 
%\beqn\label{cvgcess2}
%|c_\alpha (t)| \le  C_* \, e^{\mu _\alpha \, t}, \qquad
%\| r_\alpha(t,.) \|_{L^1_2} \le C_* \, e^{\bar\mu \, t}.
%\eeqn
%
%\item[(v)] Moreover, the self-similar profile $\bar G_\alpha$ is indeed globally attractive on bounded subset of 
%$L^1_3 \cap W^{k_*,1}$ under some explicit smallness condition on $(1-\alpha_*)$ in the following sense. For  
%any initial datum $0 \le g_{\mbox{\scriptsize{{\em in}}}}  \in L^1_3 \cap W^{k_*,1}$  there exists $\alpha_{**} \in 
%(\alpha_*,1)$ depending explicitly on the mass, energy, third moment  and $W^{k_*,1}$ norm of $g_{\mbox
%{\scriptsize{{\em in}}}}$,  such that the solution $g(t,.)$ to~(\ref{eqresc}) satisfies 
%(\ref{cvgcess}) for any $\alpha \in (\alpha_{**},1)$.
%
\item[(iv)] The self-similar profile $\bar G_\alpha$ is globally attractive on bounded 
subsets of $L^1_3$ under some smallness condition on the inelasticity in the following sense. 
For any $\rho,\EE_0, M_0 \in (0,\infty)$ there exists $\alpha_{**} \in (\alpha_*,1)$,  
$C_* \in (0,\infty)$ and $\eta \in (0,1)$, such that for any initial datum  satisfying 
$$ 
0 \le g_{\mbox{\scriptsize{{\em in}}}}  \in L^1_3 \cap \Cc_{\rho,0,\EE_0}, 
\qquad \| g_{\mbox{\scriptsize{{\em in}}}} \|_{L^1_3} \le M_0,
$$
the solution $g$ to~(\ref{eqresca}) satisfies 
\beqn \label{cvgirreg} 
\| g_t - \bar G_\alpha \|_{L^1 _2} \le e^{(1-\eta) \, \mu _\alpha \, t}.
\eeqn  

\item[(v)] Moreover, under smoothness condition on the initial datum one may 
prove a more precise asymptotic decomposition, and construct Liapunov functional 
for the equation (\ref{eqresca}). 
More precisely, there exists $k_* \in \N$ and, for any exponential weight $m$ as defined in (\ref{defdem}) 
and any $\rho,\EE_0, M_0 \in (0,\infty)$, there exists $\alpha_{**} \in (\alpha_*,1)$ and 
a constructive functional $\HH : H^{k_*} \cap L^1(m^{-1}) \to \RR$ such that, first, 
for any initial datum 
$0 \le g_{\mbox{\scriptsize{{\em in}}}}  \in H^{k_*} \cap L^1(m^{-1}) \cap \Cc_{\rho,0,\EE_0}$ satisfying 
$$ 
 \| g_{\mbox{\scriptsize{{\em in}}}} \|_{H^{k_*} \cap L^1(m^{-1})} \le M_0,
$$
the solution $g$ to~(\ref{eqresca}) satisfies 
\beqn\label{cvgcess}
g(t, \cdot) = \bar G_\alpha + c_\alpha (t) \, \phi_\alpha + r_\alpha(t, \cdot), 
\eeqn
with $c_\alpha (t) \in \R$ and $r_\alpha(t, \cdot) \in L^1_2(\R^N)$ such that 
\beqn\label{cvgcess2}
|c_\alpha (t)| \le  C_* \, e^{\mu _\alpha \, t}, \qquad
\| r_\alpha(t, \cdot ) \|_{L^1_2} \le C_* \, e^{(3/2) \, \mu_\alpha \, t}.
\eeqn 
And second when the initial datum satisfies additionally 
$$
g_{\mbox{\scriptsize{{\em in}}}} \ge M_0^{-1} \, e^{-M_0 \, |v|^8},
$$ 
the solution satisfies also
$$ %\label{liapdec}
t \mapsto \HH(g(t, \cdot)) \quad\hbox{is strictly decreasing }
$$
(up to reach the stationary state $\bar G_\alpha)$.
%And second, one can construct a Liapunov functional for the equation (\ref{eqresca}). 
%More precisely, there exists $k_* \in \N$ and, for any exponential weight $m$ and any 
%$\rho,\EE_0, M_0 \in (0,\infty)$, there exists $\alpha_{**} \in (\alpha_*,1)$ and 
%a constructive functional $\HH : H^{k_*} \cap L^1(m^{-1}) \to \RR$ such that for any initial datum $0 \le g_{\mbox{\scriptsize{{\em in}}}}  \in H^{k_*} \cap L^1(m^{-1}) \cap \Cc_{\rho,0,\EE_0}$ satisfying 
%$$ 
% \| g_{\mbox{\scriptsize{{\em in}}}} \|_{H^{k_*} \cap L^1(m^{-1})} \le M_0,
% \qquad g_{\mbox{\scriptsize{{\em in}}}} \ge M_0^{-1} \, e^{-M_0 \, |v|^8},
%$$
%the solution $g$ to~(\ref{eqresca}) satisfies 
%$$ %\label{liapdec}
%t \mapsto \HH(g(t, \cdot)) \quad\hbox{is strictly decreasing }
%$$
%(up to reach the stationary state $\bar G_\alpha)$.
\end{itemize}
\end{theo}

%{\bf LE DIRE UNE FOIS POUR TOUTE ET PUIS NE PLUS LE REPETER A CHAQUE LEMME!} 
\begin{rems} 

1) All the constants appearing in this theorem are contructive,  which means that they can be made explicit, 
and in particular that the proof does not use any compactness argument. Unless otherwise mentioned, 
these constants will depend on $b$, on the dimension $N$, and on some bounds on the initial datum 
but never on the inelasticity parameter $\alpha \in (0,1]$. 
\smallskip

2) Theorem~\ref{theo:uniq} establishes that conjectures 1 and 2 in \cite[Section~5]{MMII} holds true at least for 
weak inelastic model (that means for $\alpha$ close enough to $1$). 
\smallskip

3) In point (iv), the condition on the restitution coefficient depends on the 
mass, temperature and $L^1_3$ norm of the initial distribution, but this dependence is {\em not}
a perturbative condition of closeness to the self-similar profile. This fact relies on the 
so-called ``entropy-entropy production" estimates which yields ``overlinear" 
Gronwall-type estimates, and the decoupling of the timescales of energy dissipation 
and entropy production. 
\smallskip

4) In (\ref{cvgcess2}) one can prove $\| r_\alpha(t, \cdot) \|_{L^1_2} \le C_\zeta \, e^{\zeta\, \mu_\alpha \, t}$ 
for any $\zeta \in (1,2)$. Remark that here we do not have the decay rate $e^{\bar \lambda t}$ on the remaining 
part when one ``removes" from $g_t - \bar G_\alpha$ the projection on the energy eigenvalue, where 
$\bar \lambda<0$ would be some constant independent of $\alpha$ related to the second non-zero 
eigenvalue of $\LL_\alpha$.  
This is due to the coupling effect of the bilinear term, which mixes the different part of the spectral 
decomposition.
\smallskip

5) As a subproduct the above result provides an alternative argument to the one 
of~\cite[Section~3]{MMII} to show uniform (in time and inelasticity parameter) 
non-concentration bounds on the rescaled equation, in the case of $\alpha$ close to $1$ and 
a general initial datum $g_{in} \in L^1_3$ (whereas the proof of~\cite[Section~3]{MMII} 
was valid for all $\alpha \in (0,1)$ but for some initial datum $g_{in} \in L^1_3 \cap L^p$, $p \in (1,\infty]$). 
\smallskip

6) Our results show that no bifurcation occurs for the self-similar profile for $\alpha$ close to $1$. 
We do not know at now if some bifurcations occur for other values of the inelasticity parameter. 
Therefore we do not know if there is a continuous 
branch of self-similar profiles parametrized by $\alpha \in [0,1]$ (even if we know from~\cite{MMII} 
that self-similar profiles exist for all values of the inelasticity paramaters). 
The best one could say in terms of ``connectivity" from the estimates we have proved on the profile 
together with the classical theory of topological degree (see~\cite{Nirenb} for instance) 
is that there is a set $K \subset [0,1] \times \cal{F}$ (where 
$\cal{F}$ is for instance the set of positive functions in the Schwartz space with given mass) 
which is compact, connected, and such that for any $\alpha \in [0,1]$, the intersection 
$K \cap \{ \alpha \} \times \cal{F}$ is not empty. 
\end{rems}

%\begin{rem} More precisely $g_t - \bar G_\alpha$ splits between a part in $\Pi^e _\alpha(L^1(m^{-1}))$ 
%decreasing asymptotically like $C \, e^{-[ \rho - \OO (1-\alpha)] \, (1-\alpha)  \, t}$ 
%and a remaining part decreasing asymptotically like $C \, e^{- (1 + \eta) \, [ \rho - \OO (1-\alpha) ] \, (1-\alpha)  \, t}$ 
%with $\eta > 0$ (all the constants are explicit).

%The coefficient $(5/4)$ in the above exponential rate is due to the power of $h^1$ 
%in the integral equation for $h^2$. With more care, we can obtain a coefficient $(1+\eta)$ for any $\eta \in (1,2)$ because of the quadratic structure of the collision term. The coefficient $2$ seems to be exclude because of the use of interpolation arguments. 
%\end{rem}

%%%%%%%%%%%%%%%%%%%%%%%%%%%%%%%%%%%%%%%%%%%%%%%%%

\subsection{Coming back to the original equation}

%%%%%%%%%%%%%%%%%%%%%%%%%%%%%%%%%%%%%%%%%%%%%%%%%

When coming back to the original equation~(\ref{eqBol1}) with the help 
of (\ref{GtoF}) and (\ref{gtof}), Theorem~\ref{theo:uniq} translates into the 

\begin{theo} \label{theo:originalvar} 
There is a constructive $\alpha_* \in (0,1)$
such that for $\alpha \in [\alpha_*,1]$, and any given mass $\rho \in (0,\infty)$, we have    
\begin{itemize}
\item[(i)] Up to a translation of time there exists a unique self-similar solution $\bar F_\alpha$ 
of the equation~(\ref{eqBol1}) with mass $\rho$, and it is given by 
$$
\bar F_\alpha (t,v) = \big(1+ \tau_\alpha \,  t \big)^N \, \bar G_\alpha \big((1+ \tau_\alpha \, t) \, v\big), \qquad
\tau_\alpha = \rho \, (1-\alpha), 
$$
where $\bar G_\alpha$ was obtained in Theorem~\ref{theo:uniq}. 
More precisely, if $F_\alpha$ is a  solution of~(\ref{eqBol1}) of 
the form~(\ref{shapeselfsim}) and of mass $\rho$, there exists $t_0 \in \R$ such that 
$F_\alpha(t,v) = \bar F_\alpha(t+t_0,v)$ for 
any $t \ge \max \{ 0,-t_0\}$ and any $v \in \R^N$. 
%
%\item[(ii)] For any fixed time $t \ge 0$, the path of self-similar 
%solutions $\alpha \to \bar F_\alpha$ parametrized by the normal restitution coefficient is well-defined 
%on $[\alpha_*,1]$, and it is $C^1$ from $[\alpha_*,1]$ into $H^k(e^{A|v|})$ for 
%any $k \in \N$ and some $A = A_k>0$.  {\bf idem!} 
%
%\item[(iii)] When $\alpha$ goes to $1$ (elastic limit), the 
%self-similar solutions $\bar F_\alpha(t,v)$ converges at any fixed time $t$ to 
%$\bar G_1 = M_{\rho,0,\bar\theta}$ where $M_{\rho,0,\bar\theta}$ is the Maxwellian  distribution with mass $\rho$, 
%momentum $0$ and a ``{\em quasi-elastic self-similar temperature}" $\bar \theta$ defined in~(\ref{tempSS}).
%
\item[(ii)]  The self-similar solution $\bar F_\alpha$ is globally attractive on bounded subsets of $L^1_3$ 
under some smallness condition on $(1-\alpha_{*})$ in the following sense. For any $\rho,\EE_0, M_0 \in (0,\infty)$ 
there exists $\alpha_{**} \in (\alpha_*,1)$ and $\eta \in (0,1)$ such that for any $q \in \N$ there is $c_q \in (0,\infty)$ 
such that for any initial datum satisfying 
$$ 
0 \le f_{\mbox{\scriptsize{{\em in}}}}  \in L^1_3 \cap \Cc_{\rho,0,\EE_0}, 
\qquad \| f_{\mbox{\scriptsize{{\em in}}}} \|_{L^1_3} \le M_0,
$$
the solution $f(t, \cdot)$ to~(\ref{eqBol1}) satisfies
$$ %\label{cvgfirreg} 
%\mbox{ } \quad 
\| f(t, \cdot) - \bar F_\alpha(t, \cdot) \|_{L^1(|v|^q)} \le c_q \, \, (1 + \tau_\alpha \, t)^{(1-\eta) \, \mu_\alpha/\tau_\alpha - q} =
c_q \, (1 + \tau_\alpha \, t)^{- (1-\eta) - q + \OO(1-\alpha)}.
$$ 
\item[(iii)] 
Moreover, there exists $k_* \in \N$ and, for any exponential weight $m$ as defined in (\ref{defdem}) 
and any $\rho,\EE_0, M_0 \in (0,\infty)$, there exists $\alpha_{**} \in (\alpha_*,1)$ such that, 
for any initial datum 
$0 \le f_{\mbox{\scriptsize{{\em in}}}}  \in H^{k_*} \cap L^1(m^{-1}) \cap \Cc_{\rho,0,\EE_0}$ satisfying 
$$ 
 \| f_{\mbox{\scriptsize{{\em in}}}} \|_{H^{k_*} \cap L^1(m^{-1})} \le M_0
$$
the solution $f$ to~(\ref{eqresca}) satisfies 
\beqn \label{cvgceorig}
f(t, \cdot) = \bar F_\alpha(t,\cdot) + \tilde c_\alpha(t) \, \psi_\alpha(t,\cdot) + \tilde r_\alpha(t,\cdot)
\eeqn
where
$$
\psi_\alpha (t,v) = (1 + \tau_\alpha \, t )^N \, \phi_\alpha\Big((1+\tau_\alpha \, t) \, v\Big), \quad 
\tilde c_\alpha (t) = c_\alpha \left( \frac{\ln (1+\tau_\alpha \, t)}{\tau_\alpha}\right).
$$
In this expansion, the different terms have the following asymptotic behaviors (for any given $q \ge 0$): 
$$
\| \bar F_\alpha(t, \cdot ) \|_{L^1(|v|^q)} = (1 + \tau_\alpha \, t)^{- q} \| \bar G_\alpha \|_{L^1(|v|^q)},
$$
$$
\| \psi_\alpha (t, \cdot)  \|_{L^1(|v|^q)} = (1 + \tau_\alpha \, t)^{- q} \| \bar G_\alpha \|_{L^1(|v|^q)},
$$
$$
\left| \tilde c_\alpha(t) \right| \le C_* \, (1+ \tau_\alpha \, t)^{\mu_\alpha / \tau_\alpha} = 
C_* \,  (1+ \tau_\alpha \, t)^{-1 + \OO(1-\alpha)},
$$
$$
\exists \, C_q >0 \ ; \ \ 
\left\| \tilde r_\alpha \right\|_{L^1(|v|^q)} \le C_q \, (1 + \tau_\alpha \, t)^{(3/2) \, \mu_\alpha/\tau_\alpha - q} =
C_q \, (1 + \tau_\alpha \, t)^{- (3/2) - q + \OO(1-\alpha)}.
$$
Hence the leading term in the expansion (\ref{cvgceorig}) is, as expected, the self-similar solution, and 
the first order correction behond self-similarity is given by the second term, that is the projection  
onto the eigenspace of the ``energy eigenvalue".

\item[(iv)] We may make more precise Haff's law on 
the asymptotic behavior of the granular temperature (see~\cite{MMII}) in the following way. 
Under the assumptions of point (iii),  the solution $f=f(t,v)$ to~(\ref{eqBol1}) satisfies 
\beqn\label{preciseHaffLaw}
\EE(f(t.,)) 
= \frac{\EE(\bar G_\alpha)}{ \big(1+ \tau_\alpha \, t \big)^2} 
  +  \OO\Bigg( \frac{1}{\big(1+ \tau_\alpha \,  t \big)^{3 + \OO(1-\alpha)}} \Bigg).
\eeqn
\item[(v)] Under the assumptions of point (iii) the rescaling 
by the square root of the energy familiar to physicists is rigorously justified in the sense: 
the solution $f=f(t,v)$ to~(\ref{eqBol1}) satisfies for $t \to +\infty$
$$
\EE(f_t)^{N/2} \, f\Big(t,\EE(f_t)^{1/2} \, v\Big) 
\to \EE(\bar G_\alpha)^{N/2} \bar G_\alpha \Big(\EE(\bar G_\alpha)^{1/2} \, v \Big)  
 \quad \hbox{in}\quad L^1.
$$
%\item[(v)] If we denote by $\psi_\alpha= \psi_\alpha(t,v)$ the self-similar solution associated to 
%$\phi_\alpha=\phi_\alpha(t,v)$ by ~(\ref{gtof}) then the solution $f_t$ to~(\ref{eqBol1}) 
%satisfies the asymptotic expansion (under the assumptions of points~(iv) or~(v)) 
%$$
%f = \bar F_\alpha + c' _\alpha (t) \, \psi_\alpha  + R_\alpha
%$$
%with 
%$$
%c' _\alpha (t) \le C \, \big(1+\rho (1-\alpha)t \big)^{-\bar \lambda + \OO(1-\alpha)}, \quad 
%\| R_\alpha \|_{L^1} \le C \,\big(1+\rho (1-\alpha)t \big)^{-(1+\eta) \, \bar \lambda +\OO(1-\alpha)}. 
%$$
%Therefore $\psi_\alpha$ describes the profile of the solution beyond self-similarity. Note 
%that from the previous theorem,  for the elastic limit $\alpha \to 1$, the 
%self-similar functions $\psi _\alpha(t,v)$ converges at any fixed time $t$ to 
%$\phi_1$ which is known explicitely (see Theorem~\ref{theo:uniq}). 
 \end{itemize}
\end{theo}

\begin{rem}
We see from this theorem that the convergence towards the self-similar solution 
in indeed faster than the convergence towards the Dirac mass (hence justifying its 
interest), but also that the speed of convergence towards this 
self-similar solution degenerates to $0$ as $\alpha \to 1$ (because $\tau_\alpha \to 0$
when $\alpha\to1$). This fact is surprising, 
since the self-similar solution converges towards a stationary Maxwellian 
distribution in the elastic limit, and the latter is known to be exponentially attractive for the elastic equation 
(see~\cite{GM:04} for instance). As we shall see this is related to the fact 
that a birfurcation occurs in the spectrum of the linearized collision 
operator at $\alpha =1$ (namely the eigenvalue corresponding the kinetic 
energy vanishes at $\alpha =1$ whereas it is non-zero for $\alpha \in [\alpha_*,1)$). 
This remark may explain the fact that in the quasi-elastic limit considered -- in dimension $1$ -- 
in~\cite{CaVi}, it is proved that the rate of relaxation towards the self-similar 
solution is worse than any polynomial. 
\end{rem}

\medskip\noindent
{\sl Proof of Theorem~\ref{theo:originalvar}}.  
Except for points (i) and (v) this theorem is an obvious translation of Theorem~\ref{theo:uniq}. 
In order to prove (i), one first remarks that for two given self-similar solutions $F$ and $\tilde F$, 
there holds
$$
F(t,v) = (V_0+A \, t)^N \, G_A \big((V_0+A \, t) \, v\big), \quad
\tilde F(t,v) = (\tilde V_0+\tilde A \, t)^N \, G_{\tilde A} \big((\tilde V_0+\tilde A \, t) \, v\big), 
$$
and thus from Theorem~\ref{theo:uniq} 
$$
G_{\tilde A} (v) = \left( {A \over \tilde A} \right)^N \, G_A \left( {A \over \tilde A} \, v \right).
$$
We deduce 
$$
\tilde F(t,v) = \left(\tilde V_0 \, {A \over \tilde A} + A \, t\right)^N \, 
G_A \left( \left(\tilde V_0 \, {A \over \tilde A} + A \, t\right) \, v\right) = F(t+t_0,v)
$$
with 
$$
t_0 = {V_0 \over A} \, \left( {\tilde V_0 \over \tilde A} - 1\right).
$$
In order to prove (v), we introduce the function $\xi(t) = \EE(\bar G_\alpha)^{1/2} / [ \EE(f_t)^{1/2} \, (1+\tau_\alpha \, t)]$ and we compute
\bean
&&\left\| \EE(f_t)^{N/2} \, f\big(t,\EE(f_t)^{1/2} \, \cdot \big) -\EE(\bar G_\alpha)^{N/2} \bar G_\alpha(\EE(\bar G_\alpha)^{1/2} \,  \cdot )  \right\|_{L^1} = \vphantom{\int}\\
&&\,\,=
  \| g(\tau_\alpha^{-1} \ln (1+\tau_\alpha t), \cdot ) - \xi(t)^N \, \bar G_\alpha(\xi(t) \, \cdot) \|_{L^1} \vphantom{\int}\\
&&\,\,\le
   \| g(\tau_\alpha^{-1} \ln (1+\tau_\alpha t), \cdot ) -  \bar G_\alpha  \|_{L^1}
   + |\xi(t)^N-1| \, \| \bar G_\alpha \|_{L^1}   + \xi(t)^N \| \bar G_\alpha(\xi(t) \, \cdot) -  \bar G_\alpha \|_{L^1}. \qquad\vphantom{\int}
\eean
Using now (\ref{expanmua}),  (\ref{cvgcess}), (\ref{cvgcess2}), (\ref{preciseHaffLaw}) and the fact that $\bar G_\alpha$ is bounded in $W^{1,1}_1$ uniformly in $\alpha \in (\alpha_*,1)$ from Theorem 1.1 (ii), we deduce 
$$
\left\| \EE(f_t)^{N/2} \, f\big(t,\EE(f_t)^{1/2} \, \cdot \big) -\EE(\bar G_\alpha)^{N/2} \bar G_\alpha(\EE(\bar G_\alpha)^{1/2} \,  \cdot )  \right\|_{L^1} \le C\, (1+\tau_\alpha \, t)^{-1 + \OO(1-\alpha)},
$$
for some constant $C \in (0,\infty)$ (which depends in particular on the upper bound on 
$\| \bar G_\alpha \|_{W^{1,1}_1}$), from which (v) follows.
 \qed 

\begin{rem} 
Let us emphasize that the temperature $\bar\theta_1$ of the limit Maxwellian $\bar G_1$ is 
``universal" in the sense that it depends only on the collisional cross-section $b$ 
(through its angular momentum), and not for instance on the density distribution. 

The temperature of the self-similar solution $\bar F_\alpha=F_\alpha(t,v)$ associated to 
a self-similar profile $\bar G_\alpha$ decreases like 
  \[ \theta(\bar F_\alpha(t,\cdot)) = \frac{\theta(\bar G_\alpha)}{(1+\rho(1-\alpha)t)^2}. \]
Hence when $\alpha$ is close to $1$ (small inelasticity) we obtain 
  \[ \theta(\bar F_\alpha(t,\cdot)) \approx \frac{\bar\theta_1}{(1+\rho(1-\alpha)t)^2}. \] 
Therefore, as soon as the self-similar solutions correctly describe the 
asymptotic (at least in the framework of point (ii) of Theorem~\ref{theo:originalvar}), which 
is conjectured by physicists, generic solutions satisfy 
  \[ \theta(f_\alpha(t,\cdot)) \sim_{t \to \infty} 
      \left( \frac{\bar\theta_1}{\rho^2 (1-\alpha)^2} \right) \, t^{-2} \]
for an inelasticity coefficient $\alpha$ close to $1$. 

Hence we shall denote the universal quantity $\bar\theta_1$ as a 
``quasi-elastic self-similar temperature''. Remark that its definition 
as the temperature of $\bar G_1$ seems to depend on the 
choice of the scaling. However changing this scaling by some asymptotically 
equivalent one, as $\alpha \to 1$, would only adds a factor which would 
then disappear when coming back to the 
solution to the original equation (\ref{eqBol1}). Therefore a more ``canonical"  
way to define this quasi-elastic self-similar temperature could be 
  \[ \bar\theta_1 =  \rho^2 \,\lim_{\alpha \to 1}  \left( (1-\alpha^2) \, 
          \lim_{t \to +\infty} \theta\big(f_\alpha(t,\cdot)\big) \, t^2 \right) \]
where $f_\alpha$ denotes a generic solution with mass $\rho$ to equation (\ref{eqBol1}). 
\end{rem}

%%%%%%%%%%%%%%%%%%%%%%%%%%%%%%%%%%%%%%%%%%%%%%%%%

\subsection{Method of proof and plan of the paper}

%%%%%%%%%%%%%%%%%%%%%%%%%%%%%%%%%%%%%%%%%%%%%%%%%

The first main idea of our method is to consider the rescaled equations (\ref{eqresc}) and (\ref{eqprofGtau}) with an inelasticity dependent anti-drift coefficient $\tau_\alpha$ which exactly ``compensates" the loss of elasticity of the collision operator (in the sense that it compensates its loss of kinetic energy). This scaling allows by some technical estimates to prove uniform bounds according to $\alpha$ for the family of self-similar profiles $G_\alpha$ to the equation (\ref{profileeq3}). 
The second main idea consists in decoupling the variations along the ``energy direction" and its ``orthogonal direction". 
This decoupling makes possible to identify the limit of different objects as $\alpha \to 1$ (among them the limit 
of $G_\alpha$). 
The third main idea is to use systematically the knowledges on the elastic limit problem, once it has been 
identified thanks to the previous arguments. In particular we use the spectral study of the linearized problem and the dissipation entropy-entropy inequality for the elastic problem. This allows to argue by perturbative method. 
Let us emphasize that this perturbation is singular in the classical sense because of the addition of a (limit vanishing) 
first-order derivative operator, but also because of the gain of one more conservative quantity at the limit 
(which implies in particular at the linearized level that the ``energy eigenvalue" $\mu_\alpha$ is negative 
for $\alpha \not = 1$ but converges to $\mu_1 = 0$ in the limit $\alpha \to 0$). 

\smallskip
In Section 2, we use the regularity properties of the collision operator in order to establish on the one hand 
that the family $(G_\alpha)$ is bounded in $H^\infty \cap L^1(m^{-1})$ uniformly according to the 
inelastic parameter $\alpha$ (the key argument being the use of the entropy functional which 
provides uniform lower bound on the energy of $G_\alpha$) and on the other hand that the difference 
of two self-similar profiles in any strong norm may be bounded by the difference of these ones in weak norm 
(the key idea is a bootstrap argument). This last point shall allow to deal with the loss of derivatives and weights 
in the operator norms used in the sequel of the paper. 

\smallskip
In Section 3, we prove that $\alpha \mapsto Q^+_\alpha$ is H\"older continuous in the norm of its graph and 
is H\"older differentiable in a weaker norm. As a consequence we deduce that $G_\alpha \to \bar G_1$ 
when $\alpha \to 1$ with explicit ``H\"older" rate, which (partially) proves point (ii) Theorem~\ref{theo:uniq}. 
The cornerstone of the proof is the decoupling of the variation $G_\alpha - \bar G_1$ between the 
``energy direction" and its ``orthogonal direction".

\smallskip
In Section 4, we prove uniqueness of the profile $\bar G_\alpha$ for small inelasticity (point (i) of 
Theorem~\ref{theo:uniq}) by a variation around the implicit function theorem. We also deduce 
that $\alpha \mapsto \bar G_\alpha$ is differentiable at $\alpha=1$.

\smallskip
Section 5 is devoted to the study of the linearized operator $\LL_\alpha$, and we partially inspire from 
the method of \cite{GM:04}. We prove point (iii) of Theorem~\ref{theo:uniq} and we end the proof of point (ii) of 
Theorem~\ref{theo:uniq}. We obtain information on the localization of the spectrum and we 
establish some decay estimates on the associated semigroup. Let us emphasize that for technical 
reasons we state our results in an $L^1$ framework (because mainly we are not able to generalize 
Lemma~\ref{lem:cvgLalpha} to an $L^2$ framework), which makes the spectral analysis more intricate.  
The proof proceeds as follows (the cornerstone idea is again the decoupling of the variations in the 
``energy direction" and its ``orthogonal direction"). First, we localize the essential spectrum in the 
half plan $\Delta_{\bar\mu}^c = \{ z \in \CC, \,\, \Re e \, z \le \bar \mu < 0 \}$ with the help of Weyl's theorem, 
the compactness properties of $\LL_\alpha$ and the ``rough" (H\"older type) convergence of 
$Q^+_\alpha(\bar G_\alpha, \cdot )$ to $Q^+_1(\bar G_1, \cdot )$ in the ``good" norm of the graph. 
Second, we localize the discrete spectrum lying in $\Delta_{\bar\mu} = \{ z \in \CC, \,\, \Re e \, z \ge \bar \mu \}$ 
in the disc $\{ z \in \CC, \,\, |z| \le C \, (1-\alpha) \}$, thanks to estimates on the resolvent of $\LL_\alpha$. 
Third we establish that the spectrum $\Sigma(\LL_\alpha)$ of $\LL_\alpha$ satisfies 
$\Sigma(\LL_\alpha) \cap \Delta_{\bar\mu} = \{ \mu_\alpha \}$, where $\mu_\alpha$ has multiplicity $1$ 
(the proof mainly takes advantage of the ``precise" convergence of $Q^+_\alpha(\bar G_\alpha, \cdot)$ 
to $Q^+_1(\bar G_1, \cdot )$ in ``bad" norm, together with a regularity estimate holding on the discrete eigenspace).
Last we establish the expansion (\ref{expanmua}) using the energy equation associated to the eigenvalue $\mu_\alpha$. The decay properties of the linear semigroup are then deduced from resolvent estimates and 
the above localization of the spectrum.

\smallskip
Section 6 is devoted to the proof of points (iv) and (v) in Theorem~\ref{theo:uniq} which is split in several steps. 
First we establish a ``linearized asymptotic stability result" by decoupling the evolution equation (\ref{eqresca}) 
along the ``energy direction" and its ``orthogonal direction", and using the semigroup decay estimates and the 
quadratic structure of the collision operator. Second we establish a ``non-linear stability result" by decoupling the 
evolution equation (\ref{eqresca}), using the energy dissipation equation along the  ``energy direction"  and the 
entropy production method on its ``orthogonal direction" (let us mention that this method follows closely the physical 
idea that for small inelasticity the ``molecular" timescale of thermalization of velocity distribution 
decouples from the ``cooling" timescale of dissipation of energy). Third we prove the asymptotic decomposition 
and we exhibit a Liapunov functional for smooth initial data (point (v)) by gathering 
(and slightly modifying) the two preceding steps. Fourth and last, we prove point (iv) for general initial data, 
gathering the previous arguments with the 
decomposition of solutions between a smooth part and a small remaining part as introduced in \cite{MV**}. 

%%%%%%%%%%%%%%%%%%%%%%%%%%%%%%%%%%%%%%%%%%%%%%%

\section{{\it A posteriori} estimates on the self-similar profiles} 
\setcounter{equation}{0}
\setcounter{theo}{0}

%%%%%%%%%%%%%%%%%%%%%%%%%%%%%%%%%%%%%%%%%%%%%%%

In this section we prove various {\it a posteriori} 
regularity and decay estimates on the self-similar profiles 
(or the differences of self-similar profiles), uniform as $\alpha \to 1$, 
which shall be useful in the sequel. 

\subsection{Uniform estimates on the self-similar profiles}

%%%%%%%%%%%%%%%%%%%%%%%%%%%%%%%%%%%%%%%%%%%%%%%

For any $\alpha \in (0,1)$ we consider $\GG_\alpha$ the set of all the self-similar profiles 
of the inelastic Boltzmann equation (\ref{eqBol1}) with inelasticity coefficient $\alpha$, 
with given mass $\rho \in (0,+\infty)$ and finite energy. More precisely, we define $\GG_\alpha$
as the following set of functions 
$$
\GG_\alpha := \Big\{ 0 \le G \in L^1_2 \quad \hbox{satisfying} \quad (\ref{profileeq3}) \Big\}.
$$
For some fixed $\alpha_0\in (0,1)$, we also define 
$$
\GG  = \cup_{  \alpha \in [\alpha_0,1) }  \GG_\alpha.
$$ 
% fixed once for all (and which may be arbitrary small).  

The fact that for any $\alpha \in (0,1)$, $\GG_\alpha$ is not empty 
was proved in~\cite{MMII}, where a solution of~(\ref{profileeq3}) was 
built within the class of radially symmetric functions belonging to the Schwartz space. 
Here we show that any self-similar profile $G_\alpha \in \GG$ belongs to the Schwartz space and that 
decay estimates, pointwise lower bound and regularity estimates can be made uniform according to the inelasticity coefficient $\alpha \in [\alpha_0,1)$. Let us emphasize once again that the choice of the 
velocity rescaling parameter $\tau_\alpha = \rho \, (1-\alpha)$ in~(\ref{profileeq3}) is fundamental in order to get that 
uniformity in the limit $\alpha \to 1$. Let us also mention that our choice of scaling for 
the equation~(\ref{profileeq3}) is mass invariant, that is $G$ with density $\rho(G)$ satisfies the equation 
if and only if $G/\rho(G)$ satisfies the equation with $\rho=1$. Therefore all the estimates on the profiles 
are homogeneous in terms of the density $\rho$.  

  \begin{prop}\label{estimatesonGe} 
  Let us fix $\alpha_0 \in (0,1)$. There exists $a_1,a_2,a_3,a_4 \in (0,\infty)$ and, for any $k \in \N$,
  there exists $C_k \in (0,\infty)$ such that 
    \beqn\label{estimGeunif}\qquad
    \forall \, \alpha \in [\alpha_0,1), \,\, \forall \, G_\alpha \in \GG_\alpha, \quad
    \left\{ 
    \begin{array}{l}
   \displaystyle  \| G_\alpha \|_{L^1(  e^{a_1 \, |v|})} \le a_2, 
    \quad  \| G_\alpha \|_ {H^k(\RR^N)} \le C_k, \vspace{0.2cm} \\
     G_\alpha \ge a_3 \, e^{-a_4 \, |v|^{8}}. \vspace{0.2cm} \\ 
   
    \end{array}
    \right.  
    \eeqn
%As a consequence, the set $\GG$ is bounded in $H^k(m^{-1})$ with $m=\exp(-a|v|^s)$, 
%for any $k \in \N$, $a \in(0,\infty)$ and $s\in (0,1)$ {\bf (ou aussi $s=1$?)}: 
%there exists $C_{k,a,s} \in (0,\infty)$ such that 
%  \beqn\label{bddGG}
%  \sup_{\alpha \in [\alpha_0,1)} \sup_{G_\alpha \in \GG_\alpha} \| G_\alpha \|_{H^k(m^{-1})} \le C_{k,a,s}.
%  \eeqn
\end{prop}

We first recall the following geometrical lemma extracted (in a slightly specified form) from \cite[Lemma 2.3 \& Lemma 4.4]{MMRI}, that we shall use several times in the sequel.

\begin{lem}\label{jacobvprim} 
For any $\alpha \in (0,1]$ and $\sigma \in \Sp^{N-1}$ we define 
\bear \nonumber %\label{v*tovprim}
&&\phi^*_\alpha = \phi^*_{\alpha,v,\sigma} : \R^N \to \R^N, \quad v_* \mapsto  v' 
\\ \nonumber %\label{vtovprim}
&&\phi_\alpha = \phi_{\alpha,v_*,\sigma} : \R^N \to \R^N, \quad v \mapsto  v' 
\eear
and the Jacobian functions $J^*_\alpha = \hbox{{\em det}} \, (D \,  \phi^*_{\alpha,v,\sigma})$, 
$J_\alpha = \hbox{{\em det}} \, (D \, \phi_{\alpha,v_*,\sigma})$, as well as the cone 
  \[ 
  \Omega_\delta = \Omega_{\delta,\sigma} =  \big\{ u \in \R^N, \,\, \hat u  \cdot \sigma > \delta - 1 \big\},
  \] 
for any $\delta \in (0,2)$ and $\sigma \in \Sp^{N-1}$. 

For any $\delta \in (0,2)$, $\phi^*_\alpha$ defines a $C^\infty$-diffeomorphism 
from $v + \Omega_\delta$ onto $v + \Omega_{\omega^*(\delta)}$ 
with $\omega^*(\delta) = 1+ \sqrt{\delta/2}$ and $\phi_\alpha$ defines 
a $C^\infty$-diffeomorphism from $v_* + \Omega_\delta$ 
onto $v_* + \Omega_{\omega_\alpha(\delta)}$ with 
$$
\omega_\alpha(\delta) = 1+\frac{\delta-1+r_\alpha}{\Big(1+2 (\delta-1) r_\alpha + r_\alpha ^2\Big)^{1/2}}
$$ 
and $r_\alpha = (1+\alpha)/(3-\alpha)$. 

Moreover, there exist $C \in (0,\infty)$ such that with $C_\delta = C/\delta$
  \bear\label{vtovprim2}
  && C^{-1}_\delta \, |v-v_*| \le |\phi_\alpha (v) - v_*| \le 2 \, |v-v_*|, \\ \label{vtovprim3}
  && | \phi^{-1}_\alpha(v') - \phi^{-1}_{\alpha'}(v') | \le 
      C_\delta \, |\alpha'-\alpha| \, |v'-v_*|,  \\ \label{vtovprim4}
  && |J_\alpha | \le C_\delta, \quad |J_\alpha ^{-1}| \le 
      C_\delta, \quad  |J^{-1}_\alpha - J^{-1}_{\alpha'} | \le  C_\delta^2 \, |\alpha'-\alpha|
  \eear
on $v_*+\Omega_\delta$, uniformly with respect to the parameters 
$\alpha,\alpha' \in [0,1]$, $\sigma \in \Sp^{N-1}$ and $v_* \in \R^N$. 
The same estimate holds for $\phi^*_\alpha$ on $v+\Omega_\delta$. 
Finally, for any $\alpha,\alpha' \in [0,1]$, $\sigma \in \Sp^{N-1}$, $v_* \in \R^N$ 
and $t \in [0,1]$, there holds
  \beqn\label{vtovprim5}
  t \, \phi^{-1}_\alpha + (1-t) \,  \phi^{-1}_{\alpha'}  = \phi^{-1}_{\alpha_t} 
  \eeqn
for some $\alpha_t$ belonging to the segment with extremal points $\alpha$ and $\alpha'$. 
The same result holds for $\phi^*_\alpha$. 
\end{lem}

We will also need the following elementary result in order to estimate the 
convolution operator $L$ defined in (\ref{defL}). 

\begin{lem}\label{estimatesonL} 
For any function $g \in L^1_3(\R^N)$ there exists some constants 
$c_1, c_2 \in (0,\infty)$ such that 
  \beqn\label{estimateonL1}
  c_1 \, (1 + |v|) \le L (g) \le c_2 \, (1 + |v|).
  \eeqn
Moreover, if $g$ satisfies $\EE(g) \ge a_1 \, \rho$ and ${\bf m}_{3/2}(g) \le a_2 \, \rho$, 
for some constants $a_1, a_2>0$, we can take $c_1 = C^{-1} \, \rho$, $c_2 = C \, \rho$ 
in (\ref{estimateonL1}) for some explicit constant $C>0$ depending only on $a_1, a_2 >0$. 
\end{lem}

\smallskip\noindent{\sl Proof of Lemma \ref{estimatesonL}.} 
The upper bound in~(\ref{estimateonL1}) is immediate. 
As for the lower bound, we have, on the one hand, by Jensen's inequality,
\begin{equation}\label{minorL1}
\int_{\R^N} g_* \, |u| \, dv_* \ge \rho \, |v|. 
\end{equation}
On the other hand, by triangular inequality, 
$$
\int_{\R^N} g_* \, |u| \, dv_*  \ge {\bf m}_{1/2} - |v| \, {\bf m}_0. 
$$
By H\"older's inequality we have ${\bf m}_{1/2} \ge \EE^2 \,  {\bf m}_{3/2} ^{-1} \ge C_0 \, \rho$ for some 
explicit constant $C_0>0$ depending only on $a_1, a_2$.  
As a consequence 
\begin{equation}\label{minorL2}
\int_{\R^N} g_* \, |u| \, dv_*  \ge \rho \, (C_0 - |v|). 
\end{equation}
These two lower bounds (\ref{minorL1}, \ref{minorL2}) imply immediately that 
$$ 
\int_{\R^N} g_* \, |u| \, dv_*  \ge C ^{-1} \, \rho \, (1+|v|). 
$$
for some explicit constant $C>0$ depending only on $C_0$. \qed

\medskip\noindent{\sl Proof of Proposition \ref{estimatesonGe}.} 
We split the proof into several steps. 
In Steps~1, 2 and 3, we establish the smoothness for any profile $G_\alpha \in \GG $ 
as well as upper and lower bounds on its tail. 
In Steps~4, 5, 6, 7, 8 and 9, we show that these estimates actually are uniform with 
respect to the choice of the profile $G_\alpha \in \GG_\alpha$ and $\alpha \in [\alpha_0,1)$. 
Thanks to Steps~1, 2 and 3 the computations then performed are rigorously justified. 

We fix $\alpha \in [\alpha_0,1)$ and $G_\alpha$ a solution of (\ref{profileeq3})
for which we will establish the announced bounds. 
From now we omit the subscript ``$\alpha$'' when no confusion is possible. 

\medskip\noindent
{\it Step~1. Moment bounds. } 
From~\cite[Proposition~3.1]{MMII}, by taking $g_{\mbox{{\scriptsize in}}} = G$ 
in the evolution equation~(\ref{eqresc}), we get that $G \in L^1_k$ for any $k \in \NN$. 

\medskip\noindent
{\it Step~2. $L^2$ {\em a posteriori} bound. } 
We aim to prove that $G \in L^2$. Let us fix $A>0$ and let us 
introduce the $C^1$ function
  \[ \Lambda_A(x) := {x^2 \over 2} \, {\bf 1}_{x \le A} 
            + \left( A \, x - {A^2 \over 2}\right) \, {\bf 1}_{x > A}. \]
We multiply the equation (\ref{profileeq3}) by $\Lambda_A'(G) = \min \{ G,A \} := T_A(G)$. 
Once again we shall omit the subscript ``A'' when no confusion is possible. 
After some straightforward computation we get
\bean
&& \int_{\R^N} \Big(T(G) \, G \, L(G) + \rho \, (1-\alpha) \, N \, T(G)^2 /2 \Big) \, dv = 
 \int_{\R^N}  T(G) \, Q^+(G,G) \, dv. 
 \eean 
Since $L(G) \ge c_1 \, (1+|v|)$ thanks to Lemma~\ref{estimatesonL} and $\Lambda(G) \le G \, T(G)$ we have 
  \bear\label{TA1}
    c_1 \int_{\R^N} \Lambda(G) \, (1+|v|) \, dv &\le&  \int_{\R^N} T(G) \, G \, L(G) \, dv \\ \nonumber
  &\le& \int_{\R^N}  T(G) \, Q^+(G,G) \, dv \le I_{1} + I_{2} +I_{3} + I_{4},
  \eear
where the terms $I_k$ are defined in the following way, splitting the 
collision kernel into some smooth and non-smooth parts. 
Let $\Theta: \R \rightarrow \R_+$ be an even $C^\infty$ function such
that $\mbox{support} \, \Theta \subset (-1,1)$, and $\int_\R \Theta = 1$. 
Let $\widetilde{\Theta}: \R^N \rightarrow \R_+$ be a radial
$C^\infty$ function such that $\mbox{support} \, \widetilde{\Theta} 
\subset B(0,1)$ and $\int_{\R^N} \widetilde{\Theta} = 1$. Introduce
the regularizing sequences
  \[
  \Theta_m (z) = m \, \Theta(mz), \quad z\in \R, \qquad
  \widetilde{\Theta}_n (x)  =   n^N \widetilde{\Theta}(n x), \quad x\in\R^N.
  \]
As a convention, we shall use subscripts $S$ for ``smooth''
and $R$ for ``remainder''. We denote $\Phi(u) := |u|$. First, we set
  $$ %\label{PhiSn}
  \Phi_{S,n} = \widetilde{\Theta}_n \ast \left( \Phi \ {\bf 1}_{{\cal A}_n} \right), \qquad
  \Phi_{R,n} = \Phi - \Phi_{S,n}, 
  $$
where ${\cal A}_n$ stands for the annulus ${\cal A}_n = \left\{ x \in \R^N \ ; \ \frac{2}{n} \le |x| \le n \right\}$.
Similarly, we set
  $$ %\label{bSm}
  b_{S,m}(z) = \Theta_m \ast \left( b \ {\bf 1}_{{\cal I}_m} \right) (z), 
  \qquad b_{R,m} = b - b_{S,m},
  $$
where ${\cal I}_m$ stands for the interval ${\cal I}_m = \left\{ x
\in \R \ ; \ -1+\frac{2}{m} \le |x| \le 1-\frac{2}{m} \right\}$
($b$ is understood as a function defined on $\R$ with compact
support in $[-1,1]$). 
We then define
  $$
  I_1 =  \int_{\R^N}  T(G) \, Q_R^+(G,G) \, dv ,
  $$
where $Q_R^+$ is the gain term associated to the cross-section $B_R := |u|\, b_{R,m}$, 
  $$
  I_2 =  \int_{\R^N}  T(G) \, Q_{RS}^+(G,G) \, dv ,   
  $$
where $Q_{RS}^+$ is the gain term associated to the cross-section $B_{RS} := \Phi_{R,n}\, b_{S,m}$, 
  $$ 
  I_3 =  \int_{\R^N}  T(G) \, \Big[ \tilde Q_S^+(\chi(G),G) + \tilde Q_S^+(T(G),\chi(G) \Big] \, dv ,
  $$
where $Q_S^+$ is the gain term associated to the 
smooth cross-section $B_S :=  \Phi_{S,n} \, b_{S,m}$ and  $\chi(G) := G - T(G)$ and finally
  $$
  I_4 =  \int_{\R^N}  T(G) \, Q_S^+(T(G),T(G)) \, dv.
  $$
We estimate each term separately. We omit the subscripts $m$ and $n$ when there is no confusion. 
For $I_1$ we proceed along the line of the proof of the estimate for the term $I^r$ in 
\cite[Proof of Theorem 2.1]{MMRI}. 
Using Young's inequality $x \, T(y) \le \Lambda(x) + \Lambda(y)$ we have
  \bean
  I_1 &=& \int \!\! \int \!\! \int_{\R^N \times \R^N \times \Sp^{N-1}}
  G \, G_* \, T(G') \, b_{R,m} \, |u| \,  dv \, dv_* \, d\sigma \\
  &\le & \int \!\! \int \!\! \int_{\R^N \times \R^N \times \Sp^{N-1}}
  G \, [\Lambda(G_*) + \Lambda(G')] \, b_{R,m} \, {\bf 1}_{\hat u \cdot \sigma \le 0} 
     \, |u| \, dv \, dv_* \, d\sigma \\
  &&+ \int \!\! \int \!\! \int_{\R^N \times \R^N \times \Sp^{N-1}}
  G_* \, [\Lambda( G ) + \Lambda(G')] \, b_{R,m} \, {\bf 1}_{\hat u \cdot \sigma \ge 0} 
     \, |u| \, dv \, dv_* \, d\sigma 
  = I_{1,1} + ... + I_{1,4}.
  \eean
We just deal with the term $I_{1,2}$, the others may be handled in a similar (or even simpler) way. 
Making the change of variables $v_*\to v' = \phi^*_\alpha(v_*)$ (for some fixed $v,\sigma$) 
and using the elementary inequality $|u| \le 4 \, |v'-v|$ valid when $\sigma \cdot \hat u \le 0$, 
there holds
\bean
I_{1,2} &= & \int \!\! \int \!\! \int_{\R^N \times \R^N \times \Sp^{N-1}}
G \,  \Lambda(G')\, b_{R,m} \, {\bf 1}_{\hat u \cdot \sigma \le 0} \, |u| \, dv \, dv_* \, d\sigma \\
&=& 2^{N+2} \int \!\! \int \!\! \int_{\R^N \times \R^N \times \Sp^{N-1}}
G \,  \Lambda(G')\, b_{R,m} \, {\bf 1}_{\hat u \cdot \sigma \le 0} \, |v-v'| \, dv' \, dv \, d\sigma \\
&\le& 2^{N+2} \, \| b_{R,m}  \|_{L^1} \, \| G \|_{L^1_1} \, \int_{\R^N} \Lambda (G) \, (1+|v|) \, dv.  
\eean
Since the same estimates hold for all the terms $I_{1,k}$, we obtain
\beqn\label{I1TA}
I_1 \le \eps(m) \,  \| G \|_{L^1_1} \, \int_{\R^N} \Lambda (G) \, \langle v \rangle \, dv
\quad\hbox{ with }\quad \eps(m) \mathop{\longrightarrow}_{ m \to \infty} 0.
\eeqn

\smallskip
For $I_2$ we proceed along the line of the proof of the estimate for the term 
$I$ in \cite[Proof of Proposition 2.5]{MMII}. Using again Young's  
inequality $x \, T(y) \le \Lambda(x) + \Lambda(y)$ and the trivial 
estimate $\Phi_{R,n} \le C \, n^{-1} \, (|v|^2 + |v_*|^2)$ we get
\bean
I_2 &=& \int \!\! \int \!\! \int_{\R^N \times \R^N \times \Sp^{N-1}}
G \, G_* \, T(G') \, b_{S,m} \, \Phi_{R,n} \,  dv \, dv_* \, d\sigma \\
&\le & {C \over n} \int \!\! \int \!\! \int_{\R^N \times \R^N \times \Sp^{N-1}}
G \, |v|^2  \, [\Lambda(G_*) + \Lambda(G')] \, b_{S,m} \, dv \, dv_* \, d\sigma \\
&&+  {C \over n}  \int \!\! \int \!\! \int_{\R^N \times \R^N \times \Sp^{N-1}}
G_* \, |v_*|^2 \,  [\Lambda( G ) + \Lambda(G')] \, b_{S,m} \,  dv \, dv_* \, d\sigma 
= I_{2,1} + ... + I_{2,4}.
\eean
Because of the truncation on $b$ of frontal and grazing collisions,  
both changes of variables $v \to v' = \phi_\alpha(v)$ (for fixed $v_*,\sigma$) and $v_* \to v' = \phi^*_\alpha(v_*)$ 
(for fixed $v,\sigma$) are allowed (and the jacobian of their inverse is bounded). Hence in a similar way as 
for the term $I_1$ we obtain
  \beqn\label{I2TA}
  I_2 \le {C(m) \over n} \,  \| G \|_{L^1_2} \, \int_{\R^N} \Lambda (G)  \, dv.
  \eeqn

\smallskip
For $I_3$, using again Young's inequality, plus $T(G) \le G$ and 
the fact that both changes of variables 
$v \to v' = \phi_\alpha(v)$ (for fixed $v_*,\sigma$) and $v_* \to v' = \phi^*_\alpha(v_*)$ 
(for fixed $v,\sigma$) are allowed, we have 
\bean
I_3 &=& \int \!\! \int \!\! \int_{\R^N \times \R^N \times \Sp^{N-1}}
( T(G') + T(G'_*) ) \, \big[G \, \chi(G_*) + \chi(G) \, T(G_*)\big] \, b_{S,m} \, \Phi_{S,n} \, dv \, dv_* \, d\sigma \\
&\le & C (n) \int \!\! \int \!\! \int_{\R^N \times \R^N \times \Sp^{N-1}}
\Big\{ \chi(G_*) \, \big[\Lambda(G') + \Lambda(G'_*) + 2 \, \Lambda(G)\big]  \\
&& \hspace{5cm}+ \chi(G) \, \big[\Lambda(G') + \Lambda(G' _*) + 2 \, \Lambda(G_*)\big] \Big\} 
\, b_{S,m} \, dv \, dv_* \, d\sigma.
\eean
We deduce as before 
\beqn\label{I3TA}
I_3 \le C_{m,n} \,  \| \chi(G) \|_{L^1} \, \int_{\R^N} \Lambda (G)  \, dv 
\eeqn
for some constant $C_{m,n} >0$. 

\smallskip
Finally for $I_4$, we argue as in the proof of \cite[Proposition 2.6]{MMII} for the 
treatment of the term involving $Q^+ _S$, and we get for some $\theta \in (0,1)$
\beqn\label{I4TA}
I_4 \le C_{m,n} \, \| T(G) \|_{L^1}^{1+2\theta}  \, \| T (G)  \|_{L^2}^{2 \, (1-\theta)} , 
\eeqn
for some constant $C_{m,n} >0$. 

Gathering (\ref{TA1}),  (\ref{I1TA}), (\ref{I2TA}), (\ref{I3TA}), (\ref{I4TA}) 
and taking $m$, next $n$ and finally $A \ge A(G)$ large enough we may control the terms 
$I_1$, $I_2$ and $I_3$ by the half of the left hand side term of (\ref{TA1}) (for $I_3$ we use 
that  $\| \chi_A(G) \|_{L^1} \to 0$ when $A \to \infty$). Note that the condition 
$A \ge A(G)$ depends on the distribution $G$ (by the mean of some non-concentration bound), but shall 
play no role since we shall take the limit $A \to +\infty$ in the end. We obtain 
$$
\forall \, A \ge A(G), \qquad {c_1\over 2} \, \int_{\R^N} \Lambda_A (G) \, (1+|v|) \, dv \le 
 C_{b,\rho,\EE(G)} \, \| T_A (G)  \|_{L^2}^{2 \, (1-\theta)}
$$
for some constant $C_{b,\rho,\EE(G)}>0$ depending on the cross-section $b$ and on the 
profile $G$ {\it via} its energy.  
Using that $T_A(G)^2 /2 \le \Lambda_A (G)$ we deduce 
$$ 
\forall A \ge A(G), \qquad {c_1\over 4} \, \| T_A (G)  \|_{L^2}^{2 \, \theta} \le 
 C_{b,\rho,\EE(G)}
$$
and we then conclude that $G \in L^2$ passing to the limit $A \to \infty$ in the preceding estimate, with 
the bound 
  \begin{equation}\label{bornel2}
  \| G \|_{L^2} \le \left( \frac{4 \, C_{b,\rho,\EE(G)}}{c_1} \right)^{\frac{1}{2\theta}}. 
  \end{equation}

\begin{rem}\label{remarkL2}  
Note that the $L^2$ bound (\ref{bornel2}) only depends on the distribution $G$ by the mean of the energy $\EE(G)$ and the constant $c_1$. Therefore, thanks to Lemma~\ref{estimatesonL}, this bound only depends on a lower bound on the energy $\EE(g)$ and an upper bound on the third moment ${\bf m}_{3/2}(g)$. \end{rem}

\medskip\noindent
{\it Step 3. Smoothness and positivity. } 
Thanks to \cite[Theorem 1.3]{MMII} and \cite[Theorem 1]{BGP**}, 
taking $g_{\mbox{{\scriptsize in}}} = G$ as an initial condition 
in (\ref{eqresc}) we have that $G$ belongs to the Schwartz space 
of $C^\infty$ functions decreasing faster than any polynomials, and that 
$G \ge a_1\, e^{-a_2 \, |v|}$ for some constant $a_1,a_2 >0$. 

\smallskip
So far the estimates in Step~3 may be not uniform 
on the elasticity coefficient $\alpha \in [\alpha_0,1)$ and on the profile $G_\alpha$. 
The aim of the following steps is to prove that they actually are uniform. 
Note however that estimates of the previous steps shall ensure that 
the following computations are rigorously justified. 

\medskip\noindent
{\it Step 4. Upper bound on the energy using the energy dissipation term. } 
We prove that 
\beqn\label{UpperEnergy}
\forall \, \alpha \in (0,1] \qquad \EE \le {4 \over b_1^2} \, \rho.
\eeqn
From equation~(\ref{EnergyDEGe}) on the energy of the profile $G$ there holds
\beqn\label{dissipenergyGe}
 (1+\alpha) \, b_1  \, \int_{\RR^N} \! \int_{\RR^N} G \, G_* \, |u|^3 \, dv \, dv_* = 
 2 \, \rho \int_{\R^N} G \, |v|^2 \,  dv. 
\eeqn
From Jensen's inequality
  $$
  \int_{\R^N} |u|^3 \, G_* \, dv_* \ge \rho \, |v|^3,
  $$
and H\"older's inequality
  $$
  \int_{\R^N} |v|^3 \, G \, dv \ge \rho^{-1/2} \, \left( \int_{\R^N} |v|^2 \, G \, dv \right)^{3/2}, 
  $$
we get 
  $$
  (1+\alpha) \, b_1 \, \rho^{1/2} \, \EE^{3/2} \le 2 \, \rho \, \EE 
  $$ 
from which the bound (\ref{UpperEnergy}) follows. 

\medskip\noindent
{\it Step 5. Lower bound on the energy using the entropy.} We prove 
  \beqn\label{LowerEnergy}
  \forall \, \alpha \in (0,1] \qquad \EE \ge {N \, \alpha^4 \over 8} \, \rho.
  \eeqn
  
\begin{rem}
The choice of scaling we have made for the evolution equation in self-similar variable 
becomes clear from this computation: it is chosen 
such that the energy of the self-similar profile does not blow up nor vanishes for $\alpha \to 1$. 
The restriction $\alpha \in [\alpha_0,1)$, $\alpha_0 > 0$, is then made in order to get a uniform 
estimate from below on the energy. 
\end{rem}

\smallskip
By integrating the equation satisfied by $G$ against $\log G$ we find   
  $$
  \int_{\R^N} Q(G,G) \log G \, dv 
  - \rho \, (1-\alpha) \,  \int_{\R^N} \log G  \,\,  \nabla_v \cdot (v \, G) \, dv = 0.
  $$
Then we write the first term as in~\cite[Section~1.4]{GPV**} to find 
  \bean
  &&{1 \over 2} \int\!\! \int\!\! \int_{\R^{2N}\times S^{N-1}} G \, G_* \left( \log {G'G'_* \over G G_*} 
  - {G'G'_* \over G G_*} + 1 \right) \, B \, dv \, dv_* \, d\sigma \\
  &&\qquad + {1 \over 2} \int\!\! \int\!\! \int_{\R^{2N}\times S^{N-1}}   \left( G'G'_* -G G_* \right) \, B \, dv \, dv_* \, d\sigma  
  + \rho \,(1-\alpha) \, \int_{\R^N} v \cdot \nabla_v G  \, dv = 0. 
  \eean
If we denote 
\beqn\label{defDHa}
D_{H,\alpha} (g) = {1 \over 2} \int\!\! \int\!\! \int_{\R^{2N}\times S^{N-1}} g \, g_* \left( {g'g'_* \over g g_*} -\log {g'g'_* \over g g_*} 
     -  1 \right) \, B \, dv \, dv_* \, d\sigma \ge 0,
\eeqn
(recall that in this formula the post-collisional velocities $v',v'_*$ are computed according to the 
inelastic formula  (\ref{vprimvprim*}) with normal restitution coefficient $\alpha \in (0,1]$),
 we can write
  \beqn\label{DissipGe} 
  - D_{H,\alpha}(G) + \left( {1 \over \alpha^2} - 1 \right) \, b_2 \, \int \!\! \int_{\R^{2N}}  G \, G_* \, |u| \, dv \, dv_* 
  - (1-\alpha) \, N \, \rho^2 = 0, 
  \eeqn
with $b_2 := \| b \|_{L^1}$,  and thus we get
  $$
  \int \!\! \int_{\R^{2N}} G \, G_* \, |u| \, dv \, dv_* =  
  {\alpha^2 \over 1+\alpha}  \left( N \, \rho^2  + {1 \over 1-\alpha} D_{H,\alpha} (G) \right) 
  \ge  { N \, \alpha^2  \over 2} \, \rho^2.  
  $$
On the other hand, from Cauchy-Schwarz's inequality
\bean
&& \int \!\! \int_{\R^{2N}} G \, G_* \, |u| \, dv \, dv_* \le \\
&& \le \left( \int \!\!  \int_{\R^{2N}} G \, G_*  \, dvdv_* \right)^{1/2} \!\! 
 \left( \int \!\!  \int_{\R^{2N}} G \, G_* |u|^2 \, dvdv_* \right)^{1/2} = \sqrt{2}  \, \rho^{3/2}\, \EE^{1/2}, 
\eean
and then the bound (\ref{LowerEnergy}) follows gathering the two preceding estimates. 

\medskip\noindent
{\it Step 6. Upper bound on (exponential) moments using Povzner inequality.} 
There exists $A,C > 0$ such that 
  $$ %\label{Upperexponential}
  \forall \, \alpha \in [0,1), \qquad \int_{\R^N}  G (v) \, e^{A |v|} \, dv \le C \, \rho.
  $$

We refer to \cite{BGP**} where that bound is obtained as an immediate consequence 
of the following sharp moment estimates: there exists $X > 0$ such that 
  \beqn\label{Uppermoment}
  \forall \, \alpha \in [0,1), \qquad {\bf m}_k = \int_{\R^N} G \, |v|^k \, dv \le \Gamma(k+1/2) \, X^{k/2} \, \rho.
  \eeqn
It is worth noticing that in \cite{BGP**} the Povzner inequality used 
in order to get (\ref{Uppermoment}) is uniform in the normal restitution  
coefficient $\alpha \in [0,1]$ and that the factor $\rho$ comes from our 
choice of the scaling variables (in which $\rho$ is involved). 

\medskip\noindent
{\it Step 7. Uniform upper bound on the $L^2$ norm.} 
From (\ref{LowerEnergy}), (\ref{Uppermoment}) and Remark~\ref{remarkL2}, the $L^2$ bound (\ref{bornel2}) is uniform on  $\alpha \in [\alpha_0,1) $ and $G \in \GG_\alpha$. 

\medskip\noindent
{\it Step 8. Smoothness.}
It is enough to show some uniform bounds from above and below on the energy 
together with uniform non-concentration bounds on the 
self-similar profiles in $\GG$, in the form of upper bounds on the $L^2$ 
bounds for instance. Indeed the proofs 
of~\cite[Proposition~3.1, Proposition~3.2, Proposition~3.4, Theorem~3.5 and Theorem~3.6]{MMII} 
then apply straightforwardly (in these proofs we did not use the 
part associated with the anti-drift in the semigroup). Therefore the uniform bounds 
on the $H^k$ norms for all $k \ge 0$ follows from these results. 

\medskip\noindent
{\it Step 9. Pointwise Lower bound.}  It is a consequence of the following lemma. \qed

\begin{lem}\label{LowerBdga} Let $g \in C([0,\infty);L^1_3)$ be a solution of the rescaled equation (\ref{eqresca})
with inelasticity parameter $\alpha \in (0,1)$ and assume that for some $p > 1$ and $C,T \in (0,\infty)$
$$
\sup_{[0,T]} \| g \|_{L^p \cap L^1_3} \le C.
$$
(i) For any $t_1\in (0,T)$ there exists $a_1 \in (0,\infty)$ (depending on $C$, $\rho$ and $t_1$ but not on $T$) such that 
\beqn\label{concluLowerBdga}
 \forall \, t \in [t_1,T], \,\,\, \forall \, v \in \R^N, \quad g(t,v) \ge a_1^{-1} \, e^{-a_1 \, |v|^8}.
\eeqn
(ii) If furthermore, $g_{\mbox{\scriptsize{{\em in}}}}$ satisfies 
 $$ %\label{assumgdecomp2}
 g_{\mbox{\scriptsize{{\em in}}}} (v) \ge a_0^{-1} \, e^{-a_0 \, |v|^8},
 $$
then (\ref{concluLowerBdga}) holds with $t_0 = 0$ and some constant $a_1 \in (0,\infty)$ (depending on $C$, $\rho, a_0$ but not on $T$).
\end{lem}

\smallskip\noindent{\sl Proof of Lemma \ref{LowerBdga}.} We only prove (i), the proof of (ii) being similar. 
Let us fix $t_1 \in (0,1)$. 
We closely follow the proof of the Maxwellian lower bound for the solutions of the elastic Boltzmann equation (see \cite{Ca32,PW97}) taking advantage of some technical results established in its extension to  the solutions of the inelastic Boltzmann equation (see \cite[Theorem 4.9]{MMII}).
The starting point is again the evolution equation satisfied by $g$  written in the form
$$
\partial_t g + \tau_\alpha \, v \cdot \nabla_v g + \big(\tau_\alpha \, N + C + C \, |v|\big)  \, g = 
Q^+_\alpha(g,g) + \big(C + C \, |v| - L(g)\big)  \, g, 
$$
where the last term in the right hand side term is non-negative for some well-chosen 
numerical constant $C \in (0,\infty)$ thanks to Lemma~\ref{estimatesonL}, 
(\ref{Uppermoment}) and (\ref{LowerEnergy}). 
Let us introduce the semigroup $U_t$ associated to the operator 
$\tau_\alpha \,  v \cdot \nabla_v + \lambda(v)$, 
where $\lambda (v) := \tau_\alpha \, N + C + C \, |v|$, which action is  given by
\[ 
(U_t \, h) (v) =  h(v \, e^{-\tau_\alpha \, t}) \, \exp \left( - \int_0^t \lambda(v\,
e^{-s}) \, ds \right).
\]
Thanks to the Duhamel formula, we have
\beqn\label{DuhamelF}
\forall \, t > 0, \,\, \forall \, \tau \ge 0, \quad g(t+\tau,.) \ge  \int_0^t U_{t-s} Q^+(g(s+\tau,.),g(s+\tau,.)) \,ds.
\eeqn
Noticing that 
$$ %\label{int0tlambda}
\left( - \int_0^t \lambda(v\, e^{-s}) \, ds \right) \ge  - (C \, |v|+ \tau_\alpha \, N \, t + C \, t), 
$$
and repeating the arguments of Steps 2 and 3 in the proof of \cite[Theorem 4.9]{MMII}, 
 we get that 
\beqn\label{Getak}
\forall \, t \ge \tau, \quad g(t,.) \ge \eta \, {\bf 1}_{B(0,\delta)}(v)
\eeqn
with $\tau = \tau_1 = t_1/2$ and some constant $\eta = \eta_1> 0$, $\delta= \delta_1>1$. Let us emphasize that here we make use of  Lemma 4.6, Lemma 4.7 and Lemma 4.8 in  \cite{MMII} where the constants exhibited in these ones are uniform in $\alpha \in [\alpha_0,1)$ thanks to the uniform $L^p \cap L^1_3$ estimates assumed on $g$.  

\smallskip
Now, on the one hand, from  \cite[ Lemma 4.8]{MMII}, there exists $\kappa \in (0,\infty)$ such that 
$$
Q^+_\alpha ({\bf 1}_{B(0,1)},{\bf 1}_{B(0,1)}) \ge \kappa \, {\bf 1}_{B(0,\sqrt{5}/2)}
$$
which in turns implies 
\beqn\label{Q+delta}
\forall \, \delta > 0, \quad 
Q^+_\alpha ({\bf 1}_{B(0,\delta)},{\bf 1}_{B(0,\delta)}) \ge \kappa \, \delta^{-N-1} \, {\bf 1}_{B(0,\sqrt{5}/2 \, \delta)}.
\eeqn
On the other hand, there exists $\kappa' \in (0,\infty)$ such that 
\beqn\label{Ssdelta}
\forall \, \delta > 0, \,\,\,\,
\forall \, s \in [0,1], \quad U_s( {\bf 1}_{B(0,\delta)}) \ge \kappa' \, e^{- C \, \delta} \, {\bf 1}_{B(0,\delta)}.
\eeqn
From (\ref{Getak}) with $\eta= \eta_1$, $\delta= \delta_1$, and making use of (\ref{DuhamelF}),  (\ref{Q+delta}),  (\ref{Ssdelta}), we get that (\ref{Getak}) holds  with 
$$
\tau = \tau_2 = \tau_1 + {t_1 \over 2^2}, \quad
\delta = \delta_2 ={\sqrt{5} \over 2}  \, \delta_1 \quad\hbox{and}\quad
\eta= \eta_2 = (\tau_2 - \tau_1) \, \kappa'' \, \eta_1^2 \, e^{-C' \, \delta_1} ,
$$
where $\kappa'' = \kappa \, \kappa'$ and $C'$ depends on $C$ and $N$.  Iterating the argument we get that (\ref{Getak}) holds with $\tau = \tau_k = \tau_{k-1} + t_1 \, 2^{-k} = (1-2^{-k}) \, t_1$, $\delta = \delta_{k+1} = \left( {\sqrt{5} / 2} \right)^{k+1}$ and 
$$
\eta _{k+1} = (\kappa'' \, t_1)^{1+2+...+2^{k-1}} \, \eta_1^{2^{k}} \, e^{- C' \, (\delta_k + 2 \, \delta_{k-1} + ... + \, 2^{k-1} \, \delta_1)}
\, 2^{-[k + 2 \, (k-1) + ... + 2^{k-1} \, 1]} \ge A^{2^{k+1}},
$$
with $A := \sqrt{\kappa'' \, t_1} \, \sqrt{\eta_1} \, e^{-C' \, \delta_1}/2$.
In other words, using that $\left( {\sqrt{5} \over 2} \right)^{8} > 2$, we have proved
$$ %\label{Getak2}
\forall \, t \ge t_1, \ \forall \, k \in \N, \quad g(t,v) \ge A^{2^k} \, {\bf 1}_{B(0,2^{k/8} \, \delta_1)}(v),
$$
from which we easily conclude. \qed

%%%%%%%%%%%%%%%%%%%%%%%%%%%%%%%%%%%%%%%%%%%%%%%%
\subsection{Estimates on the difference of two self-similar profiles}

In this subsection we take advantage of the mixing effects of the collision 
operator in order to show that the $L^1$ norm of their difference of two self-similar profiles 
(corresponding to the same inelasticity coefficient) indeed controls  
the $H^k \cap L^1(m^{-1})$ norm of their difference for any $k \in \N$ and for 
some exponential weight function $m$, uniformly in terms of $\alpha \in [\alpha_0,1)$. 

 \begin{prop}\label{errorH1} 
 For any $k >0$, there is $m=\exp(-a \, |v|)$, $a \in (0,\infty)$ and 
 $C_{k}>0$ such that for any $\alpha \in [\alpha_0,1)$ 
 and any $G_\alpha,H_\alpha \in \GG_\alpha$ there holds
   \beqn\label{eq:errorH1}
   \|H_\alpha - G_\alpha  \|_{H^k \cap L^1(m^{-1})} \le C_{k} \, \|H_\alpha - G_\alpha \|_{L^1}. 
   \eeqn
 \end{prop}

\smallskip\noindent{\sl Proof of Proposition \ref{errorH1}. }
We proceed in three steps. It is worth mentioning that all the constants in the proof are uniform in terms of the  normal restitution coefficient $\alpha \in [\alpha_0,1)$, as they only depend on the 
uniform bounds of Proposition~\ref{estimatesonGe} and some uniform bounds on the 
collision kernel. 

\medskip\noindent
{\sl Step 1. Control of the $L^1$ moments.} 
We prove first that there exists $A, C \in (0,\infty)$ such that 
  $$
  \forall \, \alpha \in [\alpha_0,1), \qquad  \int_{\R^N} |H_\alpha-G_\alpha| \, e^{A \, |v|} \, dv 
      \le C  \int_{\R^N} |H_\alpha-G_\alpha| \, dv.
  $$
Let us consider some normal restitution coefficient $\alpha \in [\alpha_0,1)$ and 
two self-similar profiles $G,H \in \GG_\alpha$ (here again, we omit the subscript $\alpha$ when there is no confusion). We denote 
$D=G-H$, $S=G+H$ and $\varphi = |v|^{2p} \, \mbox{sgn}(D)$, 
$p \in {1\over2} \, \N$, $p \ge 3/2$, where $\mbox{sgn}(D)$ denotes the sign of $D$. 
The equation for $D$ reads 
  \bear
  0 &=& Q_\alpha (G,G) - Q_\alpha (H,H) - \rho \, (1-\alpha) \, \nabla_v \cdot ( v \, D) \nonumber \\ \label{QDS}  
  &=& 2 \, \tilde Q_\alpha(D,S) - \rho \, (1-\alpha) \, \nabla_v \cdot (v \, D ).
  \eear 
Multiplying equation (\ref{QDS}) by  $\varphi$, we get 
\bean
0 &=& \int_{\R^N \times \R^N \times \Sp^{N-1}} 
    B \, D \, S_* \, \Big[ \varphi'_* + \varphi' - \varphi_* - \varphi  \Big] \, dv \, dv_* \, d\sigma \\
&& \hspace{4cm} - \rho \, (1-\alpha) \, \int_{\R^N} \nabla_v(vD) \, |v|^{2p} \, \mbox{sgn}(D) \, dv \\
  &\le& \int_{\R^N \times \R^N}  |u| \, |D|\, S_* \, K_p\, dv \, dv_*+ 2 \, 
         \int_{\R^N \times \R^N} |u| \, |D| \, S_* \, |v_*|^{2p} \, dv \, dv_* \\
  && \hspace{4cm} + \rho \, \int_{\R^N} |D| \, v \cdot \nabla (|v|^{2p}) \, dv 
  \eean 
with 
$$ %\label{mp2}
K_p(v,v_*) :=  \int_{\Sp^{N-1}} (|v'|^{2p} + |v'_*|^{2p} - |v|^{2p} - |v_*|^{2p}) \, b(\sigma\cdot u) \, d\sigma.
$$
From ~\cite[Corollary 3, Lemma 2]{BGP**}, there holds
  $$ %\label{mp3}
  K_p(v,v_*) \le \gamma_p \, \Sigma_p -  (1-\gamma_p) \, ( |v|^{2p} + |v_*|^{2p})
  $$
where $(\gamma_p)_{p=3/2,2,...}$ is a decreasing sequence of real numbers such that 
  \beqn \label{mp4}
0 <  \gamma_p < \min\left\{ 1, {4 \over p+1}\right\},
  \eeqn
and $\Sigma_p$ is defined by  
  \[ \Sigma_p := \sum_{k=1}^{k_p} \pmatrix{p \cr k} \left(|v|^{2k} \, |v_*|^{2p-2k} + |v|^{2p-2k}
  \, |v_*|^{2k} \right), \]
with $k_p := [(p+1)/2]$ is the integer part of $(p+1)/2$ and 
$\pmatrix{p \cr k}$ stands for the binomial coefficient.
As a consequence, 
  \bean
  &&(1-\gamma_{3/2}) \, \int_{\R^N \times \R^N} |v|^{2p} \, |u| \, S_* \, |D| \, dv \, dv_* \le
  \gamma_p \, \int_{\R^N \times \R^N}  |u| \, |D|\, S_* \, \Sigma_p\, dv \, dv_* \\
  &&\qquad + 2 \, 
         \int_{\R^N \times \R^N} |u| \, |D| \, S_* \, |v_*|^{2p} \, dv \, dv_* 
         + 2 \, \rho \, p \,  \int_{\R^N} |D| \, |v|^{2p} \, dv .
  \eean
Using Lemma~\ref{estimatesonL} in order to estimate $L(S)$ from below, the inequality $|u| \le |v| + |v_*|$ and introducing the notations
$$
{\bf d}_k := \int_{\RR^N} |D| \, |v|^{2k} \, dv, \quad {\bf s}_k := \int_{\RR^N} S \, |v|^{2k} \, dv, 
$$
we get, for some numerical constant $C \in (0,\infty)$, 
\beqn\label{dp1/2}
{\rho \over C} \, {\bf d}_{p+1/2}  \le  \gamma_p \, S_p 
+ \left({\bf d}_0 \, {\bf s}_{p+1/2} + {\bf d}_{1/2} \, {\bf s}_p \right) + 2 \, \rho \, p \, {\bf d}_p,
\eeqn
with 
$$
S_p :=  \sum_{k=1}^{k_p} \pmatrix{p \cr k} 
\left({\bf d}_{k+1/2} \, {\bf s}_{p-k} + {\bf d}_{k} \, {\bf s}_{p-k+1/2} 
+ {\bf d}_{p-k+1/2} \, {\bf s}_{k}+ {\bf d}_{p-k} \, {\bf s}_{k+1/2}\right).
$$
From Proposition~\ref{estimatesonGe}, or more precisely (\ref{Uppermoment}), we know that ${\bf s}_k \le \rho \, \Gamma(k+1/2) \, x^k$ for any $k \ge 1$ and for some $x \in (1,\infty)$. By H\"older's inequality, we also have 
\[ 
  {\bf d}_p^{1+{1 \over 2p}} \le {\bf d}_{p+{1\over2}} \, {\bf d}_0^{{1 \over 2p}}. 
\]
Repeating the proof of \cite[Lemma~4]{BGP**}, for any $a \ge 1$, there exists $A > 0$ such that
  \[  
S_p \le A \,   \rho \, ({\bf d}_0+{\bf d}_{1/2})\, \Gamma (a\,p+a/2 + 1) \, Z_p 
  \]
with
\[ 
Z_p := \max_{k=1,..,k_p} \{ {\bf \delta}_{k+1/2} \, {\bf \sigma}_{p-k},  {\bf \delta}_{k} \, {\bf \sigma}_{p-k+1/2}, {\bf \delta}_{p-k+1/2} \, {\bf \sigma}_{k},  {\bf \delta}_{p-k} \, {\bf \sigma}_{k+1/2} \}, 
\]
and
\[
{\bf \delta}_k := {{\bf d}_k \over ({\bf d}_0+{\bf d}_{1/2})\,  \Gamma(a\,k+1/2)}, \quad
{\bf \sigma}_k := {{\bf s}_k \over \rho \, \Gamma(a\,k+1/2)}. 
\]
We may then rewrite~(\ref{dp1/2}) as
\bean
\Gamma(a\,p+1/2)^{1/2p} \, {\bf \delta}_p^{1+1/2p} \le 
A \, \gamma_p  \, {\Gamma (a\,p+a/2 + 1)  
  \over  \Gamma(ap+1/2)} \, Z_p + ({\bf \sigma}_{p+1/2} + {\bf \sigma}_p) +  2 \, \rho \, p \, {\bf \delta}_p.
\eean
On the one hand, from (\ref{mp4}), there exists $A'$ such that
  $$ %\label{mp8}
  A \, \gamma_p \, {\Gamma(ap+a/2+1) \over  \Gamma(ap+1/2)}  \le A' \, p^{a/2-1/2} \qquad 
  \forall \, p=3/2,2,\dots
  $$
On the other hand, thanks to Stirling's formula $n! \sim n^n \, 
e^{-n} \, \sqrt{2\pi n}$ when $n\to\infty$ and the estimate (\ref{mp4}), 
there exists $A'' > 0$ such that
  $$ %\label{mp9}
  (1-\gamma_p) \, \Gamma(a\,p+1/2)^{1/2p} \ge A'' \, p^{a/2}  \qquad \forall 
  \,  p = 3/2,2, \dots
  $$
Therefore,
$$
p^{a/2} \, {\bf \delta}_p^{1+1/2p} \le p^{a/2-1/2} \, Z_p + ({\bf \sigma}_{p+1/2} + {\bf \delta}_1 \, {\bf \sigma}_p) +  2 \, \rho \, p \, {\bf \delta}_p.
$$
We finally obtain
$$
{\bf d}_k \le x^k \, \Gamma(a k +1/2) \, ({\bf d}_0+{\bf d}_{1/2}),
$$
and we easily conclude as in  \cite[Proof of Theorem 1]{BGP**} or in \cite[Proof of Proposition~3.2, Step 2]{MMRI}.

\medskip\noindent
{\sl Step 2. Control of the $L^2$ norms.} 
For $k=0$, the propagation of the $L^2$ norm is immediate using the 
result~\cite[Corollary~2.3]{MMII}. Indeed one just has to split the collision 
kernel as in~\cite[Section~2.4]{MMII}. For the truncated and regularized part 
$Q^+ _S$ (we use the notation introduced in step 2 the proof of Proposition~\ref{estimatesonGe}), 
\cite[Corollary~2.3]{MMII} together with some basic interpolation yield the 
following control:
$$
\int_{\R^N} \left( Q^+ _{S} (S,D) + Q^+ _{S} (D,S) \right) D \, dv 
\le C \, \rho^{1 + 2 \theta} \, \|D\|_{L^2} ^{2-2\theta} 
$$ 
for some explicit $C>0$ and $\theta \in (0,1)$. For the remaining 
term $Q^+ _R$, we use the same control as in~\cite[Proof of Proposition~2.5]{MMII} to get 
$$ 
\int_{\R^N} \left( Q^+ _R (S,D) + Q^+ _R (D,S) \right) \, D \, dv 
\le \e \, \left( \|D\|_{L^1 _2} + \|D\|_{L^2 _{1/2}} \right) \|D\|_{L^2 _{1/2}}
$$ 
for some $\e$ which can be taken as small as wanted by the truncation. Gathering these estimates, we get
$$ %\label{estimL2Q+SD}
\forall
 \, \eps > 0 \qquad
\int_{\R^N}  \tilde Q^+ (S,D) \,  D \, dv \le \eps \, \| D \|_{L^2_{1/2}}^2 + C_\eps
$$
where $C_\eps$ depends on weighted $L^1$ and $L^2$ norms of $S$, on $L^1$ norms on $D$ and on $\eps$.  
Using equation~(\ref{eq:Ds}) with $i=0$, Lemma~\ref{estimatesonL} to treat 
the term $L(D)$, and some elementary interpolation, we deduce that 
$$ 
\|D\|_{L^2 _{1/2}} \le C \, \|D\|_{L^1 _2} 
$$ 
for some constant $C >0$, which concludes the proof for $k=0$ using the previous 
step on the $L^1$ moments. 

\medskip\noindent
{\sl Step 3. Control of the $H^k$ norms.} 
From the previous step and some interpolation, in order to conclude it is enough 
to prove~(\ref{eq:errorH1}) for any $k \in \N$ and $m \equiv 1$. 
We proceed by induction on $k$. 
For any $i \in \N^N$, the equation satisfied by $\partial^i D$ is 
\bean
&&\partial^i Q^+(S,D) + \partial^i Q^+(D,S) - \partial^i (L(D) \, S) 
- L(S) \, \partial^i D \\
&& \qquad - \sum_{0 < i' \le i} \left( \begin{array}{c} i' \\ i \end{array} \right) 
 \partial^{i'} L(S) \, \partial^{i-i'} D 
- \rho \, (1-\alpha) \, \partial^i \nabla \cdot( v \,  D) = 0.
\eean
We deduce that 
\bear \label{eq:Ds}
&& C \, \int_{\R^N} (\partial^i D)^2 \, (1+|v|) \, dv \le 
\int_{\R^N} \left( \partial^i Q^+(S,D) + \partial^i Q^+(D,S) \right) \partial^i D \, dv \\ \nonumber
&& \qquad - \sum_{0 \le i' \le i} \left( \begin{array}{c} i' \\ i \end{array} \right) 
\int_{\R^N} \partial^{i'} L(D) \, \partial^{i-i'} S \, \partial^i D \, dv \\ \nonumber
&& \qquad \qquad - \sum_{0 < i' \le i} \left( \begin{array}{c} i' \\ i \end{array} \right) 
\int_{\R^N} \partial^{i'} L(S) \, \partial^{i-i'} D \, \partial^i D \, dv 
\eear
droping the non-positive term. 

The induction is initialized by Step 2. Let us assume the induction step $k \ge 0$ 
to be proved, and let us consider some $i \in \N^N$ such that $|i|=k+1$. 
Using equation~(\ref{eq:Ds}) and~\cite[Theorem~2.5]{MMII} to estimate the gain term, we find easily 
$$ 
\| \partial^i D\|_{L^2} \le C \, \left( \|D\|_{L^1 _q} + \|D\|_{H^{k+(3-N)/2} _q} \right)  
$$ 
for some $q>0$. Therefore we obtain by interpolation (since 
$(3-N)/2 <1$ for $N \ge 2$), for another $q'$ possibly larger: 
$$ 
\| D\|_{H^{k+1}} \le C \, \left( \|D\|_{L^1 _{q'}} + \|D\|_{H^{k}} \right).
$$ 
This concludes the proof, using interpolation, the induction hypothesis $k$, and the 
Step~1 on the $L^1$ moments.
\qed

%%%%%%%%%%%%%%%%%%%%%%%%%%%%%%%%%%%%%%%%%%%%%%%%

\section{The elastic limit $\alpha \to1$}
\setcounter{equation}{0}
\setcounter{theo}{0}

%%%%%%%%%%%%%%%%%%%%%%%%%%%%%%%%%%%%%%%%%%%%%%%%

\subsection{Dependency of the collision operator according to the inelasticity}

In this subsection we show that the collision operator depends continuously 
on the inelasticity coefficient $\alpha \in [0,1]$. 
Since it is an unbounded  operator, this continuous dependency is expressed in the norm of the graph of 
the operator or in some weaker norm. We start showing that this dependency of the collision 
operator is Lipschitz, and even $C^{1,\eta}$ for any $\eta \in (0,1)$, when allowing a loss (in terms of derivatives
and  weight) in the norm they are expressed. Let define the formal derivative of the collision 
 operator according to $\alpha$ by 
 $$ %\label{defKK}
 Q'_\alpha(g,f) :=  \nabla_v \cdot \left( 
 \int_{\R^N} \!\! \int_{\Sph^{N-1}} g('v_*(\alpha)) \, f('v(\alpha))
 \, b \,  |u| \, \left( {u-|u| \, \sigma  \over 4 \,  \alpha^2} \right) \, d\sigma \, dv_*\right)  
 $$
 or by duality 
 $$ %\label{defdualKK}
 \langle Q'_\alpha(g,f) , \psi \rangle := 
 \int_{\R^N} \!\!  \int_{\R^N} \!\! \int_{\Sph^{N-1}} g_* \, f \, 
 b \,  |u| \, \left( {|u| \, \sigma - u  \over 4} \right)  \, 
 \nabla \psi (v'_\alpha) \,  d\sigma \, dv_* \, dv.
 $$

\begin{prop}\label{ContQe1} 
Let us fix a smooth exponential weight $m = \exp ( - a \, |v|^s )$, $a \in (0,+\infty)$, $s \in (0,1)$.
Then 
\begin{itemize}
\item[(i)] For any $k,q \in \N$ the exists $C \in (0,\infty)$ such that for any smooth functions $f, g$
(say in $\Ss(\R^N)$) and any $\alpha \in [0,1]$ there holds
\bear \label{QaWk1}
&&\big\|  Q^\pm_{\alpha}(g,f) \big\|_{W^{k,1}_q(m^{-1})} 
 \le C_{k,m} \, \|f \|_{W^{k,1} _{q+1} (m^{-1})} \, \| g \|_{W^{k,1} _{q+1} (m^{-1})} \\ \label{Q'aWk1}
 &&\big\|  Q'_{\alpha}(g,f) \big\|_{W^{k,1}_q(m^{-1})} 
 \le C_{k,m} \, \|f \|_{W^{k+1,1} _{q+2} (m^{-1})} \, \| g \|_{W^{k+1,1} _{q+2} (m^{-1})}.
 \eear
\item[(ii)] Moreover,  for any smooth functions $f, g$ and for any $\alpha, \alpha' \in [0,1]$, there holds
  \bear \nonumber
&& \quad \big\| Q^+_\alpha(g,f) - Q^+_{\alpha'}(g,f)- (\alpha-\alpha') \, Q'_\alpha(g,f) \big\|_{W^{-2,1}_q(m^{-1})}  \\ \label{Q'}
&& \qquad \qquad   \le |\alpha-\alpha'|^2\,    \|f \|_{L^1_{q+3} (m^{-1})} \,  \| g \|_{L^1_{q+3}(m^{-1})}.
 \eear
\item[(iii)] As a consequence, there holds
 \begin{equation}   \label{convLetoL1} 
 \mbox{ } \quad \big\| Q^+_{\alpha'}(g,f) - Q^+_{\alpha}(g,f) \big\|_{W^k _q(m^{-1})}
 \le C \, |\alpha-\alpha'| \,  \|f \|_{W^{2k +3,1}_{q+3} (m^{-1})} \,  
  \| g \|_{W^{2k +3,1}_{q+3} (m^{-1})},
 \end{equation}
% {\bf ou "seulement"
%$$
% \big\| Q^+_{\alpha'}(g,f) - Q^+_{\alpha}(g,f) \big\|_{L^1(m^{-1})}
% \le C \, |\alpha-\alpha'| \,  \|f \|_{W^{3,1}_{2}(m^{-1})}^{1/2} \, \| f \|_{L^1_{3}(m^{-1})}^{1/2} \, 
%  \| g \|_{W^{3,1}_{2}(m^{-1})}^{1/2} \,  \| g  \|_{L^1_{3}(m^{-1})}^{1/2}
%  $$
%?}
and for any $\eta \in (1,2)$, there exists $k_\eta \in \N$, $q_\eta \in \N$ and $C_\eta \in (0,\infty)$ such that 
  \bear \nonumber
&& \quad \big\| Q^+_\alpha(g,f) - Q^+_{\alpha'}(g,f)- (\alpha-\alpha') \, Q'_\alpha(g,f) \big\|_{L^1(m^{-1})}  \\ \label{expLatoL1}
&& \qquad \qquad  
  \le C_\eta \, |\alpha-\alpha'|^\eta \, 
  \|f \|_{W^{k_\eta,1} _{q_\eta}(m^{-1})} \,  \| g \|_{W^{k_\eta,1}_{q_\eta}(m^{-1})}.
 \eear
\end{itemize}
\end{prop}

\smallskip \noindent
{\sl Proof of Proposition~\ref{ContQe1}.}  
First by classical convolution-like estimates (see for instance~\cite{MV**} in the elastic case, 
and~\cite{GPV**} in the inelastic case, as well as the proof of Proposition~\ref{ContQe2} below) we easily have (\ref{QaWk1}) and (\ref{QaWk1}). 

Next, in order to prove (\ref{Q'}) we proceed  by duality. Let us consider $\varphi \in \Ss(\R^N)$ and define $\psi := \varphi \, \langle v \rangle^q \, m^{-1}$. 
We  compute
\bean
I &:=& \int_{\R^N} \left[Q^+_\alpha(g,f) - Q^+_{\alpha'}(g,f) -(\alpha - \alpha') \, Q' _\alpha(g,f)\right] \, \psi(v) \, dv \\
 & = &\int_{\RR^N \times \RR^N \times \Sph^{N-1}} |u| \, b \, g_* \, f
      \, \left[  \psi(v'_\alpha) -  \psi(v'_{\alpha'}) - (\alpha - \alpha') \, 
      \left( \frac{|u| \, \sigma -u}{4} \right) \cdot \nabla \psi (v' _\alpha) \right] 
      \, dv \, dv_* \, d\sigma.
\eean
Hence, if one denotes by $\xi_{v,v_*,\sigma}(\alpha):=\psi(v' _\alpha)$ 
(for given fixed values of $v,v_*,\sigma$), we obtain (omitting the subscripts for clarity)
\bean
|I | 
&=& \int_{\RR^N \times \RR^N \times \Sph^{N-1}} |u| \, b \, g_* \, f
      \, \left[  \xi(\alpha) - \xi(\alpha') - (\alpha - \alpha') \, \xi' (\alpha) \right] \, dv \, dv_* \, d\sigma \\
&\le& (\alpha - \alpha')^2 \,  \int_{\RR^N \times \RR^N \times \Sph^{N-1}} |u| \, b \, g_* \, f
      \, \sup_{\alpha \in (0,1)} | \xi''(\alpha)| \, dv \, dv_* \, d\sigma.
\eean
We then easily conclude that (\ref{Q'}) holds using that $\langle v' \rangle^q \, (m')^{-1} \le C \, \langle v \rangle^q \, (m)^{-1} \, \langle v_*\rangle^q \, (m_*)^{-1}$ for some constant $C \in (0,\infty)$. 

Last, we prove (\ref{convLetoL1}) by using the following interpolation on 
$J =  Q^+_\alpha(g,f) - Q^+_{\alpha'}(g,f) -(\alpha - \alpha') \, Q' _\alpha(g,f)$: 
$$
\left\| J \right\|_{W^{k,1}_q (m^{-1})} 
\le \left\| J \right\|_{W^{-2,1} _q(m^{-1})} 
\, \left\| J \right\|_{W^{2(k+1),1} _q (m^{-1})} 
$$
and using (\ref{Q'}) on the first term in the right-hand side, and (\ref{QaWk1},\ref{Q'aWk1}) on the second 
term in the right-hand side. It yields 
$$
\left\| Q^+_\alpha(g,f) - Q^+_{\alpha'}(g,f) -(\alpha - \alpha') \, Q' _\alpha(g,f) \right\|_{W^{k,1}_q (m^{-1})}  
\le C \, |\alpha - \alpha'| \,  \,  \|f \|_{W^{2 k +3,1}_{q+3} (m^{-1})} \,  \| g \|_{W^{2k +3,1}_{q+3} (m^{-1})}
$$
and (\ref{convLetoL1}) follows by using (\ref{Q'aWk1}) again.

Then the proof of (\ref{expLatoL1}) is done in the same way using suitable interpolation. \qed         

\medskip
We next state a mere (H\"older) continuity dependency on $\alpha$, which is however stronger than Proposition~\ref{ContQe1} in some sense, since it is  written in the norm of the graph of the operator for one the argument. 
 
\begin{prop}\label{ContQe2} 
For any $\alpha,\alpha' \in (0,1]$,  and any $g \in L^1_1(m^{-1})$, 
$f \in W^{1,1}_1 (m^{-1})$, there holds
 \beqn\label{convLetoL2}
 \quad \left\{
  \begin{array}{l}
  \big\|Q^+_\alpha(g,f) - Q^+_{\alpha'}(g,f) \big\|_{L^1(m^{-1})}
                 \le \eps (\alpha-\alpha') \, \| f \|_{W^{1,1}_{1}(m^{-1})} \, 
                 \|g \|_{L^1_1(m^{-1}) }, \vspace{0.3cm} \\
  \big\|Q^+_\alpha(f,g) - Q^+_{\alpha'} (f,g) \big\|_{L^1(m^{-1})} 
                 \le \eps (\alpha-\alpha') \, \| f \|_{W^{1,1}_{1}(m^{-1})} \, \|g  \|_{L^1_1(m^{-1}) }.
  \end{array}
  \right.
  \eeqn
where $\eps(r) = C \, r^{{1 \over 3 + 4 / s}}$ for some constant $C$ (depending only on $b$). 
\end{prop}

\smallskip\noindent{\sl Proof of Proposition \ref{ContQe2}.}
For any given $v,v_* \in  \R^N$, $w = v+v_* \not = 0$ and $\sigma \in S^{N-1}$ we define 
$\chi \in [0,\pi/2], \,\,\, \cos \,\chi := |\sigma \cdot \hat w|$.
Let us fix $\delta \in (0,1)$, $R \in (1,\infty)$ and let us define $\theta_\delta \in W^{1,\infty}(-1,1)$ such that 
$\theta_\delta (s) = 1$ on $(-1+2\delta,1-2\delta)$, $\theta_\delta (s) = 0$ on $(-1+\delta,1-\delta)^c$, $0 \le \theta_\delta \le 1$, $|\theta'_\delta(s)| \le 3/\delta$,  $\Theta_R(u) = \Theta(|u|/R)$ with $\Theta(x) = 1$ on $[0,1]$, $\Theta(x) = 1-x$ for $x \in [1,2]$ and $\Theta(x) = 0$ on $[2,\infty)$, $A(\delta)  := \{ \sigma \in S^{N-1}; \,\, \sin^2 \chi \ge \delta \}$,  $B(\delta) := \{ \sigma \in S^{N-1}; \,\, \cos \theta \in (-1+2\delta,1-2\delta)^c $ or $ \sin^2 \chi \le \delta \}$. 
We then split $Q^+$ in three terms, namely
$$
Q^+_\alpha = Q^{+,a}_\alpha + Q^{+,v}_\alpha + Q^{+,r}_\alpha
$$
where $Q^{+,r}_\alpha$ is defined by (\ref{Qinel}) with $b$ replaced by $b^r := b \, \theta_\delta(\sigma \cdot \hat u) \,\Theta_R(u)$, where $Q^{+,v}_\alpha$ is defined by (\ref{Qinel}) with $b$ replaced by $b^v := b \, {\bf 1}_{A(\delta)} \, (1-\Theta_R(u))$ and where
$Q^{+,a}_\alpha$ is defined by (\ref{Qinel}) with $b$ replaced by $b^a := b \, (1 - \theta_\delta(\sigma \cdot \hat u)) \, \Theta_R(u) + b \, (1-\Theta_R(u)) \, {\bf 1}_{A^c(\delta)}$. We split the proof into three steps. 

\smallskip\noindent
{\sl Step 1. Treatment of small angles. } 
There exists a constant $C \in (0,\infty)$ such that for any $\alpha \in (0,1]$ and $\delta \in (0,1)$ there holds
$$ %\label{ContQe21}
\big\|Q^{+,a}_\alpha (\psi,\varphi) \big\|_{L^1(m^{-1})} \le 
C \, \delta \,  \| \psi \|_{L^1_1(m^{-1})} \|\varphi \|_{L^1_1(m^{-1}) }.
$$
Indeed let us consider some $\ell \in L^\infty$ and let us proceed by duality. We estimate
\bean
\int_{\R^N} Q^{+,a}_\alpha (\psi,\varphi)\, \ell(v) \, m^{-1}(v) \, dv
&=& \int_{\R^N \times \R^N \times \Sp^{N-1}} |u| \,  b^a \, \psi_* \, \varphi \, \ell' \, (m')^{-1} \, dv \, dv_* \, d\sigma \\
&\le& \| b^a \|_{L^\infty_{v,v_*}(L^1(\Sp^{N-1}))} \, \| \ell \|_{L^\infty} \,  \| \psi \|_{L^1_1(m^{-1})} \|\varphi \|_{L^1_1(m^{-1}) },
\eean
and we conclude using that 
$\| b^a \|_{L^\infty_{v,v_*}(L^1(\Sp^{N-1}))} \le C \, ( \delta + \max_{v,v_*} | B(\delta) | ) \le C \, \delta$. 

\smallskip\noindent
{\sl Step 2. Treatment of large relative velocities. } 
There exists a constant $C = C_{a,s,b} \in (0,\infty)$ such that for any $\alpha \in (0,1]$ and $\delta \in (0,1)$ there holds
\beqn\label{ContQe22}
\|Q^{+,v}_\alpha (\psi,\varphi) \|_{L^1(m^{-1})} \le {C  \over R \, \delta^{2/s}} \,  \| \psi \|_{L^1_1(m^{-1})} \, 
\|\varphi \|_{L^1_1(m^{-1}) }
\eeqn
We need the following lemma, which we state below and prove at the end of the subsection. 

\begin{lem}\label{mvprim} For any $\delta > 0$ and $\alpha \in (0,1)$,  there holds
\beqn\label{mvprim1}
 \sigma \in S^{N-1}, \,\,\, \sin^2 \chi \ge \delta
\quad\hbox{implies}\quad
m^{-1}(v') \le m^{-k}(v) \, m^{-k}(v_*),
\eeqn
with $k = (1-\delta/160)^{s/2}$. 
\end{lem}

\medskip
In order to prove (\ref{ContQe22}) we fix $\ell \in L^\infty$ and we argue by duality again. 
We estimate thanks to Lemma~\ref{mvprim}
\bean
\int_{\R^N} Q^{+,v}_\alpha (\psi,\varphi)\, \ell(v) \, m^{-1}(v) \, dv
&=& \int_{\R^N \times \R^N \times \Sp^{N-1}} |u| \,  b^v \,\psi_* \, \varphi \, \ell' \, (m')^{-1} \, dv \, dv_* \, d\sigma \\
&\le& {1 \over R} \int_{\R^N \times \R^N \times \Sp^{N-1}} |u|^2 \,   b^v \, \psi_* \, \varphi \, \ell' \, (m)^{-k} \,  (m_*)^{-k} \, dv \, dv_* \, d\sigma \\
&\le &{1 \over R}\,  \| \ell \|_{L^\infty} \,  \|\psi \|_{L^1_2(m^{-k})} \|\varphi \|_{L^1_2(m^{-k}) } \\&\le &{1 \over R}\,  \| \ell \|_{L^\infty} \, \| |.| \, m^{1-k}(.) \|_{L^\infty}^2 \,  \|\psi \|_{L^1_1(m^{-1})} \|\varphi \|_{L^1_1(m^{-1}) }, 
\eean
from which we easily conclude since $x\mapsto x \, m^{1-k}(x)$ is uniformly bounded by $C_{a,s} \, (1-k)^{-1/s}$, $C_{a,s} \in (0,\infty)$.  

\smallskip\noindent
{\sl Step 3. The truncated operator. } 
Let us prove that there exists a constant $C  \in (0,\infty)$ such that for any $\delta \in (0,1)$,  
$\alpha, \alpha' \in (0,1]$ and $R \in (1,\infty)$ there holds
$$ %\label{ContQe23}
 \|Q^{+,r}_\alpha(g,f) - Q^{+,r}_{\alpha'}(g,f) \|_{L^1(m^{-1})}
                 \le C \, |\alpha-\alpha'| \, \left({R^2 \over \delta} + {R  \over \delta^3}\right)  \, \| g \|_{L^1(m^{-1})} 
                 \, \| f \|_{W^{1,1}(m^{-1}) }.
$$                
We closely follow the proof of \cite[Proposition 4.3]{MMRI}.
We consider some $\ell \in L^\infty$, $f,g \in \DD(\R^N)$, we proceed by duality and next conclude thanks to a density argument.
We have 
  \bean 
I&:=& \int_{\R^N} [Q^{+,r}_\alpha (g,f) - Q^{+,r}_{\alpha'} (g,f) ] \, m^{-1} \, \ell \, dv   \\
&=& \int_{\R^N \times \R^N \times \Sp^{N-1}} |u| \, \Theta_R(u) \, b_\delta \, g_* \, f \, 
    \left[ \ell(v'_\alpha) \, m^{-1}(v'_\alpha) -  \ell(v'_{\alpha'}) \, m^{-1}(v'_{\alpha'}) \right] \, dv \, dv_* \, d\sigma.
  \eean
With the notations of Lemma~\ref{jacobvprim} we perform the changes of variables $v \mapsto v'_\alpha = \phi_{\alpha}(v)$ and 
$v \mapsto v'_{\alpha'} = \phi_{\alpha'}(v)$ (for fixed $v_*$ and $\sigma$) with jacobians  $J_\alpha$ and $J_{\alpha'}$. Observing that without restriction we may assume $\alpha \le \alpha'$ and therefore  $\OO_{\alpha} =  v_* + \Omega_{\omega_\alpha(\delta)}\subset \OO_{\alpha'} =  v_* + \Omega_{\omega_{\alpha'}(\delta)}$ since $s \mapsto \omega_s(0)$ is an increasing function, we get 
\bean
I &=&  \int_{\RR^N \times \Sph^{N-1}} \int_{\OO_{\alpha} \backslash \OO_{\alpha'}} g_* \, \ell'  \, (m^{-1})' \,F(\phi^{-1}_\alpha) \, J^{-1}_{\alpha} \, dv' \, dv_* \, d\sigma \\ 
&+&  \int_{\RR^N \times \Sph^{N-1}} \int_{\OO_{\alpha'}}  g_* \, \ell'  \, (m^{-1})' \,F(\phi^{-1}_\alpha) 
        \, \Big[J^{-1}_{\alpha} - J^{-1}_{\alpha'} \Big] \, dv' \, dv_* \, d\sigma \\  
&+& \int_{\RR^N \times \Sph^{N-1}} \int_{\OO_{\alpha'}}  
g_* \, \ell'  \, (m^{-1})' \, \Big[F(\phi^{-1}_\alpha) - F(\phi^{-1}_\alpha)\Big] 
      \, J^{-1}_{\alpha} \, dv' \, dv_* \, d\sigma  \\
&=& I_{1} + I_{2} + I_{3},
\eean
with $F(w) := |w-v_*| \, \Theta_R(w-v_*) \, f(w) \, b_\delta(\sigma \cdot \widehat{w-v_*}) $. For the first term $I_{1}$ we use the backward change of variables $v' \mapsto v = \phi^{-1}_{\alpha}(v')$ (for fixed $v_*$ and $\sigma$) and we get 
$$
I_{1} =  \int_{\RR^N \times \Sph^{N-1}} \int_{\R^N} |u| \, 
\Theta_R(u) \, f \,  g_* \, \ell(v'_\alpha)  \, m^{-1} (v'_\alpha) \,   b_\delta \, 
{\bf 1}_{ 0 \le \hat u \cdot \sigma \le  \eta }  \, dv_* \, dv \, d\sigma
$$
with $\eta := \omega_{\alpha}^{-1} \circ \omega_{\alpha'}(\delta) \le C\, \delta^{-3/2} \, |\alpha - \alpha'|$ for some constant $C \in (0,\infty)$. Since $v \mapsto |v|^{ s/2}$ is an increasing subadditive function, we also have $|v'_\alpha|^s \le (|v|^2 + |v_*|^2)^{s/2} \le |v|^s + |v_*|^s$, which implies $m(v'_\alpha) \le C \, m^{-1} \, m_*^{-1}$ for some constant $C \in (0,\infty)$ (depending of $\zeta$).  
As a consequence,  we obtain
\bean %\label{DepeI11}
|I_{1}| \le C \, R \, \delta^{-3/2}  \, |\alpha - \alpha'| \, \| b \|_{L^\infty} \,  \| \ell \|_{L^\infty} \,  \| f \|_{L^1(m^{-1})} \, \| g \|_{L^1(m^{-1})}.
\eean

\smallskip
For the term $I_{2}$, using  the backward change of variable  $v' \mapsto v = \phi^{-1}_{\alpha'}(v')$ (for some fixed $v_*$ and $\sigma$) and using the bounds (\ref{vtovprim4}) on $J_\alpha$ and $|J^{-1}_{\alpha} - J^{-1}_{\alpha'}|$, we obtain 
\bean %\label{DepeI12}
|I_{1}| \le C \, R \, \delta^{-3} \, |\alpha - \alpha'| \, \| b \|_{L^\infty} \,  \| \ell \|_{L^\infty} \,  \| f \|_{L^1(m^{-1})} \, \| g \|_{L^1(m^{-1})}.
\eean

\smallskip
In order to estimate $I_3$, we introduce $\alpha_t := (1-t) \, \alpha + t \, \alpha'$ and, thanks to (\ref{vtovprim3})-(\ref{vtovprim2}), we get 
$$
|I_{3} | \le {C \over \delta} \,  |\alpha - \alpha'| \, \int_0^1\!\! \int_{\RR^N \times \Sph^{N-1}} \!\!\int_{\OO_{\alpha_t}} |g_* | \, 
        \, | \ell' | \, (m^{-1})' \, |v'-v| \,  \Big|\nabla_w F (\phi^{-1}_{\alpha_t}(v') )\Big| \, dv'  dv_*  d\sigma dt.
$$
Using finally the backward change of variable  $v' \mapsto v = \phi^{-1}_{\alpha_t}(v')$ and the uniform bound (\ref{vtovprim4}) on $J_{\alpha_t}$, $t \in [0,1]$, on $v_* + \Omega_\delta$, we get 
$$
|I_{3} |
\le C \, \left({R^2 \over \delta} + {R  \over \delta^2}\right) \, |\alpha-\alpha'| \, \| b \|_{W^{1,\infty}} \,  \| \ell \|_{L^\infty} \,  \|g\|_{L^1(m^{-1})} \, \|f\|_{W^{1,1}(m^{-1})}.
$$
Gathering the estimates established in Steps 1, 2 and 3, we deduce the first inequality in (\ref{convLetoL2}). The second inequality in (\ref{convLetoL2}) is proved in a similar way (using symmetric changes of variable, 
allowed by the truncation).  \qed

\smallskip 
 \smallskip\noindent
{\sl Proof of Lemma~\ref{mvprim}. } We proceed in three steps. 

\smallskip\noindent
{\sl Step 1. }  
Assume first that $(2/\sqrt{5}) \, |v_*| \le |v| \le (\sqrt{5}/2) \, |v_*|$. Using the fact that $x \mapsto x^{s/2}$ is an increasing and subadditive function, there holds
$$
|v'|^s \le (|v|^2 + |v_*|^2)^{s/2} \le (9/4)^{s/2} \, |v_*|^s,
$$ 
and then by symmetry and because $s \le 1$
$$
|v'|^s \le {1 \over 2} \,  (9/4)^{s/2} \, (|v|^s + |v_*|^s) \le {3 \over 4} \, (|v|^s + |v_*|^s).
$$ 
In that case,  (\ref{mvprim1}) holds with $k = 3/4$. 

\smallskip\noindent
{\sl Step 2. } We shall first show that for any $v,v_* \in \R^N$ and $\sigma \in S^{N-1}$, there holds
\beqn\label{mvprim2}
\quad |v'|^2, \ |v'_*|^2  \le |v|^2 + |v_*|^2 - \frac{1+\alpha}8 \, \sin^2\chi \, |v+v_*|^2.
\eeqn
We recall the formula 
$$
v' := {v + v_* \over 2} + \frac12 \, \Bigg[ \frac{1-\alpha}2 \, u + \frac{1+\alpha}2 \, |u| \, \sigma \Bigg], \qquad
v'_*:= {v + v_* \over 2} + \frac12 \, \Bigg[ \frac{1-\alpha}2 \, u - \frac{1+\alpha}2 \, |u| \, \sigma \Bigg].
$$
Straightforward computations yield (denoting $S=v+v_*$)
$$
|v'|^2 \le \frac{|S|^2}4 + \frac14\, \Bigg[ \frac{1+\alpha^2}2 \, |u|^2 
             + \frac{1-\alpha^2}2 \, |u|^2 \, \cos \theta \Bigg] 
              + \frac{1-\alpha}4 \, (S \cdot u) + \frac{1+\alpha}4 \, |S| \, |u| \, \cos \chi. 
$$
We deduce the bound from above 
$$
|v'|^2 \le \frac{|S|^2}4 + \frac{|u|^2}4
              + \frac{1-\alpha}4 \, |S| \, |u| + \frac{1+\alpha}4 \, |S| \, |u| \, \cos \chi. 
$$
Then by applying twice Young's inequality
\bean
|v'|^2 &\le& |S|^2 \, \Big( \frac14 +\frac{1-\alpha}8 \Big) + |u|^2  \, \Big( \frac14 + \frac{1-\alpha}8 \Big)  
                  + \frac{1+\alpha}4 \, |S| \, |u| \, \cos \chi, \\
&\le& |S|^2 \, \Big( \frac14 +\frac{1-\alpha}8 \Big) 
           + |u|^2  \, \Big( \frac14 + \frac{1-\alpha}8 + \frac{1+\alpha}8 \Big)  
                  + \frac{1+\alpha}8 \, |S|^2 \, \cos^2 \chi, \\
&\le& \frac{|S|^2}2 + \frac{|u|^2}2   
                  + \frac{1+\alpha}8 \, |S|^2 \, (\cos^2 \chi -1), 
\eean
from which we deduce (\ref{mvprim2}). 
%In the same way we prove that (\ref{mvprim2}) holds with $v'$ replaced by $v'_*$. 

\smallskip\noindent
{\sl Step 3. } Assume that $\sin ^2 \chi \ge \delta$ and that either $(2/\sqrt{5}) \, |v_*| \ge |v|$ or $|v| \ge (\sqrt{5}/2) \, |v_*|$. In the first case, we have
$$
|v+v_*| \ge \big(1-(2/\sqrt{5})\big) \, |v_*| + (2/\sqrt{5}) \, |v_*| - |v| \ge \big(1-(2/\sqrt{5})\big) \, |v_*|,
$$ 
which then implies 
$$
|v+v_*| \ge  \big(1-(2/\sqrt{5})\big) \, (\sqrt{5}/2) \, |v| \ge \big(1-(2/\sqrt{5})\big) \, |v|. 
$$ 
The same inequalities are proved in a similar way in the second case. 
We deduce 
$$
|v+v_*|^2 \ge  {1 \over 2} \big(1-(2/\sqrt{5})\big) \, (|v|^2 + |v_*|^2). 
$$
We then deduce from (\ref{mvprim2}) that $|v'|^2 \le (1-\delta/160) \, (|v|^2 + |v_*|^2)$ and we conclude that (\ref{mvprim1}) holds as in Step~1. \qed

%%%%%%%%%%%%%%%%%%%%%%%%%%%%%%%%%%%%%%%%%%%%%%

\subsection{Quantification of the elastic limit $\alpha \to1$}

We begin with a simple consequence of Proposition~\ref{ContQe1}. 

\begin{cor}\label{DHeDH1} There exists $k_0, q_0 \in \N$ such that for any $a_i \in (0,\infty)$ $i=1, \, 2, \, 3$, 
there exists an explicit constant $C \in (0,\infty)$ such that for any function $g$ satisfying 
$$ %\label{gbddinfsup}
 \| g \|_{H^{k_0} \cap L^1_{q_0}} \le a_1, \quad g \ge a_2 \, e^{-a_3 \, |v|^8},
$$
there holds
$$
\big| D_{H,\alpha}(g) - D_{H,1}(g) \big| \le C \, (1-\alpha),
$$
where we recall that $D_{H,\alpha}$ is defined in (\ref{defDHa}). 
\end{cor}

\medskip
\noindent{\sl Proof of Corollary \ref{DHeDH1}.} We write 
\bean
D_{H,\alpha}(g) - D_{H,1}(g) &=&
\int\!\!\!\int\!\!\!\int b \, |u| \, \Big[g'_\alpha\, g'_{* \alpha} - g'\, g'_*] \, dv \, dv_* \, d\sigma \quad (=: I_1) \\
&+& \int\!\!\!\int\!\!\!\int b \, |u| \, g \, g_* \,  
\big[ \log g'_\alpha+ \log g'_{*\alpha} - \log g' - \log g'_* \big] \, dv \, dv_* \, d\sigma \quad (=: I_2). 
\eean
For the first term, thanks to Proposition~\ref{ContQe1}, we have
\bean
|I_1| &\le& \big\| Q^+_\alpha(g,g) - Q^+_1(g,g) \big\|_{L^1} 
\le C \, (1-\alpha) \, \| g \|_{W^{3,1}_{3}}^2. 
\eean
For the second term, we write 
\bean
|I_2| 
&=&2 \,  \Big|\big \langle (Q^+_\alpha(g,g) - Q^+_1(g,g)) \, \langle v \rangle^8, \langle v \rangle^{-8} \, \log g \big \rangle \Big| \\
&\le& 2 \, \big\| Q^+_\alpha(g,g) - Q^+_1(g,g) \big\|_{L^1_8} \, \| \langle v \rangle^{-8} \, \log g \|_{L^\infty} \\
&\le& C \, (1-\alpha) \, \| g \|_{W^{3,1}_{11}}^2 \,  
\big( |\log \, \| g\|_{L^\infty}| + |\log a_2| + a_3  \big) \\
&\le& C \, (1-\alpha) \, a_1 ^2 \, \big( |\log a_1 | + |\log a_2| + C \, a_3 \big), 
\eean 
thanks to Proposition~\ref{ContQe1} and the bounded embedding 
$H^{k_0} \cap L^1_{q_0} \subset L^\infty \cap W^{3,1}_{11}$ for $k_0,q_0$ large enough (see Proposition~\ref{interpolineg}).  We  conclude the proof gathering these two estimates. \qed

\medskip
Let us recall now two famous inequalities, namely the Csisz\'ar-Kullback-Pinsker inequality (see \cite{Csisz,Kul59}) and
the so-called entropy-entropy production inequalities (the version we present here is established in \cite{VillCerc}) that we will
 use several time in the sequel. 

\smallskip
  \begin{theo}\label{CSK&EEP} 
\begin{itemize}
\item[(i)] For a given function $g \in L^1_2$, let us denote by $M[g]$ the Maxwellian function 
with the same mass, momentum and temperature as $g$.
For any $0 \le g \in L^1_2(\R^N)$, there holds
\beqn\label{CKineg}
\big\| g - M[g] \big\|_{L^1}^2 \le 2 \, \rho(g) \,  \int_{\R^N} g \ln { g \over M[g]} \, dv.
\eeqn
\item[(ii)] For any $\eps > 0$ there exists $k_\eps, \, q_\eps \in \N$  and for any $A \in (0,\infty)$ there exists $C_\eps = C_{\eps,A} \in (0,\infty)$ such that for any $g \in H^{k_\eps} \cap L^1_{q_\eps}$ such that 
$$ %\label{EEPcv}
g(v) \ge A^{-1} \, e^{-A \, |v|^8}, \quad \| g \|_{ H^{k_\eps} \cap L^1_{q_\eps}} \le A,
$$
there holds
\beqn\label{EEPineg}
C_\eps \, \rho(g)^{1-\eps} \, \left( \int_{\R^N} g \ln { g \over M[g]} \, dv \right)^{1+\eps} \le D_{H,1}(g).
\eeqn
\end{itemize}
\end{theo}

\medskip
We have then the following estimate on the distance between $G_\alpha$ and $\bar G_1$ for any self-similar profile $G_\alpha$. 
\smallskip
  \begin{prop}\label{ConvGexplicit} 
  For any $\eps > 0$ there exists $C_\eps$ (independent of the mass $\rho$) such that 
    \beqn\label{estimateonGe2}
    \forall \, \alpha \in [\alpha_0,1) \qquad  \sup_{G_\alpha \in \GG_\alpha} 
    \, \|G_\alpha - \bar G_1 \|_{L^1_2} \le C_\eps \, \rho \, (1-\alpha)^{{1 \over 2+\eps}}
    \eeqn
where we recall that  $\bar G_1$ is the Maxwellian function defined by (\ref{Maxlim})--(\ref{tempSS}).
  \end{prop}

\medskip
\noindent{\sl Proof of Proposition \ref{ConvGexplicit}.} 
On the one hand, for any inelasticity coefficient $\alpha \in [\alpha_0,1)$ and 
profile $G_\alpha$, there holds from (\ref{DissipGe}) together 
with Corollary~\ref{DHeDH1} and the uniform estimates of Proposition~\ref{estimatesonGe}
\beqn\label{DH1<DHa}
D_{H,1}(G_\alpha)  \le D_{H,\alpha}(G_\alpha) + \rho^2 \, \OO(1-\alpha) \le \rho^2 \, \OO(1-\alpha).
\eeqn
On the other hand, introducing the Maxwellian function $M_\theta$ with the same mass, momentum and temperature 
as $G_\alpha$, that is $M_\theta$ given by (\ref{defMrut}) with $u=0$ and 
$\theta = \EE(G_\alpha)/\rho$, and gathering (\ref{DH1<DHa}), (\ref{EEPineg}), (\ref{CKineg}) 
with the uniform estimates of Proposition~\ref{estimatesonGe} and interpolation inequality, 
we obtain that for any $q,\eps > 0$ there exists $C_{q,\eps}$ such that 
\beqn\label{estimGeMtheta}
\forall \alpha \in [\alpha_0,1) \qquad \| G_\alpha - M_\theta \|_{L^1_q}^{2+\eps} \le C_{q,\eps} \, 
\rho^{2+\eps} \, (1-\alpha).
\eeqn
Next, from (\ref{dissipenergyGe}), we  have 
$$
 b_1  \, \int_{\RR^N} \! \int_{\RR^N} G_\alpha \, G_{\alpha*} \, |u|^3 \, dv \, dv_* - \rho \, \int_{\R^N} G_\alpha \, |v|^2 \,  dv = 
 (1-\alpha)  \,  {b_1 \over 2} \, \int_{\RR^N} \! \int_{\RR^N} G_\alpha \, G_{\alpha*} \, |u|^3 \, dv \, dv_*
$$
and then 
\beqn\label{ineqtheta1}
\left|  \Psi(\theta) \right| \le C_1 \, \| G_\alpha - M_\theta \|_{L^1_3} + C_2 \, \rho^2 \, (1-\alpha),
\eeqn
where we have used that $G_\alpha$ and $M_\theta$ are bounded thanks to Proposition~\ref{estimatesonGe} and we have defined 
\beqn\label{defPsi}
\Psi(\theta) = \rho \, \int_{\R^N} M_\theta \, |v|^2 \,  dv  
- b_1 \, \int_{\RR^N} \! \int_{\RR^N} M_\theta \, M_{\theta*} \, |u|^3 \, dv \, dv_*.
\eeqn
By elementary changes of variables, this formula simplifies into 
$$
\Psi(\theta) = k_1 \, \theta - k_2 \, \theta^{3/2}
$$
with $k_1 = \rho^2 \, N$ and, using (\ref{MMu3}),
$$
k_2 =   \rho^2 \, b_1 \, \int_{\RR^N \times \RR^N} M_{1,0,1} \, (M_{1,0,1}) _* \, |u|^3 \, dv \, dv_* =
2^{3/2} \, \rho^2 \, b_1 \, {\bf m}_{3/2}(M_{1,0,1}).
$$ 
We next observe that $\Psi \in C^\infty(0,\infty)$ and $\Psi$ is strictly concave. 
It is also obvious that the equation $\Psi(\theta) = 0$ for $\theta > 0$ has a unique solution 
which is $\bar\theta_1$ defined in~(\ref{tempSS}), and that we have 
$$
\Psi(\theta) \le \Psi'(\bar\theta_1) \, (\theta-\bar\theta_1) = - k_1 \, (\theta-\bar\theta_1)/2
$$ 
as well as
\beqn\label{propPsi}
\Psi(\theta) = \theta \, [k_1 - k_2 \, \theta^{1/2}] = k_2 \, \theta \,  [\bar\theta_1^{1/2} -  \theta^{1/2}].
\eeqn
Plugging this expression for $\Psi$ into~(\ref{ineqtheta1}) and using the lower bound (\ref{LowerEnergy})  on the temperature $\theta$ and the estimate (\ref{estimGeMtheta}) we obtain that for any $\eps > 0$ there is $C_\eps \in (0,\infty)$ such that 
\beqn\label{thetaa-b}
\forall \,  \alpha \in (\alpha_0,1) \qquad 
 \left| \theta^{1/2} -  \bar\theta_1^{1/2} \right| ^{2+\eps} \le C_\eps \, (1-\alpha). 
\eeqn
Namely, we have thus proved that the temperature of $\bar G_\alpha$ converge (with rate) to the expected temperature $\bar\theta_1$. In order to come back to the norm of $G_\alpha - \bar G_1$, we first write, using Cauchy-Schwarz's inequality,
\bear \nonumber
\| G_\alpha - \bar G_1 \|_{L^1_{-N}} 
&\le& \| G_\alpha - M_\theta \|_{L^1_{-N}} +  \| M_\theta - \bar G_1 \|_{L^1_{-N}} \\ \label{Ga-G1L1}
&\le&  \| G_\alpha - M_\theta \|_{L^1} + C_N \, \| M_\theta - \bar G_1 \|_{L^2}, 
\eear
and we remark that 
\beqn\label{Mtheta-G1L2}
\| M_\theta - \bar G_1 \|_{L^2}^2 %= {1 \over (2 \, \theta)^{N/2}} - {2 \over (\theta + \bar\theta)^{N/2}} 
\le C \, \rho^2 \, \big|\theta^{1/2} - \bar\theta_1^{1/2}\big|. 
\eeqn
Gathering (\ref{Ga-G1L1}) with  (\ref{Mtheta-G1L2}), (\ref{thetaa-b}) and (\ref{estimGeMtheta}) we deduce 
that for any $\eps > 0$ there is $C_\eps \in (0,\infty)$ such that 
\bean
\forall \,  \alpha \in (\alpha_0,1) \qquad  
\| G_\alpha - \bar G_1 \|_{L^1_{-N}}^{2+\eps} \le C_\eps \, \rho^{2+\eps} \, (1-\alpha), 
\eean
and (\ref{estimateonGe2}) follows by interpolation again. 
\qed
 
%%%%%%%%%%%%%%%%%%%%%%%%%%%%%%%%%%%%%%%%%%%%%%%

\section{Uniqueness and continuity of the path of self-similar profiles}
\setcounter{equation}{0}
\setcounter{theo}{0}

\subsection{The proof of uniqueness}

\begin{theo}\label{uniqueness1} 
There exists a constructive $\alpha_1 \in (0,1)$ 
such that the solution $G_\alpha$ of (\ref{profileeq3}) is unique for any $\alpha \in [\alpha_1,1]$. 
We denote by $\bar G_\alpha$ this unique self-similar profile. 
\end{theo}

That is an immediate consequence of the following result. 

\begin{prop}\label{uniqueness2} 
There is a constructive constant $\eta \in (0,1)$ such that 
$$ %\label{equniq2}
\quad \left.
  \begin{array}{l}
  G, \, H \in \GG_\alpha, \,\, \alpha \in (1-\eta,1) \vspace{0.3cm} \\
  \|G - \bar G_1 \|_{L^1_2} \le \eta, \,\,\,
    \|H - \bar G_1 \|_{L^1_2} \le \eta 
    \end{array}
  \right\} \quad\hbox{implies}\quad G = H. 
$$
\end{prop}

\smallskip\noindent{\sl Proof of Theorem~\ref{uniqueness1}. }
Let us assume that Proposition~\ref{uniqueness2} holds. 
Then Proposition~\ref{ConvGexplicit} implies that there is 
some explicit $\e \in (0,1)$ such that for $\alpha \in (1-\e,1]$ 
one has 
$$
\sup_{G_\alpha \in \GG_\alpha} \| G_\alpha - \bar G_1 \|_{L^1 _2} \le \eta
$$
where $\eta$ is defined in the statement of Proposition~\ref{uniqueness2}. Up to reducing 
$\eta$, it is always possible to take $\eta \le \e$, and the proof 
is completed by applying Proposition~\ref{uniqueness2}. 
\qed

\smallskip\noindent{\sl Proof of Proposition~\ref{uniqueness2}. }
Let us consider any exponential weight function  $m$ with $s \in (0,1)$, $a \in (0,+\infty)$, 
or with $s=1$ and $a \in (0,\infty)$ small enough. 
Let us also define $\OO= \Cc_{0,0,0}  \cap \Ll^1(m^{-1})$ the subvector space of  
of $\Ll^1(m^{-1})$ of functions with zero energy, $\psi = C \, (|v|^2 - N) \, M_{1,0,1}$ 
such that $\EE(\psi) = 1$, and $\Pi$ the following projection 
$$
\Pi : \Ll^1(m^{-1}) \to \OO, \qquad \Pi(g) = g - \EE(g) \, \psi.
$$
Finally, let us introduce $\Phi$ the following non-linear functional operator 
$$
\Phi : [0,1) \times ( W^{1,1}_1(m^{-1}) \cap \Cc_{\rho,0} ) \,\, \to \,\,  \R \times  \OO,
$$
and
$$
\Phi(1, \cdot ) : ( L^1_1(m^{-1}) \cap \Cc_{\rho,0} ) \,\, \to \,\,  \R \times  \OO,
$$ 
by setting 
\bean
\Phi(\alpha,g) = \left(  
(1+\alpha) \, D_\EE(g) - 2 \, \rho \, \EE(g), \,  \Pi \Bigl[Q_\alpha(g,g) 
        - \tau_\alpha \, \hbox{div}_v (v \, g) \Bigr]  \right). &&
\eean
It is straightforward that $\Phi(\alpha,G_\alpha) =0$ for any $\alpha \in [\alpha_0,1]$ and 
$G_\alpha \in \GG_\alpha$, and that the equation 
$$
\Phi(1,g) = (0,0)
$$
has a unique solution, given by $g = \bar G_1 = M_{\rho,0,\bar\theta_1}$ defined in (\ref{Maxlim}), (\ref{tempSS}). 

The function $\Phi$ is linear and quadratic in its second argument by inspection, and easy computations 
yield the following formal differential according to the second argument at the point $(1,\bar G_1)$: 
\beqn\label{D2Phi}
  D_2 \Phi(1,\bar G_1) \, h = A \, h := \Bigg(  4 \,  \tilde D_\EE( \bar G_1, h) 
    - 2 \rho \, \EE(h),  \ 2 \, \tilde Q_1(\bar G_1,h) \Bigg) 
\eeqn
where $\tilde Q_\alpha$ is defined in (\ref{Qinelsym}) and 
$$
\tilde D_\EE (g,h) := {b_1} \, \int\!\!\int_{\R^N \times \R^N} g \, h_*\, |u|^3 \, dv \, dv_*.
$$
Notice that we can remove the projection on the last argument in (\ref{D2Phi}) since the elastic collision 
operator always has zero energy. 

\smallskip
Then we have the
\begin{lem}\label{lemDPhiM} The linear functional 
\bean %\label{lemDPhiM1} 
A \ : \ \Ll^1_1(m^{-1}) & \to &  \R \times  \OO \\
h & \mapsto & A \, h = D_2 \Phi(1,\bar G_1) \, h \nonumber
\eean
is invertible: it is bijective with $A^{-1}$ bounded with explicit estimate. 
\end{lem}

\medskip\noindent
{\sl Proof of Lemma~\ref{lemDPhiM}. } 
Since the spectrum of the linear operator $\LL_1$ defined on $L^1(m^{-1})$ (with domain $L^1_1(m^{-1})$) includes $0$ as a discrete  eigenvalue associated with the eigenspace 
$\hbox{Ker} \LL_1 = \hbox{Span}\{ \bar G_1,v_1 \, \bar G_1,\dots,v_N \, \bar G_1,|v|^2 \, \bar G_1\}$ 
by~\cite[Theorem~1.3]{GM:04} and since moreover $\OO \cap \hbox{Ker} \LL_1 = \{ 0 \}$, 
we deduce that it is invertible from $\OO \cap \Ll^1_1(m^{-1})$ onto $\OO$. 
Moreover the work~\cite[Section~4]{GM:04} provides explicit estimates 
on the norm of its inverse. We deduce immediately that $\LL_1^{-1}$ maps $\OO$ onto itself with explicit bound.

For any $h \in \Ll^1(m^{-1})$, we decompose 
$$
h = h_1 \, \phi_1  + h^\bot, \quad \mbox{with} \quad h_1:= {\EE(h) \over \EE(\phi_1)} \in \R, \,\,h^\bot \in \OO,
$$
where we recall that $\phi_1$ is defined in (\ref{phiaTOphi1}). 
Then, using the characterization (\ref{profileeq3}) of $\bar G_1$, 
  \bean 
  && A \, h = \Bigg( 
     b_1 \, \int_{\R^N \times \R^N} \bar G_1(v) \, h^\bot (v_*) \, |u|^3 \, dv \, dv_* \\
  && \hspace{2cm} + h_1 \, \left[ b_1 \, \int_{\R^N \times \R^N} |v|^2 
             \, \bar G_1 \, \bar G_{1*} \, |u|^3 \, dv \, dv_* 
                        - 2 \rho \, \int_{\R^N}  \bar G_1 \, |v|^4 \, dv \right], \ \LL_1( h^\bot) \Bigg). 
  \eean
The claimed invertibility follows from the fact that $C^* = 2 \, N \, \rho^2 \, \bar\theta_1^2 \not= 0$. Indeed, from (\ref{Mv4}) and (\ref{MMv2u3}) there holds
\bean
C^* &:= & b_1 \, \int_{\R^N \times \R^N} |v|^2 
             \, \bar G_1 \, \bar G_{1*} \, |u|^3 \, dv \, dv_* 
                        - 2 \rho \, \int_{\R^N}  \bar G_1 \, |v|^4 \, dv \\
&= & \rho^2 \, \bar\theta_1^2 \, \Bigg[ 
    b_1 \, \bar\theta_1^{1/2} \, \int_{\R^N \times \R^N} M_{1,0,1} \, 
          (M_{1,0,1})_* \, |v|^2 \, |u|^3 \, dv \, dv_* 
    -  2 \, \int_{\R^N} M_{1,0,1} \, |v|^4 \, dv \Bigg] \\
&=&  \rho^2 \, \bar\theta_1^2 \, \Bigg[ 
         b_1 \, \bar\theta_1^{1/2} \, \sqrt{2} \, (2N+3) \, {\bf m}_{3/2}( M_{1,0,1}) -  2 \, N (N+2)  \Bigg],
\eean
and we conclude thanks to formula~(\ref{tempSS}).
\qed

\medskip\noindent
Let us come back to the proof of Proposition~\ref{uniqueness2}. We write 
\bear
G_\alpha - H_\alpha &=& A^{-1} \, \Big[ A \, G_\alpha - \Phi(\alpha,G_\alpha) 
      + \Phi(\alpha,H_\alpha) - A \, H_\alpha  \Big]  \nonumber \\
\label{DA-1I}
&=& A^{-1} \, (I_1, I_2) 
\eear
with (recall that the bilinear operators $\tilde D_\EE$ and $\tilde Q_\alpha$ are symmetric)
\[
\left\{
\begin{array}{l}
%I_1 := 0 \vspace{0.3cm} \\
I_1 :=  4 \, \tilde D_\EE(\bar G_1, G_\alpha - H_\alpha) 
- (1+\alpha) \, D (G_\alpha) + (1+\alpha) \, D (H_\alpha) \vspace{0.3cm} \\
I_2 := \Pi \, I_{2,1}+  \Pi \, I_{2,2}
\end{array}
\right.
\]
and
\[
\left\{
\begin{array}{l}
I_{2,1} :=  2 \, \tilde Q_1(\bar G_1,G_\alpha-H_\alpha) 
                - Q_\alpha(G_\alpha,G_\alpha) + Q_\alpha(H_\alpha,H_\alpha)  \vspace{0.3cm} \\
I_{2,2} :=  \rho \, (1-\alpha)  \, \nabla_v \cdot \big(v \, (H_\alpha - G_\alpha) \big). 
\end{array}
\right.
\]
On the one hand, 
$$
I_1 = 2 \, D\big(2 \bar G_1 - (G_\alpha + H_\alpha), G_\alpha - H_\alpha\big) 
+ (1-\alpha)  \, D (G_\alpha + H_\alpha, G_\alpha - H_\alpha)
$$
so that 
\bear \nonumber
|I_1| &\le& C_3 \, 
\Big( \|\bar G_1 - G_\alpha\|_{L^1 _3} + \|\bar G_1 - H_\alpha \| _{L^1 _3} \\ \nonumber
&& \hspace{4cm} + (1-\alpha) \, 
 \|G_\alpha \|_{L^1 _3} + (1-\alpha) \, \|H_\alpha \|_{L^1 _3} \Big) \, \|G_\alpha - H_\alpha \|_{L^1 _3}
\nonumber \\
\label{I2} &\le& \eta_1(\alpha) \, \|G_\alpha - H_\alpha \|_{L^1_1(m^{-1})}
\eear
with $\eta_1(\alpha) \to 0$ when $\alpha \to1$ (with explicit rate, for instance $\eta_1(\alpha) = C_1 \, (1-\alpha)^{1/3}$) 
because of Propositions~\ref{estimatesonGe} and~\ref{ConvGexplicit}.

On the other hand, 
\bean
I_{2,1} &=&  Q_1(\bar G_1, G_\alpha - H_\alpha) -  Q_\alpha (\bar G_1, G_\alpha - H_\alpha) 
+ Q_1(G_\alpha - H_\alpha, \bar G_1) - Q_\alpha (G_\alpha - H_\alpha,\bar G_1) \\
&& \hspace{3cm}
+ Q_\alpha(\bar G_1 - G_\alpha, G_\alpha - H_\alpha) + Q_\alpha (G_\alpha - H_\alpha,\bar G_1 - H_\alpha). 
\eean  
From Proposition~\ref{ContQe2} there holds 
  \[ \| Q_1(\bar G_1, G_\alpha - H_\alpha) -  Q_\alpha(\bar G_1, G_\alpha - H_\alpha) \|_{L^1(m^{-1})} \le 
           \eps(\alpha) \, \| G_\alpha - H_\alpha \|_{L^1 _1(m^{-1})} \]
  \[ \| Q_1(G_\alpha - H_\alpha, \bar G_1) - Q_\alpha (G_\alpha - H_\alpha,\bar G_1) \|_{L^1(m^{-1})} \le 
           \eps(\alpha) \, \| G_\alpha - H_\alpha \|_{L^1 _1(m^{-1})} \]
with $\eps(\alpha) \to 0$ as $\alpha \to 1$  (with again explicit rate, for instance 
$\eps(\alpha) = C_1' \, (1-\alpha)^{1/12}$ if $s=1/2$ in the formula of $m$).
From elementary estimates in $L^1(m^{-1})$ we have 
\bean 
&& \|  Q_\alpha(\bar G_1 - G_\alpha, G_\alpha - H_\alpha) 
+ Q_\alpha (G_\alpha - H_\alpha,\bar G_1 - H_\alpha) \|_{L^1(m^{-1})} \\ 
&& \hspace{3cm} 
\le C_4 \, \left( \| G_\alpha - \bar G_1 \|_{L^1 _1(m^{-1})} + \| H_\alpha -\bar G_1 \|_{L^1 _1(m^{-1})} \right) 
           \, \| G_\alpha - H_\alpha \|_{L^1 _1(m^{-1})}.
\eean
Together with Propositions \ref{ConvGexplicit} we thus obtain
\beqn \label{I31}
\|I_{2,1} \|_{L^1(m^{-1})} \le \eta_2(\alpha) \, \|G_\alpha - H_\alpha \|_{L^1_1(m^{-1})}
\eeqn 
for some $\eta_2(\alpha) \to 0$ as $\alpha \to 1$. Here we can take for instance 
(when $s=1/2$ in the formula of $m$)
$\eta_2(\alpha) = C_2 \, (1-\alpha)^{1/12}$ 
for some $C_2 \in (0,\infty)$ by picking a suitable $\varepsilon$ and interpolating. 

Finally from Proposition~\ref{errorH1} there holds 
\beqn \label{I32}
\|I_{2,2} \|_{L^1(m^{-1})} \le C_5 \, (1-\alpha) \, \|G_\alpha - H_\alpha \|_{L^1_1(m^{-1})}.
\eeqn 
Gathering~(\ref{I2}), (\ref{I31}) and (\ref{I32}) we obtain from (\ref{DA-1I}) and Lemma~\ref{lemDPhiM}
$$
\| G_\alpha - H_\alpha \|_{L^1_1(m^{-1})} \le \eta(\alpha) \, \|A^{-1}\| \, \|G_\alpha - H_\alpha \|_{L^1_1(m^{-1})}
$$
for some function $\eta$ such that  $\eta(\alpha) \to 0$ as $\alpha \to 1$ (with explicit rate). 
Hence choosing $\alpha_1$ close enough to $1$ 
we have $\eta(\alpha) \, \|A^{-1}\|  \le 1/2$ for any $\alpha \in [\alpha_1,1)$. 
This implies $G_\alpha = H_\alpha$ and concludes the proof. \qed 

%%%%%%%%%%%%%%%%%%%%%%%%%%%%%%%%%%%%%%%%%%%%%%%

\subsection{Differentiability of the map $\alpha \mapsto \bar G_\alpha$ at $\alpha=1$}
\label{derivee=1Ge}

\begin{lem}\label{Contderivee=1Ge}  
The map $[\alpha_1,1] \to L^1(m^{-1})$, $\alpha \mapsto \bar G_\alpha$ is continuous 
on $[\alpha_1,1]$ and differentiable at $\alpha= 1$. More precisely, 
there exists $\bar G'_1 \in L^1(m^{-1})$ and for any $\eta \in (1,2)$ 
there exists a constructive $C_\eta \in (0,\infty)$ such that 
\beqn\label{devGalpha1}
\| \bar G_\alpha - \bar G_1 - (1-\alpha) \, \bar G'_1 \|_{L^1(m^{-1})} \le C_\eta \, (1-\alpha)^\eta
\qquad \forall \, \alpha \in (\alpha_0,1).
\eeqn
\end{lem}

%\begin{rem} {\bf PLUS TARD}
%Using the same arguments as in the proof below, it could be proved that the map 
%$[\alpha_1,1] \to L^1(m^{-1})$, $\alpha \mapsto \bar G_\alpha$ is lipschitzian (and $C^1$??) 
%on the whole interval $[\alpha_1,1]$, but this point would not be constructive.
%\end{rem}
%% CM : Voir s'il y en a besoin plus tard.

\medskip\noindent
{\sl Proof of Lemma~\ref{Contderivee=1Ge}.} We split the proof into four steps. 

\medskip
\noindent
{\em Step 1.} For the continuity we use a classical stability argument. Let us consider a sequence $(\alpha_n)_{n \ge 0}$ such that $\alpha_n \in [\alpha_1,1]$ and $\alpha_n \to \alpha$. From the uniform bound (\ref{estimGeunif}), we may extract  a subsequence $(\bar G_{\alpha_{n'}})$ which strongly converges in $L^1(m^{-1})$ to a function $G_\alpha$. Passing to the limit in the equations (\ref{profileeq3}) associated to the normal restitution coefficient $\alpha_n$ and written for $G_{\alpha_{n'}}$,  we deduce that $G_\alpha$ satisfies (\ref{profileeq3}) associated to the normal restitution coefficient $\alpha$.  
From the uniqueness of the solution proved in Theorem~\ref{uniqueness1}, there holds $G_\alpha = \bar G_\alpha$ and thus the whole sequence $\bar G_{\alpha_n}$ converges to $\bar G_\alpha$. 

\medskip\noindent
{\sl Step 2.} 
We next prove that there exists an explicit constant $C$ such that 
$$
\forall \, \alpha \in [\alpha_1,1] \qquad \| \bar G_\alpha - \bar G_1 \|_{L^1(m^{-1})} \le C \, (1-\alpha).
$$
We write 
\bear
\bar G_\alpha - \bar G_1 &=& A^{-1} \, [ A \, \bar G_\alpha - \Phi(\alpha,\bar G_\alpha) + \Phi(1,\bar G_1) - A \, \bar G_1] \nonumber \\
\label{DA-1J}&=& A^{-1} \, (J_1, J_2) 
\eear
with 
\[
\left\{
\begin{array}{l}
%J_1 := 0 \vspace{0.3cm} \\
J_1 :=  4 \, \tilde D_\EE(\bar G_1, \bar G_\alpha - \bar G_1) 
+ 2 \,  \tilde D_\EE(\bar G_1, \bar G_1) - (1+\alpha) \, \tilde D_\EE(\bar G_\alpha, \bar G_\alpha) \vspace{0.3cm} \\
J_2 := \Pi \, J_{2,1}+  \Pi \, J_{2,2}
\end{array}
\right.
\]
and
\[
\left\{
\begin{array}{l}
J_{2,1} :=  Q_1(\bar G_1,\bar G_\alpha) +Q_1(\bar G_\alpha,\bar G_1) 
                - Q_\alpha(\bar G_\alpha,\bar G_\alpha)  \vspace{0.3cm} \\
J_{2,2} :=  \rho \, (1-\alpha)  \, \nabla_v \cdot \big(v \, (\bar G_\alpha)\big). 
\end{array}
\right.
\]
On the one hand, 
$$
J_1 = - 2 \, \tilde D_\EE( \bar G_1 - \bar G_\alpha,  \bar G_1 - \bar G_\alpha) + (1-\alpha) \, D (\bar G_\alpha,\bar G_\alpha)
$$
so that 
\bean
|J_1| \le C \, \|\bar G_1 - \bar G_\alpha\|_{L^1 _3}^2 + C \, (1-\alpha).
\eean

\medskip\noindent
On the other hand, 
\bean
J_{2,1} &=&  - Q_1(\bar G_1-\bar G_\alpha, \bar G_1 - \bar G_\alpha) +  Q_1(\bar G_\alpha, \bar G_\alpha) - Q_\alpha(\bar G_\alpha, \bar G_\alpha).
\eean  
Hence using Propositions~\ref{errorH1}, \ref{ContQe1}, and the bound~(\ref{estimGeunif}), we deduce 
$$
| J_{2,1} | \le \left\| \bar G_\alpha - \bar G_1 \right\|_{L^1 _2} ^2  + C \, (1-\alpha)
$$ 
and we also have  straightforwardly $J_{2,2} = \OO(1-\alpha)$.
Gathering all these estimates, we thus obtain from (\ref{DA-1J}) 
$$
\big\| \bar G_\alpha - \bar G_1 \big\|_{L^1_1(m^{-1})} \le \|A^{-1}\| \, \Big[ \big\|\bar G_\alpha - \bar G_1 \big\|_{L^1_1(m^{-1})}^2 
+ C \, (1-\alpha) \Big].
$$
Using then the explicit result of quantification of the elastic limit in Proposition~\ref{ConvGexplicit}, we have that for some $\alpha_2 \in [\alpha_1,1)$ close enough to $1$: 
$$ 
\forall \, \alpha \in [\alpha_2,1] 
\qquad \|A^{-1}\| \,   \big\|\bar G_\alpha - \bar G_1 \big\|_{L^1_1(m^{-1})} < \frac12
$$ 
and thus we get 
$$
\forall \, \alpha \in [\alpha_2,1], \quad \| \bar G_\alpha - \bar G_1\|_{L^1_1(m^{-1})} \le 2 \, C \,  \|A^{-1}\| \, (1-\alpha)
$$
which implies the claimed estimate. 

\medskip\noindent
{\sl Step 3.} In order to prove the differentiability we must slightly improve the estimate established in the preceding step. 
On the one hand we exhibit what should be the derivative of $\bar G_\alpha$ at $\alpha=1$, and denote it by $R$. 
Formally differentiating equation~(\ref{profileeq3}) at $\alpha = 1$ we have
$$
Q'_1(\bar G_1,\bar G_1) + 2 \tilde Q_1(R,\bar G_1) + \rho \, \nabla_v \cdot \big(v \, \bar G_1\big) = 0.
$$
%where $K$ is defined in the weak sense in the following way
%\bean
%\langle K,\psi \rangle 
%&=& {d \over d\alpha} \Big|_{\alpha=1} \left(  { 1 \over 2}  \int \!\! \int \!\! \int_{\R^N \times \R^N \times S^{N-1}}
%b \, |u| \, \bar G_1 \, \bar G_{1*} \, \Bigl[ \psi(v'_\alpha) + \psi(v'_{*\alpha}) \Bigr] \, dv \, dv_* \, d\sigma \right) \\
%&=& { 1 \over 8} \, \int \!\! \int \!\! \int_{\R^N \times \R^N \times S^{N-1}}
%b \, |u| \, \bar G_1 \, \bar G_{1*}  \, (|u| \, \sigma - u) \, \Bigl[  (\nabla \psi) (v') - 
%(\nabla \psi) (v'_*)  \Bigr] \, dv \, dv_*\, d\sigma,
%\eean
%where $v'_\alpha$ and $v'_{*\alpha}$ denote the post-collisional velocities (\ref{vprimvprim*}) 
%associated to the normal restitution coefficient $\alpha$, while $v'$ and $v'_{*}$ denote 
%the elastic post-collisional velocities ({\em i.e.}, defined thanks to (\ref{vprimvprim*}) 
%with $\alpha=1$). Performing the prepost-collisional
%change of variable $(v,v_*,\sigma) \to (v',v'_*,k)$ with $k = \hat u$ (see for instance 
%\cite[Chapter 1, Section 4.4]{CVhand}), we see that $K$ is indeed the smooth function 
%$$
%K (v) = \nabla_v \cdot \left(  
%{ 1 \over 4} \, \int \!\! \int_{\R^N \times S^{N-1}}
%b \, |u| \, \bar G_1 \, \bar G_{1*}  \, (|u| \, \sigma - u) \, dv_* \, d\sigma \right) = Q'_1(\bar G_1,\bar G_1). 
%$$
On the other hand, we may compute 
\bear\nonumber
\langle Q'_\alpha(\bar G_1,\bar G_1), |.|^2 \rangle 
&=& { 1 \over 4} \, \int \!\! \int \!\! \int_{\R^N \times \R^N \times S^{N-1}}
b \, |u| \, \bar G_1 \, \bar G_{1*}  \, (|u| \, \sigma - u) \cdot \Bigl(  |u| \, \sigma  \Bigr) \, dv \, dv_* \, d\sigma \\\label{Kv2}
&=& 2 \, D_\EE(\bar G_1). 
\eear
Next, diving the  equation (\ref{EnergyDEGe}) on the energy of $G_\alpha$ by $(1-\alpha)$ and formally differentiating the resulting expression we get 
$$
2  \, \rho \, \EE(R) - \tilde D_\EE(\bar G_1,\bar G_1) -  4 \,  \tilde D_\EE(R,\bar G_1) = 0.
$$
We now rigorously define $R$ in the following way
$$
\bar G'_1 = R := A^{-1} \Big( - \tilde D_\EE(\bar G_1, \bar G_1) , - F \Big), \qquad 
F := Q'_\alpha(\bar G_1,\bar G_1) + \rho \, \nabla_v \cdot \big(v \, \bar G_1\big).
$$
Note that $R$ is well-defined since $\EE(F) = 0$ because of (\ref{Kv2}) and the definition of $\bar G_1$. 

\medskip\noindent 
{\sl Step 4.} 
We finally come back to the Step 2 and we shall construct a Taylor expansion of order $1$. 
We want to estimate 
$$
\bar G_\alpha - \bar G_1 + (\alpha-1) \, \bar G'_1 = 
A^{-1} \, \Big(J_1 - (\alpha-1) \, \tilde D_\EE(\bar G_1, \bar G_1), J_2-(1-\alpha) \, F \Big). 
$$
On the one hand
$$
J_1 - (\alpha-1) \, \tilde D_\EE(\bar G_1, \bar G_1) =  
-2 \, \tilde D_\EE( \bar G_1 - \bar G_\alpha,\bar G_1 - \bar G_\alpha) 
+ (1-\alpha) \, \Big(D (\bar G_\alpha,\bar G_\alpha) - \tilde D_\EE(\bar G_1,\bar G_1)\Big)
$$
so that we obtain straightforwardly
$$
\big|J_1 - (\alpha-1) \, \tilde D_\EE(\bar G_1,\bar G_1)\big| \le C \, (1-\alpha)^2. 
$$
On the other hand, 
$$
J_2 - (1- \alpha) \, F := \Pi \, J_{2,1}+  \Pi \, J_{2,2}
$$
with 
\bean
J_{2,1} &=& 
- Q_1(\bar G_1 - \bar G_\alpha,\bar G_1 - \bar G_\alpha) +  Q_1(\bar G_\alpha - \bar G_1,\bar G_\alpha) -
Q_\alpha(\bar G_\alpha - \bar G_1,\bar G_\alpha) \\
&&+ Q_1(\bar G_1, \bar G_\alpha -\bar G_1) -  Q_\alpha(\bar G_1, \bar G_\alpha -\bar G_1)  
+ (1-\alpha)  \, \nabla_v \cdot \big(v \, (\bar G_\alpha-\bar G_1)\big)
\eean
and 
$$
J_{2,2} =  Q_1(\bar G_1, \bar G_1) - Q_\alpha(\bar G_1, \bar G_1) - (1-\alpha) \, K.
$$
It is clear from Propositions~\ref{ContQe1}, the bound of Step~2, and some interpolation with 
the uniform bounds~(\ref{estimGeunif}), that 
$$
\| J_{2,1} \|_{L^1(m^{-1})}, \, \| J_{2,2} \|_{L^1(m^{-1})} \le C_k \, (1-\alpha)^k
$$ 
for any $k \in (1,2)$.
\qed

%%%%%%%%%%%%%%%%%%%%%%%%%%%%%%%%%%%%%%%%%%%%%%%%%%

\section{Study of the spectrum and semigroup of the linearized problem}
\setcounter{equation}{0}
\setcounter{theo}{0}

%%%%%%%%%%%%%%%%%%%%%%%%%%%%%%%%%%%%%%%%%%%%%%%%%%

In this section we shall obtain the geometry of the spectrum of the linearized rescaled inelastic collision operator 
for a small inelasticity, as well as  estimates on its resolvent and on the associated linear semigroup. 
This is based on the properties of the elastic linearized operator and some perturbation arguments again. 
In order to do so, one needs some common functional ``ground" for the the linearized operators in the limit
of vanishing inelasticity. This common functional setting is given by  the study~\cite{GM:04} in which the spectral study of the elastic linearized operator is made in $L^1$ spaces with exponential weights $e^{a \, |v|^s}$, $a \in (0,+\infty)$, $s \in (0,1)$. 

We thus consider the operator 
  \[ g \mapsto Q_\alpha(g,g) - \tau_\alpha \, \nabla_v \cdot \big(v \, g\big) \]
and some fluctuations $h$ around the self-similar profile $\bar G_\alpha$:  
that means $g = \bar G_\alpha  +  h$ with $h \in L^1(m^{-1})$ where $m$ is a fixed smooth exponential weight function, 
as defined in (\ref{defdem}).  
The corresponding linearized unbounded operator $\LL_\alpha$ acting on $L^1(m^{-1})$  
with domain $\hbox{dom}(\LL_\alpha) = W^{1,1}_1(m^{-1})$ if $\alpha \not= 1$ and $\hbox{dom}(\LL_1) = L^1_1(m^{-1})$, is defined in (\ref{defLa}) (it is straightforward to check that it is closed in this space).
Since the equation in self-similar variables preserves mass and the zero momentum, the correct spectral 
study of $\LL_\alpha$ requires to restrict this operator to zero mean and centered distributions (which are preserved as well), that means to work in $\Ll^1(m^{-1})$. When restricted to this space, the operator $\LL_\alpha$ is denoted 
by $\hat \LL_\alpha$. We denote by $R(\hat \LL_\alpha)$ the resolvent 
set of $\hat \LL_\alpha$, and by $\Rr_\alpha(\xi) = (\hat \LL_\alpha - \xi)^{-1}$ its 
resolvent operator for any $\xi \in R(\hat \LL_\alpha)$. 

Let us recall that for the linearized elastic hard spheres Boltzmann equation the spectrum and the 
asymptotic stability have been studied by many authors since the pioneering works
by Hilbert~\cite{Hilb:EB:12}, Carleman~\cite{Carl:57} and Grad \cite{Grad63}, 
and we refer for instance to \cite{GM:04} for 
more references. The result established for $\LL_1$ (and translated straightforwardly 
to $\hat \LL_1$) in~\cite{GM:04} is the following:

\begin{theo}\label{thL1} 
\begin{itemize}
\item[(i)] There exists a decreasing sequence of real discrete eigenvalues  $(\mu_n)_{n \ge 1}$  
(that is: eigenvalues isolated and with finite multiplicity) of $\hat \LL_1$, with ``energy" eigenvalue $\mu_1 = 0$ of multiplicity $1$ and ``energy" eigenvector $\phi_1$ (defined in (\ref{phiaTOphi1})), 
$\mu_2 < 0$ and $\lim \mu_n = \mu_\infty \in (- \infty,0)$ such that the spectrum 
$\Sigma(\hat\LL_1)$ of $\hat \LL_1$ in $\Ll^1(m^{-1})$ writes
$$
\Sigma(\hat\LL_1) = (-\infty,\mu_\infty]  \cup \{ \mu_n \}_{n \in \N}.
$$
In particular,  $\hat \LL_1$ is onto from $\OO \cap \Ll^1_1(m^{-1})$ onto $\OO$. 
\item[(ii)] The resolvent  $\Rr_1(\xi)$ has a sectorial property for the spectrum substracted from the 
``energy" eigenvalue, namely there is a constructive $\mu_2 < \lambda <0$ such that  
$$
\forall \, \xi \in \Aa, \quad \| \Rr_1 (\xi) \|_{\Ll^1(m^{-1})} \le a + \frac{b}{|\xi + \lambda|}, 
$$
with 
$$
\Aa = \left\{ \xi \in \mathbb{C}, \quad \mbox{{\em arg}}(\xi +\lambda) 
\in \left[-\frac{3 \pi}{4}, \frac{3 \pi}{4} \right] \ \mbox{ and } \ \Re e \, \xi \le \frac{\lambda}2 \right\}.
$$

\item[(iii)] The linear semigroup $S_1(t)$ associated to $\hat \LL_1$ in $\Ll^1(m^{-1})$ writes
$$
\forall \, t \ge 0 \qquad S_1(t) = \Pi_1+  R_1(t),
$$
where $\Pi_1$ is the projection on the eigenspace associated to $\mu_1$ and $R_1(t)$ is a semigroup which satisfies 
$$
\forall \, t \ge 0 \qquad  \|R_1(t) \|_{\Ll^1(m^{-1})} \le C \, e^{ \mu_2 \, t }
$$
with explicit constant $C$. 
\end{itemize}
\end{theo}

\medskip
The main result proved in this section is a perturbation result which extends 
Theorem~\ref{thL1} in the following way. Let us define for any $x \in \R$ 
the half-plane $\Delta_x$ by
$$
\Delta_x = \{ \xi \in \CC, \,\, \Re e \, \xi \ge x \}.
$$

\begin{theo}\label{thLalpha} 
Let us fix $\bar \mu \in (\mu_2,0)$, $k,q \in \N$ and $m$ a smooth weight  exponential function with $s \in (0,1)$. 
Then there exists $\alpha_2 \in (\alpha_1,1)$ such that for any 
$\alpha \in [\alpha_2,1]$ the following holds:

\begin{itemize}
\item[(i)] The spectrum $\Sigma(\hat  \LL_\alpha)$ of $\hat \LL_\alpha$ 
in $\Ww^{k,1} _q(m^{-1})$ writes
$$
\Sigma(\hat \LL_\alpha) = E_\alpha \cup \{ \mu_\alpha \},% \cup \{ 0 \},
\qquad E_\alpha \subset \Delta_{\bar \mu} ^c,
$$
where $\mu_\alpha$ is a $1$-dimensional real eigenvalue which does not depend 
on the choice of the space $\Ww^{k,1} _q(m^{-1})$  and  satisfies (\ref{expanmua}).
 
\item[(ii)] The resolvent  $\Rr_\alpha(\xi)$  in $\Ww^{k,1}_q(m^{-1})$ is holomorphic on a 
neighborhood of $\Delta_{\bar \mu} \backslash \{ \mu_\alpha \}$ and 
there are explicit constants $C_1,C_2$ such that 
$$
\sup_{z \in \CC, \,\, \Re e \, z = \bar \mu} 
\| \Rr_\alpha(z)|\|_{\Ww^{k,1} _q(m^{-1}) \to \Ww^{k,1} _q(m^{-1})} \le C_1
$$
and 
$$
\| \Rr_\alpha(\bar \mu + i s) \|_{\Ww^{k+1,1} _{q+1}(m^{-1}) \to \Ww^{k,1} _q(m^{-1})} \le 
\frac{C_2}{1+|s|}.
$$

\item[(iii)] The linear semigroup $S_\alpha(t)$ associated to $\hat \LL_\alpha$ in 
$\Ww^{k,1} _q (m^{-1})$ writes
$$
S_\alpha(t) = e^{\mu_\alpha \, t}  \, \Pi_\alpha + R_\alpha(t),
$$
where $\Pi_\alpha$ is the projection on the ($1$-dimensional) eigenspace associated to $\mu_\alpha$ 
and where $R_\alpha(t)$ is a semigroup which satisfies 
\beqn\label{estimResta}
\|R_\alpha(t) \|_{\Ww^{k+2,1} _{q+2}(m^{-1}) \to \Ww^{k,1} _q(m^{-1})}
 \le C_{k} \, e^{\bar \mu \, t } 
\eeqn
with explicit bounds. 
\end{itemize}
\end{theo}

\begin{rem} 
Note that we do not claim that the resolvent $\Rr_\alpha$ is  sectorial for $\alpha < 1$ 
and it is likely that indeed it is not (because of the contribution of the drift term). 
Moreover, it is not clear how to make the spectral study in the Hilbert setting $L^2(m^{-1})$
with convenient  weight function $m$. 
In particular, we are not able to prove Proposition~\ref{ContQe2} in an $L^2$ framework.
In such a situation the spectral study and 
the obtaining of constructive rate of decay on the semigroup become tricky.  Let us 
emphasize also that (as most of the results established in this paper) this result 
is not an easy consequence of perturbation theory of unbounded operator since the elastic 
limit $\alpha \to 1$ is strongly bad-behaved (for instance neither the relative bound nor 
the operator gap of~\cite{Kato} go to $0$) because of the anti-drift term. 
\end{rem}

%%%%%%%%%%%%%%%%%%%%%%%%%%%%%%%%%%%%%%%%%%%%%%%%%%

\subsection{Recalls and improvments of technical tools from~\cite{GM:04}}\label{subsec:impGM}

 \begin{prop}\label{thL1Wk}
 In the statement  of Theorem~\ref{thL1} one can replace  
 everywhere $L^1(m^{-1})$ by $W^{k,1}_q(m^{-1})$, $k,q \in \N$.
 \end{prop}

%Thus from these two last propositions and following exactly the same 
%proof as in~\cite{GM:04} one can 

Let us first recall the key decomposition of $\hat \LL_1$ in~\cite[Section~2]{GM:04} 
(re-written within the notation of this paper):

Let ${\bf 1}_E$ denote the usual indicator function of the set $E$, 
let $\Theta: \R \to \R_+$ be an even $C^\infty$ function with 
mass $1$ and support included in $[-1,1]$ and 
$\tilde{\Theta}: \R^N \to \R_+$ a radial $C^\infty$ function 
with mass $1$ and support included in $B(0,1)$. We define the 
following mollification functions ($\epsilon >0$):
  \[ \left\{
     \begin{array}{ll}
     \Theta_\epsilon (x) = \epsilon^{-1} \, \Theta(\epsilon^{-1} x), \ \ \ (x \in \R) \vspace{0.2cm} \\
     \tilde{\Theta}_\epsilon (x) = \epsilon^{-N} \, \tilde{\Theta}(\epsilon^{-1} x), \ \ \ (x \in \R^N).
     \end{array}
     \right. \]
Then we consider the decompositions 
$$
 \LL_1(g) = \LL^c_1 (g) - \LL^\nu (g) \ \ \ \mbox{ with } \ \ \ 
 \LL^\nu  (g) := \nu \, g
 $$
where $\LL^c _1$ splits between a ``gain'' part $\LL^+ _1$ (denoted so because it 
corresponds to the linearization of $Q^+$) and a convolution part $\LL^*$ (not 
depending on $\alpha$) as 
 $$
 \LL^c _1(g) = \LL^+ _1(g) - \LL^*(g) 
 \ \ \ \mbox{ with } \ \ \ 
 \LL^* (g) :=  M \, \big[ g * \Phi \big],
 $$
(we do not write the subscript $1$ when there is no dependency on $\alpha$). 
Then for any $\delta \in (0,1)$ we set 
  $$ %\label{eq:auxdefronL^+delta}
  \LL^+ _{1,\delta} (g) = {\cal I}_\delta(v) \, 
  \int_{\R^N \times \Sph^{N-1}} \Phi(|v-v_*|) \, b_\delta (\cos \theta) \, 
  \big[ g' M' _* + M' g' _* \big] \, dv_* \, d\sigma,  
  $$ 
where  
  $$ %\label{eq:defI}
  {\cal I}_\delta = \tilde{\Theta}_\delta * {\bf 1}_{\{|\cdot| \le \delta^{-1}\}}, 
  $$ 
and 
  $$ %\label{eq:defbdelta}
  b_\delta(z) = \left( \Theta_{\delta^2} * {\bf 1}_{\{-1+2\delta^2 \le z \le 1 - 2\delta^2 \}} 
                \right) \, b (z). 
  $$
This approximation induces $\LL_{1,\delta} = \LL^+ _{1,\delta} - \LL^* - \LL^\nu$.
Then the key result is that this approximation converges (in the norm of the graph) 
to the original linearized operator $\LL_1$ as $\delta \to 0$, first in the small classical 
linearization space $L^2(\bar G_1^{-1})$ (this technical result was in fact 
mostly already included in Grad's results~\cite{Grad63}), 
and second most importantly in the larger space 
$L^1(m^{-1})$. On the basis of this approximation result the spectrum 
is then proved to be the same in both functional spaces, and then the norm 
of the resolvents within these two functional spaces are related by an explicit control. 

\medskip
Hence the keys elements of the proof which are to be extended are, on the one hand, 
the approximation argument (which has to be extended from an $L^1(m^{-1})$ setting 
to an $W^{k,1}_q(m^{-1})$ setting), and, on the other hand the explicit control on the 
resolvent in the space $L^2(M^{-1})$ provided by the self-adjointness structure of the 
collision operator in this space and the explicit estimates on the spectral gap (see~\cite{BM:04}), 
which has to be extended to an $H^k(M^{-1})$ setting. 
Then the rest of the proof of~\cite{GM:04}
would extend as well (up to minor technical modifications) to $W^{k,p}(m^{-1})$. 
\medskip

Therefore for the first point let us prove the 
 \begin{prop}\label{prop:cvgronLp}
For any $k,q \in \N$ and $g \in W^{k,1}_{q+1}(m^{-1})$, we have 
   \[ \left\| \left( \LL^+ _1 - \LL^+ _{1,\delta} \right) (g) \right\|_{W^{k,1}_q(m^{-1})} \le 
       \eps (\delta) \, \|g\|_{W^{k,1} _{q+1} (m^{-1})} \]
 where $\eps (\delta) >0$ is an explicit constant 
 going to $0$ as $\delta$ goes to $0$. 
 \end{prop}
  
\smallskip\noindent
{\sl Proof of Proposition~\ref{prop:cvgronLp}. } 
The case $k=q=0$ is provided by Proposition~\ref{ContQe2}. 
Then higher-order derivatives follows by differentiation, and the incoporation 
of a polynomial weight is trivial. 
\qed

Concerning the second point let us prove the 

\begin{prop}\label{prop:resolHk} 
The spectrum $\Sigma(\LL_1)$ of $\LL_1$ in $L^2(M^{-1})$ is the same in any $H^k(M^{-1})$, 
$k \in \N$. Moreover the control on the resolvent, which was (self-adjoint operator) 
$$
\| \Rr_1(\xi) \|_{L^2(M^{-1})} \le \frac{1}{\mbox{{\em dist}}(\xi,\Sigma(\LL_1))}
$$
in the space $L^2(M^{-1})$, extends into 
$$
\forall \, \xi \in \Aa, \quad \| \Rr_1(\xi) \|_{H^k(M^{-1})} \le \frac{C_k}{\mbox{{\em dist}}(\xi,\Sigma(\LL_1))}, 
$$
with 
$$
\Aa = \left\{ \xi \in \mathbb{C}, \quad \mbox{arg}(\xi +\lambda) 
\in \left[-\frac{3 \pi}{4}, \frac{3 \pi}{4} \right] \ \mbox{ and } \ \Re e \, \xi \le \frac{\lambda}2 \right\},
$$
for any $k \in \N$ and some explicit constant $C_k >0$, 
\end{prop} 

\smallskip\noindent
{\sl Proof of Proposition~\ref{prop:resolHk}. } 
A quick way to prove the result for instance is the following. 
It is easy to prove by induction on $k \in \N$ the following estimate 
on the Dirichlet form: 
$$
\sum_{|s| \le k} a_s \, \langle \nabla^s \LL_1(g), \nabla^s g \rangle_{L^2(M^{-1})} 
\le - \tau_k \, \left( \sum_{|s| \le k} \| \bar \Pi ( \nabla^s g) \|^2 _{L^2(M^{-1})} \right)
$$
for some explicit $\tau_k >0$ and $a_s > 0$, $|s|\le k$, and where $\bar \Pi$ denotes 
the orthogonal projection in $L^2(M^{-1})$ onto the functions with zero mass, momentum and 
energy. Therefore we deduce on $\hat \LL_1$ that its semigroup satisfies 
$$
\forall \, k \in \N, \quad \| e^{t \, \hat \LL_1} \|_{H^k(M^{-1})} \le C_k
$$
and that obviously the same is true on the stable subspace of functions with zero energy. 
Then by interpolation with the rate of decay of the semigroup for functions with zero 
energy in $L^2(M^{-1})$, we deduce that 
$$
\forall \, \varepsilon >0, \ k \in \N, \quad \| e^{t \, \hat \LL_1} \Pi \|_{H^k(M^{-1})} \le 
C_{\varepsilon,k} \, e^{-(\mu_2 - \varepsilon) \, t} 
$$
for some explict $C_{\varepsilon,k} >0$, and where $\Pi$ is the orthogonal 
projection in $L^2(M^{-1})$ onto functions with zero energy. This implies on the resolvent that 
for any $k \in \N$, 
$$
\forall \, \xi \in \Aa, \quad \| \Rr_1(\xi) \|_{H^k(M^{-1})} \le C' _k, 
$$
with 
$$
\Aa = \left\{ \xi \in \mathbb{C}, \quad \mbox{arg}(\xi +\lambda) 
\in \left[-\frac{3 \pi}{4}, \frac{3 \pi}{4} \right] \ \mbox{ and } \ \Re e \, \xi \le \frac{\lambda}2 \right\},
$$
for some explicit $C' _k >0$. Then the result follows by straightforward interpolation 
with the estimates on the resolvent in $L^2(M^{-1})$. 
\qed
\medskip

Then we can conclude to the following extension of point (ii) of Theorem~\ref{thL1}:
\begin{prop}\label{prop:elastderiv}
We have 
$$
\forall \, \xi \in \Aa, \quad \| \Rr_1 (\xi) \|_{\Ww^{k,1}_q(m^{-1})} \le a_{k,q} + \frac{b_{k,q}}{|\xi + \lambda|}, 
$$
with 
$$
\Aa = \left\{ \xi \in \mathbb{C}, \quad \mbox{{\em arg}}(\xi +\lambda) 
\in \left[-\frac{3 \pi}{4}, \frac{3 \pi}{4} \right] \ \mbox{ and } \ \Re e \, \xi \le \frac{\lambda}2 \right\}.
$$
for any $k,q \in \N$ and some explicit constant $a_{k,q}, b_{k,q} >0$. 
\end{prop}

\subsection{Decomposition of $\hat \LL_\alpha$ and technical estimates} 

We fix once for all some $\bar \mu \in (\mu_2,0)$ and we split the proof of Theorem~\ref{thLalpha} into 
four steps, detailed in the following four subsections. 

Let us introduce the operator 
$$
P_\alpha = \LL_1 - \LL_\alpha = \LL^+ _1 - \LL^+ _\alpha + \tau_\alpha \, \nabla_v \cdot (v \, \cdot).
$$

Our first step in this subsection is to estimate the convergence to $0$ of the first part of 
this operator in suitable norm. Namely we prove 

\begin{lem}\label{lem:cvgLalpha}
(i) For any $k,q \in \N$, there exists $C = C_{k,q,m}$ such that %for any $g \in W^{k,1}_{q+1}(m^{-1})$ 
$$
 \left\|  \LL^+ _\alpha (g) \right\|_{W^{k,1}_q(m^{-1})} \le C \, \|g\|_{W^{k,1}_{q+1} (m^{-1})}, 
 \quad
 \left\|  \LL_\alpha (g) \right\|_{W^{k,1}_q(m^{-1})} \le C \, \|g\|_{W^{k+1,1}_{q+1} (m^{-1})}.
$$
(ii) For any $k,q \in \N$, there is a constructive function $\eps : (0,\infty) \to (0,\infty)$ satisfying $\eps(\alpha) \to 0$ 
as $\alpha$ goes to $1$ and such that for any $g \in W^{k,1}_{q+1}(m^{-1})$
   \[ 
   \left\| \left( \LL^+ _1 - \LL^+ _\alpha \right) (g) \right\|_{W^{k,1}_q(m^{-1})} \le 
       \eps (\alpha) \, \|g\|_{W^{k,1}_{q+1} (m^{-1})}.
       \]
        
 \noindent
(iii) There exists $C \in (0,\infty)$ such that for any  $g \in W^{3,1}_{3}(m^{-1})$, we have 
\[ 
\left\| \left( \LL _1 - \LL_\alpha \right) (g) \right\|_{L^1(m^{-1})} \le 
       C \, (1 - \alpha) \, \|g\|_{ W^{3,1}_{3}(m^{-1})}.
\]
 \end{lem}

\smallskip \noindent
{\sl Proof of Lemma~\ref{lem:cvgLalpha}. } 
The case $k=q=0$ is proved in Proposition~\ref{ContQe2}. 
Then higher-order derivatives are obtained 
from the $L^1(m^{-1})$ estimates by straightforward differentiation, and 
the incoporation of polynomial weights is trivial. 
\qed

Now let us consider some $\xi \in \CC$ and let us define 
$$
A_\delta = \LL^+ _{1,\delta} - \LL^*
$$
and 
$$
B_{\alpha,\delta}(\xi) = \nu + \xi + \big(\LL^+ _{1,\delta} - \LL^+ _1 \big) + P_\alpha
$$
(let us recall that the approximation $\LL^+ _{1,\delta}$ was defined in the beginning of 
Subsection~\ref{subsec:impGM}. 
It yields the decomposition 
$$
\LL_\alpha - \xi = A_\delta - B_{\alpha,\delta}(\xi).
$$

Then we have the 
\begin{lem}\label{Bdeltaalpha} 
Let us consider any $k,q \in \N$ and $\xi$ such that $\Re e \, \xi \ge - \min \nu$. Then
\begin{itemize}
\item[(i)] For any $\delta >0$, the operator 
$A_\delta : L^1 \to W^{\infty,1} _\infty(m^{-1})$ is a bounded linear operator 
(more precisely it maps functions of $L^1$ into $C^\infty$ functions with compact support).
 
\item[(ii)] For $\delta \in [0,\delta^*]$ and $\alpha \in [\alpha_2,1]$ for some constructive 
$\delta^* >0$ and $\alpha_2 \in (\alpha_1,1)$ (depending on a lower 
bound on $\mbox{{\em dist}}(\xi,\nu(\R^N))$),  
the operator 
$$
B_{\alpha,\delta} : W^{k+1,1} _{q+1}(m^{-1}) \to W^{k,1}_q (m^{-1})
$$ 
is invertible

\item[(iii)] The inverse operator $B_{\alpha,\delta}(\xi) ^{-1}$ satisfies for 
$\delta \in [0,\delta^*]$ and $\alpha \in [\alpha_3,1]$:
$$
\left\| B_{\alpha,\delta}(\xi) ^{-1} \right\|_{W^{k,1}_q (m^{-1}) \to W^{k,1}_q (m^{-1})} 
\le \frac{C_1}{\mbox{{\em dist}}(\Re e \, \xi,\nu(\R^N))} 
$$
and 
$$
\left\| B_{\alpha,\delta}(\xi) ^{-1} \right\|_{W^{k+1,1}_{q} (m^{-1}) \to W^{k,1}_{q}(m^{-1})}
\le \frac{C_2}{\mbox{{\em dist}}(\xi,\nu(\R^N))}
$$
for some explicit constants $C_1,C_2>0$ depending on $k,q,\delta^*,\alpha_2$ and a lower bound on 
$\mbox{{\em dist}}(\Re e \, \xi,\nu(\R^N))$.
\end{itemize}
\end{lem}

\smallskip\noindent
{\sl Proof of Lemma \ref{Bdeltaalpha}. } 
For $\xi \in \nu(\R^N)^c$, it was proved 
in~\cite[Proposition~4.1, Theorem~4.2]{GM:04} 
the convergence to $0$ of $\big( \LL^+ _{1,\delta} - \LL^+ _{1} \big)$ 
as $\delta \to 0$ (which was done in $L^1_1(m^{-1}) \to L^1(m^{-1})$ in~\cite{GM:04} and 
is extended in any $W^{k,1}_{q+1}(m^{-1}) \to W^{k,1}_q(m^{-1})$ by Proposition~\ref{prop:cvgronLp}), 
we deduce as in~\cite{GM:04} that for $\delta$ small enough (depending on a lower bound on the 
coercivity norm of $\nu+\xi$, that is on a lower bound on $\mbox{dist}(\xi,\nu(\R^N))$), we have 
$$
\big\| \big(\LL^+ _{1,\delta} - \LL^+ _{1}\big) g\big\|_{W^{k,1}_q(m^{-1})} \le 
\frac12 \, \| (\nu + \xi) \, g \|_{W^{k,1} _{q+1}}. %\qquad g \in L^p _1(m^{-1}).
$$
%and so that this operator is coercive and therefore invertible with explicit bound on its inverse. 

It was also proved that $A_\delta$ maps functions of $L^1$ into $C^\infty$ functions 
with compact support (with explicit estimates).

Let now consider $B_{\alpha,\delta}(\xi)$ only in the case $k=q=0$ (estimates for 
higher-order derivatives and weights are obtained by straightforward differentiation and 
computations). 
From Lemma~\ref{lem:cvgLalpha} we have for 
$\alpha$ close enough to $1$ (depending on a lower bound on $\mbox{dist}(\xi,\nu(\R^N))$), 
$$
\left\| \left( \LL^+ _1 - \LL^+ _\alpha \right) g \right\|_{L^1(m^{-1})} \le 
       \frac12 \, \|(\nu + \xi) \, g\|_{L^1 (m^{-1})}.
$$
By considering the semigroup on $L^1(m^{-1})$ of $B_{\alpha,\delta}(\xi)$ and computing 
the evolution of the norm in symmetric form 
using the formula for the differentiation of the complex modulus of a function
$$
\nabla |h| = \frac{\nabla h \, \bar h + h \, \nabla \bar h}{2 \, |h|},
$$
it is easily seen that
$$
\big\| e^{B_{\alpha,\delta}(\xi) \, t} g \big\|_{L^1(m^{-1})} \ge \| (\nu + \Re e \, \xi) \, g\|_{L^1 (m^{-1})}
-  \frac12 \, \|(\nu + \Re e \xi) \, g\|_{L^1 (m^{-1})}
$$
and therefore for $\alpha$ close enough to $1$ 
(depending on a lower bound on $\mbox{dist}(\xi,\nu(\R^N))$), we deduce that 
$$
\big\| e^{B_{\alpha,\delta}(\xi) \, t} g \big\|_{L^1(m^{-1})} 
\ge \frac12 \,  \| (\nu + \Re e \, \xi) \, g\|_{L^1 (m^{-1})} 
$$
and thus that the operator is invertible with its inverse bounded by 
$$
\big\| B_{\alpha,\delta}(\xi) ^{-1} \big\|_{L^1(m^{-1})} \le \frac{2}{\mbox{dist}(\Re e \, \xi,\nu(\R^N))}.
$$
Moreover by computing separately the evolution of the $L^1(m^{-1})$ norm 
in non-symmetric form (thus keeping $\nu+\xi$ but creating a term of the form 
$\OO(1-\alpha)$ times a $W^{1,1}_1(m^{-1})$ norm) and the evolution 
of the $W^{1,1}_1(m^{-1})$ norm in symmetric form: it yields easily 
$$
\big\| e^{B_{\alpha,\delta}(\xi) \, t} g \big\|_{W^{1,1}(m^{-1})} 
\ge \frac12 \,  \| (\nu + \xi) \, g\|_{L^1 (m^{-1})} 
+ \frac12 \,  \| (\nu + \xi) \, \nabla_v g\|_{L^1 (m^{-1})}
$$
which implies the result, by droping the second term.
\qed

\subsection{Geometry of the essential spectrum and estimates on the eigenvalues} 

First concerning the geometry of the spectrum, following the same strategy 
as in~\cite[Subsection~3.2]{GM:04}  we can prove the 

\begin{prop} \label{prop:weyl} 
Let us pick any $k,q \in \N$ and $m$ a smooth exponential weight function (as defined in (\ref{defdem})). 
Then for any $\alpha \in [\alpha_2,1]$,  the spectrum of $\LL_\alpha$ in 
$\Ww^{k,1}_q(m^{-1})$ is composed of 
a part included in $\Delta_{-\nu_0}^c$ (containing 
all possible essential spectrum), and a remaining part included 
in $\Delta_{-\nu_0}$ exclusively composed 
of discrete eigenvalues. 
\end{prop} 

\smallskip\noindent
{\sl Proof of Proposition \ref{prop:weyl}. } 
We follow the same method as in the proof of \cite[Proposition~3.4]{GM:04}. 
One uses the decomposition 
$$
\LL_\alpha = A_\delta - B_{\alpha,\delta}(0),
$$
the compactness of the first part $A_\delta$ and the coercivity 
$$
\| B_{\alpha,\delta}(0) \|_{L^1(m^{-1})} \ge 
%\Big( \nu_0-c_0 \, \rho \, (1-\alpha)\Big) \, 
\| \nu \, g \|_{L^1 (m^{-1})} - \varepsilon(\delta) \, \| \nu \, g \|_{L^1 (m^{-1})}
$$
of the second part (where $\epsilon(\delta) \to 0$ as $\delta \to 0$). 
Then one applies Weyl's theorem and show that (for any $\delta >0$)
$\Delta_{-\nu_0 + \epsilon(\delta)}$ has to be a Fredhom set with indices $(0,0)$ 
(except possibly for a countable family of points) 
since $[a,+\infty)$ is included in the resolvent set for $a$ big enoug.
\qed

\smallskip
Second concerning the discrete part of the spectrum, that is the isolated eigenvalues with finite multiplicity, 
following the same strategy as in~\cite[Proof of Proposition~3.5]{GM:04} we can prove the 
\begin{prop} \label{prop:regvp} 
Let us fix $\bar \mu \in (\mu_2,0)$. 
%Then there is $\alpha_3 \in [0,1)$ such that 
Then for any $\alpha \in [\alpha_2,1]$ (where $\alpha_2$ is obtained from Lemma~\ref{Bdeltaalpha} 
for this choice of $\bar \mu$), for any $\mu \in \Delta_{\bar \mu}$ and $\phi \in W^{1,1} _1$ satisfying 
$$
\LL_\alpha(\phi) = \mu \, \phi
$$
in $L^1$, we have 
$$
\| \phi \|_{W^{k,1}(m^{-1})} \le C_{k,m} \, \|\phi\|_{L^1_2}
$$
for any $k \in \N$ and $m=\exp(-a \, |v|^s)$, $a>0$, $s \in (0,1)$, where the constant 
$C_{k,m}$ depends on $k$, $m$ and a lower bound on $\bar \mu-\mu$.
\end{prop} 

\smallskip\noindent
{\sl Proof of Proposition \ref{prop:regvp}. } 
Let us sketch the idea of the proof. We use the decomposition 
$$
0 = \LL_\alpha \phi - \mu \, \phi = A_\delta \phi - B_{\alpha,\delta}(\mu) \phi
$$
and the fact that for the choices made for 
$\mu$ and $\alpha$ in the assumptions we have 
(adjusting $\delta$ as in Lemma~\ref{Bdeltaalpha})
$B_{\alpha,\delta}(\mu)$ is invertible in any $W^{k,1}(m^{-1})$ with 
explicit bound, and $A_\delta$ maps $L^1$ into $C^\infty$ 
functions with compact support. \qed

\smallskip
\begin{rem}  An alternative proof could be to adapt the proof of Proposition~\ref{errorH1}.
\end{rem}

\subsection{Estimate on the resolvent and global stability of the spectrum}

\begin{lem}\label{Ralpha}  
Let us pick $k,q \in \N$ and $m$ a smooth exponential weight function (as defined in (\ref{defdem})) 
and consider the operator $\LL_\alpha$ in $\Ww^{k,1} _q(m^{-1})$. Then
\begin{itemize}
\item[(i)] For any $\xi \in R(\LL_1)$, there is $\alpha_\xi \in [\alpha_2,1)$ such that $\xi \in R(\LL_\alpha)$ for any $\alpha \in [\alpha_\xi,1]$.
\item[(ii)] More precisely, the resolvent $\Rr_\alpha(\xi)$ satisfies the two following estimates 
for $\alpha \in [\alpha_2,1)$:
\bean %\label{estim1Ra}
&&
\| \Rr_\alpha(\xi) \|_{\Ww^{k,1}_q(m^{-1})} \le 
\frac{C_1 + C_2 \, \| \Rr_1(\xi) \|_{W^{k+1,1}_{q+1}(m^{-1})}}
{1- C_3 \, (1-\alpha) \, \| \Rr_1(\xi) \|_{W^{k+1,1}_{q+1}(m^{-1})}} \\ %\label{estim2Ra}
&& \| \Rr_\alpha(\xi) \|_{\Ww^{k+1,1} _{q+1}(m^{-1}) \to \Ww^{k,1}_q(m^{-1})} \le 
{1 \over \delta(\xi)}  \, \frac{C_1' + C_2' \, \| \Rr_1(\xi) \|_{W^{k+1,1}_{q+1}(m^{-1})}}
{1- C_3 \, (1-\alpha) \, \| \Rr_1(\xi) \|_{W^{k+1,1}_{q+1}(m^{-1})}}
\eean
with $\delta(\xi) := \mbox{{\em dist}}(\xi,\nu(\R^N))$ and where the  constants $C_i, C'_i$, $i=1,2,3$ depend on a positive lower bound on  $\mbox{{\em dist}}(\Re e \, \xi,\nu(\R^N))$.
\item[(iii)] Finally, for any compact set $K \subset \rho(\LL_1)$ there exists 
$\alpha_K \in [\alpha_2,1)$, $C_K \in (0,\infty)$ such that 
\bean %\label{estim3Ra}
&&
\forall \, \xi \in K, \,\, \alpha \in (\alpha_K,1] \qquad \| \Rr_\alpha(\xi) \|_{\Ww^{k,1}_q(m^{-1})} \le C_K, \\
%\label{estim4Ra}
&&
\forall \, \xi \in K, \,\, \alpha,\alpha' \in (\alpha_K,1] \qquad 
\| \Rr_\alpha(\xi) \, h  - {\Rr_{\alpha'}}(\xi) \, h \|_{\Ll^1(m^{-1})} \le C_K \, (1-\alpha) \, \| h \|_{W^{3,1}_{3}}.
\eean
\end{itemize}
\end{lem}

\smallskip\noindent
{\sl Proof of Lemma \ref{Ralpha}. } We split the proof into three steps. 

\smallskip\noindent{\sl Step 1. }
Let us consider the following operator defined from $\Ww^{k,1}_q(m^{-1})$ to $W^{k+1,1}_{q+1}(m^{-1})$ 
(which is seen to be well-defined at a glance)  
$$
I_{\alpha,\delta}(\xi) := - B _{\alpha,\delta} (\xi) ^{-1} 
+ \mathcal{R} _1(\xi) \, A_\delta \, B _{\delta,\alpha} (\xi) ^{-1}. 
$$
Some straightforward computations show that 
$$
(\LL_\alpha - \xi) \, I_{\alpha,\delta}(\xi) = - A_\delta \, B _{\alpha,\delta} (\xi) ^{-1} + \mbox{Id} 
+ \Big[ \mbox{Id} - P_\alpha \,  \mathcal{R} _1(\xi) \Big] \, A_\delta \, B _{\alpha,\delta} (\xi) ^{-1}
$$
which simplifies into 
$$
(\LL_\alpha - \xi) \, I_{\alpha,\delta}(\xi) =: J_{\alpha,\delta}(\xi) := 
\mbox{Id} - P_\alpha \,  \mathcal{R} _1(\xi) \, A_\delta \, B _{\alpha,\delta} (\xi) ^{-1} 
=: \mbox{Id} - K_{\alpha,\delta}(\xi).
$$
First using that 
$$
\| P_\alpha h \|_{\Ww^{k,1}_q(m^{-1})} \le C \, (1-\alpha) \, \| h \|_{W^{k+1,1}_{q+1}(m^{-1})},
$$
the control of $\mathcal{R}_1(\xi)$ in $\Ww^{k+1,1}_{q+1}(m^{-1})$ 
and the regularization property of $A_\delta$ we deduce that 
$$
K_{\alpha,\delta}(\xi) = 
P_\alpha \,  \mathcal{R} _1(\xi) \, A_\delta \, B _{\delta,\alpha} (\xi) ^{-1} = \OO(1-\alpha)
$$
in the norm of bounded operators on $W^{k,1}_q(m^{-1})$, 
and therefore for $(1-\alpha)$ small enough (with explicit bound) we get that 
$$
\| K_{\alpha,\delta}(\xi) \|_{W^{k,1}_q(m^{-1})} \le C_3 \, (1-\alpha) \, \| \mathcal{R} _1(\xi) \|_{\Ww^{k+1,1}_{q+1}(m^{-1})} <1
$$
and $\mbox{Id}-K_{\alpha,\delta}(\xi)$ is invertible in  $W^{k,1}_q(m^{-1})$. 
As a consequence 
$$
(\LL_\alpha - \xi) \, I_{\alpha,\delta}(\xi) \, (\mbox{Id} - K_{\alpha,\delta}(\xi))^{-1} =\mbox{Id}_{\Ww^{k,1}_q(m^{-1})}
$$
and we have proved that $\LL_\alpha - \xi$ admits a right-inverse, namely 
so that $I_\alpha(\xi) \, (\mbox{Id} - K_{\alpha,\delta}(\xi))^{-1} $. This proves that 
the operator $\LL_\alpha - \xi$ is onto. 

\smallskip\noindent{\sl Step 2. }
In order to show that $\LL_\alpha - \xi$ is invertible and that we have identified the 
resolvent it remains to prove that it is one-to-one.  
Let us consider the eigenvalue equation $(\LL_\alpha - \xi)\, h =0$ which writes 
$$
(\LL_1 - \xi)h = P_\alpha h
$$
from which we deduce (using Proposition~\ref{prop:regvp}  to get regularity bounds 
on $h$)
\begin{eqnarray*}
 \| h \|_{\Ww^{k,1}_q(m^{-1})} &\le& \| \mathcal{R} _1(\xi) \|_{\Ww^{k,1}_q(m^{-1})} \, 
 \| P_\alpha h \|_{\Ww^{k,1}_q(m^{-1})} \\
 &\le& C \, (1-\alpha) \, \| \mathcal{R} _1(\xi) \|_{\Ww^{k,1}_q(m^{-1})} \, \| h \|_{\Ww^{k+1,1}_{q+1}(m^{-1})} \\
 &\le& C' \, (1-\alpha) \, \| \mathcal{R} _1(\xi) \|_{\Ww^{k,1}_q(m^{-1})} \, \| h \|_{\Ww^{k,1}_q(m^{-1})}. 
\end{eqnarray*}
Therefore for $(1-\alpha)$ small enough (depending on the norm of $\mathcal{R} _1(\xi)$) 
we have that necessarily $h=0$, and thus the operator $(\LL_\alpha - \xi)$ is one-to-one. 

For $\alpha$ satisfying all the previous conditions, the operator $(\LL_\alpha - \xi)$ is bijective 
from $\Ww^{k+1,1} _{q+1}(m^{-1})$ to $\Ww^{k,1} _q (m^{-1})$ and its inverse is given by 
$$
\mathcal{R}_\alpha (\xi) = I_{\alpha,\delta} (\xi) \, J_{\alpha,\delta}(\xi) ^{-1} 
$$
from which we get the desired bound on the resolvent thanks to the study 
of $B_{\alpha,\delta}(\xi)^{-1}$ in Lemma~\ref{Bdeltaalpha}.  At this point we have 
proved points (i), (ii) and the first estimate in (iii). 

\smallskip\noindent{\sl Step 3. } The second estimate in point (iii) is obtained 
from the resolvent identity
$$
\Rr_\alpha(\xi) - \Rr_1(\xi) = \Rr_\alpha(\xi) \, [\LL_1 - \LL_\alpha] \, \Rr_1(\xi),
$$
together with the previous estimates on the resolvent and point (iii) in Lemma~\ref{lem:cvgLalpha}.
\qed

\medskip
Remark that this lemma proves the point (ii) in Theorem~\ref{thLalpha}. 
Moreover, as a consequence of this estimate on the resolvent $\Rr_\alpha(\xi)$, 
we may go one step further in the localization of the spectrum of $\hat \LL_\alpha$ around $0$.

\begin{cor} \label{coro:cvspec}
Let us fix $\bar \mu \in (\mu_2,0)$. In any $\Ww^{k,1}_q(m^{-1})$ there is
some constant $C \in (0,\infty)$ such that 
$$
\forall \, \alpha \in [\alpha_2,1], \quad
\Sigma(\hat \LL_\alpha) \cap \Delta_{\bar \mu} \subset B(0, C \, (1- \alpha) ). 
$$ 
\end{cor}

\smallskip\noindent
{\sl Proof of Corollary~\ref{coro:cvspec}}. The proof follows from the 
estimates in point (ii) of Lemma~\ref{Ralpha}, together with the fact that 
(Proposition 4.1 of \cite{GM:04} in $\Ll^1(m^{-1})$ extended to $\Ww^{k,1} _q (m^{-1})$ 
by the previous discussion): 
$$
\forall \, \xi \in \Delta_{\bar \mu}, \quad \| \Rr_1(\xi) \|_{\Ww^{k,1}_q(m^{-1})}
\le a + \frac{b}{|\xi|}
$$ 
for some explicit constants $a,b>0$. 
We get thus that $ \| \Rr_\alpha(\xi) \|_{\Ww^{k,1}_q(m^{-1})} < \infty$ if 
$\xi \in \Delta_{\bar\mu}$ and $|\xi| \ge C \, (1-\alpha)$, which concludes the proof. \qed

\subsection{Fine study of spectrum close to $0$}

Let us fix $r \in (0,|\bar\mu|]$ and let us choose any $\alpha_r \in [\alpha_2,1)$ such 
that $C \, (1- \alpha_r) < r$ (with the notations of Corollary~\ref{coro:cvspec}) in such a way that 
$\Sigma(\hat \LL_\alpha) \cap \Delta_{\bar \lambda} \subset B(0, r)$ for any $\alpha \in [\alpha_r,1]$. 
We may then define the  spectral projection operator (see~\cite{Kato})
\beqn\label{defPia}
\Pi_\alpha := - {1 \over 2 \, \pi \, i} \int_{S(0,r)} \Rr_\alpha(\zeta) \, d\zeta 
%= - {1 \over 2 \, \pi \, i} \int_{\zeta = - r} \Rr_\alpha(\zeta) \, d\zeta (?) 
\eeqn
in any $\Ww^{k,1}_q(m^{-1})$, with $S(0,r) := \{ \xi \in \CC, \, |\xi|=r \}$. The operator $\Pi_\alpha$ is the projection operator on the sum of  eigenspaces associated to eigenvalues lying in the half plane $\{ \xi \in \CC, \, \Re e \, \xi \ge - r \}$, see \cite{Kato}.  In particular the operator  $\Pi_1$ is the projection on the energy eigenline $\R \, \phi_1$, where we recall that  $\phi_1$ is the energy eigenfunction defined by (\ref{phiaTOphi1}).

\smallskip
\begin{lem}\label{Pia} 
The operator $\Pi_\alpha$ satisfies 
\begin{itemize}
\item[(i)] For any $k \in \N$  and any exponential weight function $m$ (as defined in (\ref{defdem})), 
it is well-defined and bounded in $\Ww^{k}_q(m^{-1})$.
\item[(ii)] Moreover there is a constant $C>0$ (depending on $m$) such that 
\begin{equation}\label{Pia-Pia'}
\forall \, \alpha, \, \alpha' \in [\alpha_r,1] \qquad \left\| \Pi _\alpha - \Pi _{\alpha'} \right\|_{W^{3,1} _{3}(m^{-1}) \to L^1(m^{-1})} \le 
C \, |\alpha' -\alpha|.
\end{equation}
\end{itemize} 
\end{lem}

\smallskip\noindent{\sl Proof of Lemma~\ref{Pia}.} It is a straightforward 
consequence of (\ref{defPia}) and Lemma~\ref{Ralpha}. 
\qed

\smallskip
\begin{cor}\label{spa3} 
There exists $\alpha_3 \in [\alpha_2,1)$ such that for any $\alpha \in [\alpha_3,1)$ there holds 
$$
\Sigma (\hat \LL_\alpha) \cap \Delta_{\bar\mu} = \{ \mu_\alpha \}
\quad\hbox{and the eigenspace associated to $\mu_\alpha \in \R$ is 1-dimensional.}
$$
This eigenvalue is called the energy eigenvalue. We may  furthermore remark that 
Corollary~\ref{coro:cvspec} implies
\beqn\label{estim1mua}
\forall \, \alpha \in [\alpha_3,1) \qquad |\mu_\alpha| \le C \, (1-\alpha).
\eeqn
 \end{cor}

\smallskip\noindent
{\sl Proof of Corollary~\ref{spa3}. } 
We already know that $\Sigma (\hat \LL_\alpha) \cap \Delta_{\bar\mu}$ is entirely 
composed of discrete spectrum. Therefore we have to prove that it is of dimension $1$. 
Indeed once this is proved, the fact that $\mu_\alpha \in \R$ is trivial since the operator is real, 
and the control (\ref{estim1mua}) is trivial from Corollary~\ref{coro:cvspec}.

Let us define the space $X_\alpha := \Pi_\alpha(L^1(m^{-1})) +\Pi_1(L^1(m^{-1}))$ 
endowed with the norm $\| \cdot \|_{L^1(m^{-1})}$. 
From Proposition~\ref{prop:regvp}, there exists a constant $C_1 >0$ such that 
$$
\forall \, \psi \in X_\alpha, \quad \| \psi \|_{W^{3,1}_{3}(m^{-1})} \le C_1 \,  \| \psi \|_{L^1(m^{-1})}.
$$
Thanks to the definition of $\Pi_\alpha$ and $\Pi_1$ and to Lemma~\ref{Pia}, we then get
\bean
\| \Pi_\alpha -  \Pi_1 \|_{X_\alpha \to X_\alpha} 
&\le& C_2 \, \sup_{\psi \in X_\alpha }\sup_{z \in S(0,r)} { \| (\Rr_{\alpha}(z) - \Rr_1(z)) \, \psi  \|_{L^1(m^{-1})} 
\over \| \psi \|_{L^1(m^{-1})} }\\
&\le& C_2 ' \, (1-\alpha) \, \sup_{\psi \in X_\alpha}  {\| \, \psi  \|_{W^{3,1}_{3}(m^{-1}))} \over \| \psi \|_{L^1(m^{-1})}} \\
&\le& C_2 '' \, (1-\alpha)  < 1,
\eean
for $(1-\alpha)$ small enough. By classical operator theory (see for instance the arguments 
presented in \cite[Chap 1, paragraph 4.6]{Kato} in order to prove \cite[Lemma 4.10]{Kato}) one deduces 
that $\hbox{dimension} (\Pi_\alpha) = \hbox{dimension} (\Pi_1)$. Since $\hbox{dimension} (\Pi_1)=1$ 
(as recalled in Theorem~\ref{thL1}), this concludes the proof. 
\qed

\medskip
Let us  introduce for any $\psi \in L^1$ the decomposition
$$
\psi = \Pi_1 \psi + \Pi_1^\perp \psi = (\pi_1 \psi) \, \phi_1 + \Pi_1^\perp \psi,
$$
where $\pi_1 \psi \in \R$ is the coordinate of $\Pi_1 \psi $ on $\R \, \phi_1$ (defined 
thanks to the projection $\Pi_1$). 
For any $\alpha \in [\alpha_3,1)$ we denote by $\phi_\alpha$ the unique eigenfunction
associated to $\mu_\alpha$ such that $\| \phi_\alpha \|_{L^1_2} = 1$ and $\pi_1 \phi_\alpha \ge 0$. 

\medskip
We can now establish a first order approximation of the eigenfunction $\phi_\alpha$.

\begin{lem} \label{lem:fp1}
For  any $k,q \in \N$ and any exponential weight function $m$ (as defined in~(\ref{defdem})), 
there exists $C$ such that 
\beqn\label{phi1-phia}
\forall \, \alpha \in [\alpha_3,1] 
\qquad \|\phi _\alpha - \phi _1 \|_{W^{k,1}_q (m^{-1})} \le C \, (1-\alpha).
\eeqn
\end{lem}

\smallskip
\begin{rem} \label{rem:fp1} We immediately deduce from Lemma~\ref{lem:fp1}  that $\phi_\alpha(0) < 0$ for
$\alpha$ close enough to $1$, and therefore, we get that this definition of $\phi_\alpha$ coincides with the definition 
in Theorem~\ref{theo:uniq}.
\end{rem}

\smallskip\noindent
{\sl Proof of Lemma \ref{lem:fp1}. } On the one hand, from the normalization conditions, we have
\bean
\| \phi_1 - \Pi_1 \phi_\alpha \|_{L^1_2}
&=& | 1 - \pi_1 \phi_\alpha | =  \left| \, \| \phi_\alpha \|_{L^1_2} - \| \Pi_1 \phi_\alpha \|_{L^1_2} \, \right| \\
&\le&  \| \phi_\alpha  - \Pi_1 \phi_\alpha \|_{L^1_2}  = \| \Pi_1^\perp \phi_\alpha \|_{L^1_2}.
\eean
We then deduce
\beqn\label{fp1:1}
\| \phi_1 - \phi_\alpha \|_{L^1_2} \le \| \Pi_1 \phi_\alpha - \phi_\alpha \|_{L^1_2} + \| \Pi_1^\perp \phi_\alpha \|_{L^1_2}
\le 2 \, \| \Pi^\perp_1 \phi_\alpha \|_{L^1_2}. 
\eeqn

On the other hand, the eigenfunction $\phi_\alpha$ satisfies 
$$
\hat\LL_1(\phi _\alpha) = [\hat\LL_1(\phi _\alpha)-\hat\LL_\alpha(\phi _\alpha)] -  \mu _\alpha \, \phi _\alpha. 
$$ 
Recall that from Proposition~\ref{prop:regvp} one has uniform bounds in $W^{\infty,1} _\infty(m^{-1})$ 
on $\phi_\alpha$ in terms of its $L^1_2$ norm which has been fixed to $1$, so that 
for any $\alpha \in [\alpha_3,1]$, $\| \phi_\alpha \|_{W^{k,1}_q(m^{-1})} \le C$. 
Using Proposition~\ref{ContQe1} and Proposition~\ref{prop:regvp} we get 
$$
\| \hat\LL_1 \phi _\alpha \|_{L^1 (m^{-1})} = \OO(1-\alpha).
$$
Using that $\hat\LL_1$ is invertible from $\Pi_1^\perp \Ll^1 _{1}(m^{-1})$ to $\Ll^1 (m^{-1})$ 
we deduce that 
\beqn\label{fp1:2}
\| \Pi^\perp_1 \phi_\alpha \|_{L^1 (m^{-1})}   = \OO(1-\alpha).
\eeqn
We conclude the proof of (\ref{phi1-phia}) holds for the $L^1_2$ norm gathering (\ref{fp1:1}) and (\ref{fp1:2}):
$$
\forall \, \alpha \in [\alpha_3,1] 
\qquad \|\phi _\alpha - \phi _1 \|_{L^1 _2} \le C \, (1-\alpha).
$$

Let now consider the eigenfunctions $\Phi_\alpha$ associated to $\mu_\alpha$ for $\alpha \in [\alpha_3,1]$ such
that $\pi_1 \Phi_\alpha > 0$ with the normalization condition $\| \Phi_\alpha \|_{W^{k,1}(m^{-1})} = 1$. 
Proceeding similarly as before (by working in the space $W^{k,1}(m^{-1})$), we can get 
$$
\| \Phi_\alpha - \Phi_1\|_{W^{k,1}(m^{-1})}  = \OO(1-\alpha).
$$
Because the eigenspace associated to $\mu_\alpha$ is of dimension $1$, we have $\Phi_\alpha = c_\alpha \, \phi_\alpha$
for some constant $c_\alpha \in (0,\infty)$. Then 
\bean
|c_1 - c_\alpha| &=& \| c_\alpha \, \phi_\alpha - c_1 \,  \phi_\alpha \|_{L^1_2} \le 
\| \Phi_\alpha -\Phi_1 \|_{L^1_2} + |c_1| \, \|  \phi_1 -  \phi_\alpha \|_{L^1_2} = \OO(1-\alpha).
\eean  
We then easily conclude that (\ref{phi1-phia}) holds for any $W^{k,1}(m^{-1})$ norm.
\qed

\medskip
We now use the linearized energy dissipation equation to get 
a second order expansion of the eigenvalue. 

\begin{lem} \label{lem:vp2}
For $\alpha \in [\alpha_3,1]$, the eigenvalues $\mu _\alpha$ 
satisfies (with explicit bound)
$$
\mu _{\alpha} = - \rho \, (1-\alpha) + \OO(1-\alpha)^2.
$$
\end{lem}

\smallskip\noindent
{\sl Proof of Lemma \ref{lem:vp2}. } 
By integrating the eigenvalue equation
$$
\hat\LL_\alpha \phi_\alpha = \mu_\alpha \, \phi_\alpha
$$
 against $|v|^2$ and dividing it by $(1-\alpha)$, we get
$$
\frac{\mu_{\alpha}}{1-\alpha}  \, \EE(\phi _{\alpha})= 
2 \, \rho \, \EE(\phi_{\alpha}) - 2 \, (1+\alpha) \, 
\tilde D ( \bar G_\alpha, \phi_{\alpha}).
$$
Using the rate of convergence of
$\bar G_\alpha \to \bar G_1$ and  $\phi_\alpha \to \phi_1$ established 
in Lemma~\ref{Contderivee=1Ge} and Lemma~\ref{lem:fp1} we deduce that 
\beqn\label{vp2:1}
\frac{\mu_{\alpha}}{1-\alpha}  \, \EE(\phi _{1}) = 
2 \, \rho \, \EE(\phi_{1}) - 4 \, 
\tilde D ( \bar G_1, \phi_1) + \OO(1-\alpha).
\eeqn
Then we compute thanks to (\ref{Mv2}) and (\ref{Mv4})
\beqn\label{vp2:2}
\EE(\phi_1) = 2 \, N \, c_0 \, \rho \, \bar\theta_1^2, 
\eeqn
where $c_0$ is still the normalizing constant in (\ref{phiaTOphi1}) such that $\| \phi_1 \|_{L^1_2} = 1$. Similarly, using (\ref{MMu3}),  (\ref{MMv2u3}) and the relation (\ref{tempSS}) which make a link between $b_1$ and $\bar\theta_1$, 
we find
\beqn\label{vp2:3}
\tilde D ( \bar G_1, \phi_1) = {3 \over 2} \, N \, c_0 \, \rho^2 \, \bar \theta_1 ^2.
\eeqn
We conclude gathering (\ref{vp2:1}), (\ref{vp2:2}) and (\ref{vp2:3}). 
\qed

%%%%%%%%%%%%%%%%%%%%%%%%%%%%%%%%%%%%%%%%%%%%%%%%
\subsection{The map $\alpha \mapsto \bar G_\alpha$ is $C^1$}

The fact that the path of self-similar profiles $\alpha \mapsto \bar G_\alpha$ is $C^0$ 
on $[\alpha_3,1]$ and $C^1$ at $\alpha =1$ was already proved in Lemma~\ref{Contderivee=1Ge}. 
Therefore we have to prove that it is $C^1$ for $\alpha \in [\alpha_3,1)$. 

Let us define the functional
$$
(\alpha,g) \mapsto \Psi(\alpha,g) := Q_\alpha(g,g) - \tau_\alpha \, \nabla_v(v \, g).
$$
The map $\Psi$ is $C^1$ from $\R \times (W^{1,1}_1(m^{-1})\cap \Cc_{\rho,0})$ into $\Ll^1(m^{-1})$ 
and it is such that for any $\alpha \in [\alpha_1,1)$, the equation  $\Psi(\alpha,g) = 0$ has 
only one solution which is the profile $\bar G_\alpha$. 
Moreover, for any $\alpha \in [\alpha_3,1)$, the linearized operator 
$D_2 \Psi(\alpha,\bar G_\alpha) = \LL_\alpha$ is invertible from $\Ww^{1,1}(m^{-1})$ into $\Ll^1(m^{-1})$ 
because of the spectral properties of $\LL_\alpha$ established in 
Theorem~\ref{thLalpha} (i) \& (ii) (note that here there is no eigenvalue approaching 
$0$ at $\alpha$).  Then using the same strategy as in 
Subsection~\ref{derivee=1Ge} based on the implicit function theorem 
we easily conclude that $\alpha \mapsto \bar G_\alpha$ is $C^1$ 
from $[\alpha_3,1)$ into $L^1(m^{-1})$. That ends the proof of Theorem~\ref{theo:uniq} (ii).

%%%%%%%%%%%%%%%%%%%%%%%%%%%%%%%%%%%%%%%%%%%%%%%%
\subsection{Decay estimate on the semigroup}\label{subsec:decaysemi}

We start with a lemma on non sectorial semigroups in Banach spaces. This 
result is a tool for deriving constructive decay rate on non sectorial 
semigroups, from the knowkedge on the resolvent of their generator. We do not 
try to prove such a decay rate for the semigroup in the norm of the Banach space 
but instead in a weaker norm (corresponding to the norm of the graph of some 
power of its generator), which shall be sufficient for our study of the linearized 
stability of the non-linear equation (\ref{eqresca}). 

\begin{lem}\label{Lemsemigroup} 
Let $A$ be a closed unbounded operator on a Banach space $E$ with dense domain 
$\hbox{{\em dom}}(A)$. We denote by $S(t)$ the associated semigroup, 
by $R(A)$ the associated resolvent set and by $\Rr=\Rr(\xi)$ the resolvent 
operator defined on $R(A)$. Assume that we have a sequence of Banach spaces 
$E_2 \subset E_1 \subset E_0=E$ decreasing for inclusion 
(in most cases this sequence shall be provided by $E_k = \hbox{dom}(A^k)$ 
endowed with the norm of the graph of $A^k$). We assume on the operator that:
\begin{itemize}
\item[(i)] the resolvent set $R(A)$ contains the half plan $\Delta_a$ for some $a \in \R$, 
together with the estimates
$$
\sup_{s \in \R} \| \Rr (a+i \, s) \|_{E_0 \to E_0} \le C_1
$$
and 
$$
\forall \, s \in \R, \quad \| \Rr (a+i \, s) \|_{E_1 \to E_0}, \  
\| \Rr (a+i \, s) \|_{E_2 \to E_1} \le \frac{C_2}{1+|s|}
$$
for some constants $C_1,C_2>0$;

\item[(ii)] the semigroup $S(t)$ satisfies
\beqn\label{Stle1}
\forall \, t \ge 0, \quad \| S(t) \|_{E_2 \to E_0}  \le C_3 \, e^{b \, t}
\eeqn
for some constants $C_3, b > 0$. 
\end{itemize}

Then for any $a' > a$, there exists a constant $C_4$ depending only 
on $a,b,a',C_1,C_2,C_3$ such that 
\beqn\label{Stforall}
\forall \, t \ge 0, \quad \| S(t) \|_{E_2 \to E_0}  \le C_4 \, e^{a' \, t}.
\eeqn
\end{lem}

\medskip\noindent
{\sl Proof of Lemma~\ref{Lemsemigroup}. } We split the proof into two parts. 

\smallskip\noindent{\sl Step 1. } 
The first bound on the resolvent implies that for any $x \in E_0$
$$ \| \Rr( a + i s) x \|_{E_0} \to 0, \quad |s| \to \infty. $$
%, et donc \'egalement
%$$ \lim_{|s| \to \infty} \int_{a \pm i s} ^{b \pm is} e^{z t} \, R( z) x \, dz =0. $$
Indeed we first consider $x \in \hbox{dom}(A)$ and then argue by density (since the domain $\hbox{dom}(A)$ 
is dense). When $x \in \hbox{dom}(A)$ the result is proved by the relation 
$$  
\Rr(z) x = z^{-1} \, \big[ - \mbox{Id} + \Rr(z) \, A \big] \, x. 
$$

\smallskip\noindent{\sl Step 2. } 
Then consider the following integral of $\Rr(z) x$ on a vertical segment with 
real part $a$ (for some $M>0$)
$$
I_M(x) := \int_{a-iM} ^{a+iM} e^{zt} \, \Rr(z) x \, dz. 
$$
The function $z \to \Rr(z)$ is differentiable on this segment and we can perform 
an integration by part: 
$$
I_M(x) = \frac{e^{(a+iM) \, t}}{t} \, \Rr(a + i M) x- \frac{e^{(a-iM) \, t}}{t} \, \Rr(a - i M) x
         - \int_{a-iM} ^{a+iM} \frac{e^{zt}}{t} \, \Rr(z)^2 x \, dz
$$
where we have used $R'(z) = R(z)^2$. Now we estimate the $E_0$ norm of this quantity: 
$$
\| I_M(x) \|_{E_0} \le \left\| \frac{e^{(a+iM) \, t}}{t} \, \Rr(a + i M) x \right\|_{E_0} 
+ \left\| \frac{e^{(a+iM) \, t}}{t} \, \Rr(a + i M) x \right\|_{E_0} 
$$
$$
+ C_2 ^2 \, \frac{e^{a \,t}}{t} \, \left( \int_{-\infty} ^{+\infty} \frac{1}{(1+|s|)^2} \, ds\right) \, \|x\|_{E_2}.
$$
Therefore the integral is semi-convergent and we can pass to the limit $M \to +\infty$ and 
use (see~\cite{Yao1995,Blake2001}) that 
$$
S(t) x = \frac{1}{2 i \pi} \, \lim_{M \to \infty} \int_{a-iM} ^{a+iM} e^{zt} \, \Rr(z) x \, dz =
\frac{1}{2 i \pi} \, \lim_{M \to \infty} I_M
$$
to obtain (the two boundary terms go to $0$ as $M \to +\infty$ from the first step) 
\beqn\label{Stge1}
\| S(t) x \|_{E_0} \le C_2' \, \frac{e^{a \,t}}{t}  \, \|x\|_{E_2}, \quad\hbox{with}\quad
C'_2 =  C_2 ^2 \,  \, \left( \int_{-\infty} ^{+\infty} \frac{1}{(1+|s|)^2} \, ds\right).
\eeqn
Using (\ref{Stle1}) for $t \le 1$ and (\ref{Stge1}) for $t \ge 1$, we conclude that (\ref{Stforall}) holds with
$C_4 = \max(C'_2, C_3 \, e^{b-a'})$.  
\qed

\medskip\noindent
{\sl Proof of point (iii) in Theorem~\ref{thLalpha}. } 
The point (ii) of Theorem~\ref{thLalpha}  was proved in Lemma~\ref{Ralpha} 
and it shows that the operator $\bar \LL_\alpha = (\mbox{Id} - \Pi_\alpha) \, \hat\LL_\alpha$ 
together with the sequence of Banach spaces  
$E_i = \Ww^{k+i,1}_i(m^{-1})$, $i=0,1,2$, for any fixed $k \in \N$ and any exponential weight function $m$ 
(as defined in~(\ref{defdem})), satisfies the assumption (i) of Lemma~\ref{Lemsemigroup} for 
any $a \in (\mu_2,0)$. Moreover it is trivial to prove that it satisfies the assumption (ii) of 
Lemma~\ref{Lemsemigroup} for some explicit $b>0$ from the decomposition 
$\LL_\alpha = A_\delta - B_{\alpha, \delta}(\xi)$ already introduced. 
\qed

%%%%%%%%%%%%%%%%%%%%%%%%%%%%%%%%%%%%%%%%%%%%%%%%
\section{Convergence to the self-similar profile}
\setcounter{equation}{0}
\setcounter{theo}{0}

%{\bf In this section we shall obtain some stability and convergence result of the rescaled 
%solution towards the self-similar profile. The idea is to connect the nonlinear tools of the elastic 
%case with the linearized theory. Indeed for a small inelasticity the $H$ theorem of the 
%nonlinear elastic theory is ``almost satisfied" and, together with the energy equation, 
%one naturally expects that it drives 
%the solution into some neighborood of the self-similar profile. This will be shown in the second  
%subsection. But a first difficulty arises already in the linearized theory, 
%since the stability result degenerates 
%as $\alpha \to 1$. However this degeneracy is due only to the effect of the 
%``energy" eigenvalue, and the nonlinear effects are equally degenerating to 
%$0$ along this eigenspace. Hence by balancing 
%the two estimates, one can solve this difficulty, as we shall see in the first subsection.  }

\smallskip
In this section, we consider the nonlinear rescaled equation (\ref{eqresca}) and we prove the convergence of its solutions to the self-similar profile. As a preliminary let us recall some result on propagation and appearance of moments and regularity which is picked up from \cite[Proposition 3.1, Theorem 3.5, Theorem 3.6]{MMII}.

\begin{lem}\label{estimMSg} 
Let us consider $g_{\mbox{\scriptsize{{\em in}}}} \in L^1_3 \cap \Cc_{\rho,0}$ and 
the associated solution $g \in C([0,\infty);L^1_3)$ to the rescaled equation (\ref{eqresca}). Then
\begin{itemize}
\item[(i)] 
For any exponential moment weight  $m$ (as defined in (\ref{defdem})) 
with exponent $s \in (0,1/2)$ and any time $t_0 \in (0,\infty)$, there exists a constant $M_1 = M_1(t_0)$ such that 
\beqn\label{appearL1m}
\sup_{[t_0,\infty)} \| g(t, \cdot) \|_{L^1(m^{-1})} \le M_1.
\eeqn
Moreover, if $ g_{\mbox{\scriptsize{{\em in}}}} \in L^1(m^{-1})$ for some polynomial 
or exponential (with exponent $s \in (0,1)$) moment weight $m$ then (\ref{appearL1m}) 
holds (for this weight $m$) with $t_0 = 0$ and some constant $M_1 = M_1( \| g_{in} \|_{L^1(m^{-1})} )$. 

\smallskip
For the two following points we now assume that for some constants $c_1, T \in (0,\infty)$ there holds
\beqn\label{lowerbdenergy}
 \inf_{[0,T]} \EE(g(t,\cdot)) \ge c_1,
 \eeqn
 and we state some smoothness properties of the solution $g$ which depend on $c_1$ but not on $T$ nor $\alpha$. 
 
\smallskip
\item[(ii)] Assume (\ref{lowerbdenergy}). Then for any $k_0 \in \N$ there is 
$q_0 = q_0(k_0) \in \N$ such that if 
$\| g_{\mbox{\scriptsize{{\em in}}}} \|_{ H^{k_0} \cap L^1_{q_0}} \le C_0$ holds,  
then for any $c_1 \in (0,\infty)$ there exists $C_1 = C_1(C_0,c_1) \in (0,\infty)$ 
such that for any time $T \in (0,\infty)$, we have 
\beqn\label{propagHk}
\forall \, t \in [0,T], \quad \| g(t, \cdot) \|_{H^{k_1}} \le C_1,
 \eeqn
 with $k_1 = 0$ if $k_0 = 0$ and $k_1 = k_0 - 1$ if $k_0 \in \N^*$.

\smallskip
\item[(iii)] Assume (\ref{lowerbdenergy}) and that $g_{\mbox{\scriptsize{{\em in}}}} \in L^2$, 
with $\| g_{\mbox{\scriptsize{{\em in}}}} \|_{L^2 \cap L^1_3} \le M_1 \in (0,\infty)$. 
Then there exists $\lambda \in (-\infty,0)$ and for any exponential weight function $m$ with exponent $s \in (0,1/2)$ and any $k \in \N$, there exists a constant $K$ (which depends on $\rho, c_1,M_1,k,m$) such that we may split $g = g^S + g^R$ with 
\beqn\label{concludegdecomp1}
\forall \, t \in [0,T], \quad  \| g^S(t,\cdot) \|_{H^k \cap L^1(m^{-1})} \le K, 
\quad \| g^R(t,\cdot) \|_{ L^1_3} \le K \, e^{\lambda \, t}.
\eeqn
\end{itemize}
\end{lem} 

\begin{rem}
It is worth mentioning that these estimates are uniform with respect to the inelasticity parameter 
$\alpha \in (0,1)$. Indeed, one the one hand, this was already the case for the moment estimate 
(\ref{appearL1m}) in \cite[Proposition 3.1]{MMII}. On the other hand (\ref{propagHk}) and 
(\ref{concludegdecomp1}) from \cite[Theorem 3.5, Theorem 3.6]{MMII} were (partially) based 
on the use of the damping effect of the anti-drift term (whose coefficient was fixed to $\tau=1$). 
Here the damping effect of the anti-drift term vanishes ($\tau_\alpha \to 0$) but it is replaced 
(as for the elastic Boltzmann equation) by the lower bound on the energy (\ref{lowerbdenergy}) 
which allows for a control from below on the convolution term $L(g)$ appearing in the loss term 
of the collision operator (see Lemma~\ref{estimatesonL}), which is enough to conclude also in this case.
\end{rem}

%%%%%%%%%%%%%%%%%%%%%%%%%%%%%%%%%%%%%%%%%%%%%%%

\subsection{Local linearized asymptotic stability}

Let us first consider the nonlinear evolution equation~(\ref{eqresc}) in $L^1(m^{-1}) \cap H^k$, 
and the associated equation on the fluctuation $h$ of a solution $g$ around the unique equilibrium 
$\bar G_\alpha$: $g = \bar G_\alpha + h$ and 
$$ %\label{eqNLh}
\partial_t h = \LL_\alpha h + Q_\alpha(h,h). 
$$

Let us start by stating an inequality that we shall need in the sequel. 

\begin{lem} \label{unifNL} 
For any exponential weight function $m$ (as defined in (\ref{defdem})), 
there is a constant $C \in (0,\infty)$ such that 
for any $h \in W^{3,1}_{3}(m^{-1})$ and any $\alpha \in (0,1)$, 
$$
\| \Pi _\alpha Q_\alpha(h,h) \|_{L^1(m^{-1})} 
\le C \, (1-\alpha) \, \| h \|_{W^{3,1}_{3}(m^{-1})}^2.
$$
\end{lem}

\medskip
\noindent{\sl Proof of Lemma~\ref{unifNL}.} We write 
$$
\Pi _\alpha Q_\alpha(h,h) =  \Pi _\alpha (Q_\alpha(h,h) -Q_1(h,h)) + 
(\Pi _\alpha - \Pi _1) Q_1(h,h).
$$
On the one hand, from Lemma~\ref{Pia} (i) and (\ref{convLetoL1}), there is $C \in (0,\infty)$ such that 
$$
\|  \Pi _\alpha (Q_\alpha(h,h) -Q_1(h,h)) \|_{L^1(m^{-1})} \le C \, (1-\alpha) 
\| h \|_{W^{3,1}_{3} (m^{-1})} ^2.
$$
On the other hand, from (\ref{Pia-Pia'}) and (\ref{QaWk1}), we get 
$$
\| (\Pi  _\alpha - \Pi  _1 ) Q_1(h,h) \|_{L^1(m^{-1})} \le 
C \, (1-\alpha) \, \| h \|^2_{W^{3,1} _{3} (m^{-1})}.  
$$
The proof of the lemma is immediate by gathering the two previous  estimates. 
\qed

\medskip
We now state a first local linearized stability result. 

\begin{prop}\label{domaineunif} 
For any $\alpha \in [\alpha_3,1)$, the self-similar profile $\bar G_\alpha$ 
is locally asymptotically stable, with domain of stability uniform according to $\alpha \in [\alpha_3,1)$. 

More precisely, let us fix $\rho \in (0,\infty)$ and some exponential weight function $m$ as in~(\ref{defdem}). 
There is $k_1,q_1 \in  \N^*$ such that for any $M_0 \in (0,\infty)$ there exists $C,\eps \in (0,\infty)$ such 
that for any $\alpha \in [\alpha_3,1]$, for any $g_{\mbox{\scriptsize{{\em in}}}} \in H^{k_1} \cap L^1(m^{-q_1})$ 
with mass $\rho$, momentum $0$ satisfying
\beqn\label{hypginGa}
\| g_{\mbox{\scriptsize{{\em in}}}} \|_{H^{k_1} \cap L^1(m^{-q_1})} \le M_0,
\qquad  \|g_{\mbox{\scriptsize{{\em in}}}} - \bar G_\alpha \|_{L^1(m^{-1})} \le \eps,  
\eeqn
the solution $g$ to the rescaled equation (\ref{eqresca}) with initial datum $g_{\mbox{\scriptsize{{\em in}}}}$ satisfies 
\beqn \label{energyGa}
\forall \,  t \ge 0, \quad \| \Pi_\alpha \, (g_t - \bar G_\alpha ) \|_{L^1(m^{-1})} 
\le C \, \| g_{\mbox{\scriptsize{{\em in}}}} - \bar G_\alpha \|_{L^1(m^{-1})} \, 
e^{\mu_\alpha  \, t},
\eeqn
\beqn\label{lingtoGa}
\forall \,  t \ge 0, \quad \| (\mbox{{\em Id}} - \Pi_\alpha) \, (g_t - \bar G_\alpha ) \|_{L^1(m^{-1})} 
\le C \, \| g_{\mbox{\scriptsize{{\em in}}}} - \bar G_\alpha \|_{L^1(m^{-1})} \, 
e^{(3/2) \, \mu_\alpha  \, t}.
\eeqn
\end{prop}

\medskip
\noindent{\sl Proof of Proposition~\ref{domaineunif}. Step 1. }
 Let us first denote by $c_1$ the constant given in Step 5 of Proposition~\ref{estimatesonGe} such that 
$$
\forall \, \alpha \in [\alpha_1,1), \quad \EE(\bar G_\alpha) \ge 2 \, c_1.
$$ 
We may then fix $\eps_0 \in (0,\infty)$ in such a way that 
\beqn\label{L2implEE}
\| g - \bar G_\alpha \|_{L^1(m^{-1})} \le \eps_0
\qquad\hbox{implies}\qquad \EE(g) \ge c_1,
\eeqn  
and define
$$
T_* := \sup \big\{T, \,\, \EE(g_t) \ge c_1 \,\,\, \forall \, t \in [0,T] \big\} \in (0,\infty].
$$
From Lemma~\ref{estimMSg} (i) \& (ii),  there exists $M \in (0,\infty)$ (depending on 
$\rho, c_1, k_1, q_1,M_0$) such that for any $T \in (0,\infty)$ there holds
\beqn\label{bdHkqT*}
 \sup_{t \in [0,T_*]} \| g \|_{H^{k_1} \cap L^1(m^{-q_1})} \le M.
\eeqn

Let us now consider the fluctuation $h_t = g_t - \bar G_\alpha$. Thanks to the mass and momentum 
conservations, it satisfies 
$h_t \in \Cc_{0,0}$ for all times, as well as the bound (\ref{bdHkqT*}).
We define the following decomposition on $h$: 
$$
h^1 = \Pi _\alpha h \ \mbox{ and } \ h^2 = (\mbox{Id} - \Pi _\alpha) h =: \Pi^{\bot} _\alpha h. 
$$
Since the spectral projection $\Pi_\alpha$ commutes with the linearized 
operator $\LL_\alpha$, the equation on $h^1$ writes 
$$
\partial_t h^1  =\mu_\alpha \, h^1 + \Pi _\alpha Q_\alpha (h,h).
$$
Multiplying that equation by $(\hbox{sign} \, h) \, m^{-1}$ and integrating in the velocity variable, we deduce  
thanks to Lemma~\ref{unifNL} and to (\ref{hWk1qm2}), (\ref{bdHkqT*}) that on $(0,T_*)$ the following holds
\bear \nonumber
\frac{d}{dt} \| h^1 \|_{L^1(m^{-1})}  &\le& \mu_\alpha\,  \| h^1 \|_{L^1(m^{-1})} 
+ C \, (1-\alpha) \,  \| h  \|_{W^{3,1}_{3}(m^{-1})}^2 
\\ \label{estimh1}
&\le& (1-\alpha) \, \Big[ C_1 \,   \| h^1 \|^{3/2}_{L^1_2} + C_1  \| h^2 \|^{3/2}_{L^1_2}- C_2 \, \| h^1 \|_{L^1(m^{-1})} \Big],
\eear
for some constants $C_1$ depending on $M$ and the possible choice $C_2 = \rho/2$ for $C_2$. 
For the second part $h^2$ we have the following equation 
$$ 
\partial_t h^2  = \Pi^\bot _\alpha \, \LL_\alpha \,  h^2 + \Pi_\alpha^\bot Q_\alpha (h,h).  
$$
Since the linearized operator $\LL_\alpha$ restricted to $\Pi^\bot _\alpha$ generates the semigroup $R_\alpha(t)$ defined in point (iii) of Theorem~\ref{thLalpha}, the Duhamel formula reads
$$
h^2(t) = R_\alpha(t)  \, h_{\mbox{\scriptsize{in}}} + \int_0^t R_\alpha(t-s) \,   \Pi_\alpha^\bot Q_\alpha (h,h)(s) \, ds.
$$
From (\ref{estimResta}) and (\ref{QaWk1}) we have 
$$
\| h^2(t) \|_{L^1(m^{-1})} \le C \, e^{\bar\mu \, t} \, \| h_{\mbox{\scriptsize{in}}} \|_{L^1(m^{-1})} +  C \int_0^t  e^{\bar\mu \, (t-s)} \, \| h(s) \|_{W^{2,1}_2(m^{-1})} ^2 \, ds.
$$
We deduce 
\begin{equation} \label{estimh2}
\| h ^2 _t \|_{L^1(m^{-1})} \le C_3 \, e^{\bar\mu \, t} \, \| h ^2 _{\mbox{\scriptsize{in}}} \| + 
C_4 \, \int_0 ^t e^{\bar\mu \, (t-s)} \, \Big( \| h ^1 _s \|_{L^1(m^{-1})}^{3/2} +   \| h ^2 _s \|_{L^1(m^{-1})}^{3/2} \Big) \, ds.
\end{equation}
with $C_4$ depending on $M$ thanks to (\ref{hWk1qm2}) and (\ref{bdHkqT*}).
It is then easy to show by comparison arguments from (\ref{estimh1}) and (\ref{estimh2}) 
that there are $0< \eps_2 \le \eps_1 \le \eps_0$ (one can take for instance $\eps_1 \le \eps_0/2$ 
satisfying $2 \, C_1 \, \eps_1^{1/2} < C_2$ and $2 \, C_4 \, \eps_1^{1/2} < 1/2$ and next $\eps_2 \le \eps_1$ 
satisfying $C_3 \, \eps_2 < \eps_1/2$) such that 
\beqn\label{h1+h2}
\| h ^1 _{\mbox{\scriptsize{in}}} \|_{L^1(m^{-1})} + \| h ^2 _{\mbox{\scriptsize{in}}} \|_{L^1(m^{-1})} \le \eps_2
\quad\hbox{implies}\quad
\sup_{t \in [0,T_*]} \max \Big\{ \| h ^1 _t \|_{L^1(m^{-1})}, \| h ^2 _t \|_{L^1(m^{-1})} \Big\} \le \eps_1.
\eeqn
Gathering (\ref{L2implEE}) and (\ref{h1+h2}) we deduce that there exists $\eps \in (0,\eps_2)$ such that 
under condition (\ref{hypginGa}) there holds $T_* = \infty$ as well as 
$$
\sup_{t \in (0,\infty)} \| g-\bar G_\alpha \|_{L^1(m^{-1})} \le 2 \, \eps_1 \le \eps_0.
$$

\medskip
\noindent{\sl Step 2. }
In a second step, coming back to (\ref{estimh2})  and to the integral version of (\ref{estimh1}) and setting $y(t) =  \| h^1 \|  + |\mu_\alpha| \,  \| h^2 \|$, we obtain
\beqn\label{ybdd1}
y(t) \le C_5 \, e^{\mu_\alpha \, t} \, y(0) + C_6 \,|\mu_\alpha| \,  \int_0^t e^{\mu_\alpha \, (t-s)} \, y(s)^{3/2} \, ds.
\eeqn
%We then follow the proof of \cite[Lemma 4.5]{GM:04}. Observing that we already know $0 \le y(t) \le C_7 \, \eps_1$ and taking $\eps_1$ smaller if necessary in order to get $C_6 \, C_7^{1/2} \, \eps_1^{1/2} \le 1/6$, we first deduce 
%\beqn\label{ybdd2}
%y(t) \le C_5 \, e^{\mu_\alpha \, (5/6) \, t} \, y(0)  \qquad \forall \, t \ge 0.
%\eeqn
Then we have to the following variant of the Gronwall lemma whose proof is the same that the one 
of \cite[Lemma~4.5]{GM:04} and is therefore skipped:
\begin{lem}\label{lem:GronModif}
  Let $y=y(t)$ be a nonnegative continuous function on $\R_+$ such that 
  for some constants $a$, $b$, $\theta$, $\mu>0$, 
    $$ %\label{eq:gron}
    y(t) \le a \, e^{-\mu t} X + b \, 
     \left( \int_0 ^t e^{-\mu(t-s)} y(s)^{1+\theta} \, ds \right) 
    $$
  (as compared to~\cite[Lemma~4.5]{GM:04}, $X$ needs not 
  necessarily be $y(0)$). 
Then if $X$ and $b$ are small enough, we have  
    $$    
    y(t) \le C \, X \, e^{- \mu t}.
    $$
for some explicit constant $C>0$. 
  \end{lem} 
Thanks to the uniform smallness estimate on $y(t)$ we can apply the lemma 
with $\theta = 1/4$ for instance, and we get 
$$ %\label{ybdd3}
y(t) \le  C_7 \, y(0) \, e^{\mu_\alpha \, t} 
$$
%Putting (\ref{ybdd2}) into (\ref{ybdd1}) we get 
%\beqn\label{ybdd3}
%y(t) \, e^{|\mu_\alpha| \, t} \le C_5  \, y(0) + C_6 \,|\mu_\alpha| \, (C_5 \, y(0))^{3/2} \,  \int_0^t e^{\mu_\alpha \, s /4} \, ds \le C_8 \, y(0),
%\eeqn
from which we deduce the estimate (\ref{energyGa}) for the $h^1$ part of $g-\bar G_\alpha$.
Finally, we may insert that estimate on $h^1$ in (\ref{estimh2}) and we get 
$$
\| h^2(t) \|_{L^1(m^{-1})} \le C'_3 \, (e^{\bar \mu \, t} 
+ e^{(3/2) \, \mu_\alpha \, t}) \, \| h_{\mbox{\scriptsize{in}}} \|_{L^1(m^{-1})} 
+ C_4 \int_0^t  e^{\bar\mu \, (t-s)} \, \| h^2(s) \|_{L^1(m^{-1})} ^{3/2} \, ds.
$$
The same kind of computation yields to 
$$
\ \| h ^2 _t \| \le C_8 \, e^{(3/2) \, \mu_\alpha \, t} \, \| h(0) \|
$$
from which (\ref{lingtoGa}) follows. 
\qed

%%%%%%%%%%%%%%%%%%%%%%%%%%%%%%%%%%%%%%%%%%%%%%%%%%
\subsection{Nonlinear stability estimates}

In this subsection we shall prove that when the inelasticity is small, depending on the size of 
the initial datum (but not necessary close to the self-similar profile), the equation~(\ref{eqresc}) 
is stable. This relies mainly on the fact that the entropy production timescale 
is of a different order (much faster) that the energy dissipation timescale as $\alpha \to 1$. 
This point is familiar to physicists (see for instance~\cite{BPlivre}) which separate, 
for granular gases with small inelasticity, the molecular 
timescale (the level where entropy production effects dominate) and the cooling timescale 
(much slower than the molecular timescale). 

  \begin{prop}\label{Attractiveg} 
  Define $k_2 := \max\{ k_0, k_1 \}$, $q_2 := \max \{ q_0,q_1, 3 \}$, where $k_i$ and $q_i$ are 
  defined in Theorem~\ref{CSK&EEP} and Corollary~\ref{DHeDH1}. For any $\rho, \, \EE_0, \, M_0$ 
  there exists $\alpha_4 \in [\alpha_3,1)$, $c_1 \in (0,\infty)$ and for any $\alpha \in [\alpha_4,1]$ 
  there exist $\varphi=\varphi(\alpha)$ with $\varphi(\alpha) \to 0$ as $\alpha \to 1$ 
  and $T=T(\alpha)$ (possibly blowing-up as $\alpha \to 1$) such that any initial datum 
  $0 \le g_{\mbox{{\em \scriptsize{in}}}} \in L^1 _{q_2} \cap H^{k_2} \cap \Cc_{\rho,0,\EE_0}$ with 
  $$
  \| g_{\mbox{{\em \scriptsize{in}}}}\|_{ L^1 _{q_2} \cap H^{k_2} } \le M_0,
  $$
the solution $g$ associated to the rescaled equation~(\ref{eqresca}) satisfies
  $$
  \forall \, t \ge 0, \quad \EE (g_t) \ge c_1 
  $$
  and for all $\alpha' \in [\alpha_4,1)$ and then all $\alpha \in [\alpha',1]$
  \beqn\label{estimateonGe3}
  \forall \, t \ge T(\alpha'), \quad  \left\| g_t - \bar G_\alpha \right\|_{L^1 _2} \le \varphi(\alpha').
  \eeqn
\end{prop}

\medskip
\noindent{\sl Proof of Proposition \ref{Attractiveg}.} 
Let us consider a solution $g \in C([0,\infty);L^1_{q_2} \cap H^{k_2})$ to the rescaled equation~(\ref{eqresca}) 
with given initial datum $ g_{\mbox{\scriptsize{in}}}$, whose existence has been established in \cite{MMRI,MMII}.
We split the proof of the Proposition into five steps.
\smallskip

\noindent
{\em Step~1}. From the propagation and appearance of uniform moment 
bounds~\cite[Proposition~3.1, (iii)]{MMII}, which it is worth noticing have been obtained 
uniformly with respect to the elastic coefficient (see also \cite{BGP**}), there exists $C_1 \in (0,\infty)$ such that 
\bear\label{unifL1mg}
&& \sup_{t \ge 0} \| g \|_{L^1_{q_2}} \le C_1.
\eear
Let us define $c_1 := \min \{ \EE(\bar G_1),\EE_0 \} /4$, and
\beqn\label{defT*}
T_* := \sup \big\{ T\ ; \ \forall \, t \in [0,T], \ \EE(g(t, \cdot )) \ge c_1 \big\}.
\eeqn
Next from the equation on the evolution of energy 
\beqn\label{eqEg}
\EE'(t) = - (1-\alpha^2) \, b_1 \, D_\EE(g) + (1-\alpha) \,2 \,  \rho \, \EE
\eeqn
and (\ref{unifL1mg}) there holds
$$ %\label{bd1Eg}
|\EE'(t)| \le C_2 \, (1-\alpha) \qquad \forall \, t \ge 0
$$
(take for instance $C_2 = 2 \, b_1 \, C_1^2 + C_1$), from which we deduce that 
we necessarily have 
$$ %\label{T*lowerbd}
T_* \ge C_3 \, (1-\alpha)^{-1}
$$
(take for instance $C_3 = (3/4) \, \EE_0/C_2$).

\smallskip\noindent
{\em Step~2}. 
%Observe that in~\cite[Section~3]{MMII}  the uniform $H^k$ bound on the solution $g$,  which was based on the use of the anti-drift  damping effect, can be obtained directly knowing some  lower bound on the energy. (Indeed the lower bound on the energy allows one to estimate from below the loss term in the collision operator and then to do not use the anti-drift  damping term).  
From point (ii) of Lemma~\ref{estimMSg}, we have for some constant $C_5 \in (0,\infty)$ 
\beqn\label{bdHkg}
\forall \, t \in [0,T_*] \qquad \| g_t \|_{H^{k_2}} \le C_4.
\eeqn
Moreover from Lemma~\ref{LowerBdga}, for any time $t_1 \in (0,T_*)$, there exists some constant 
$C_5 = C_5(\rho,C_4,t_1)$ such that 
\beqn\label{lowerbdg}
\forall \, t \in [t_1,T_*] \qquad g(t,v) \ge C_5^{-1} \, e^{-C_5 \, |v|^8}.
\eeqn

\smallskip\noindent
{\em Step~3}. With the notations of Theorem~\ref{CSK&EEP}, we compute the evolution of the 
relative entropy of $g(t, \cdot)$ with respect to the associated Maxwellian $M[g(t,\cdot)]$, and we obtain
\bean
\frac{d}{dt} H(g| M[g]) &=& \frac{d}{dt} H(g) - \frac{d}{dt} \int_{\R^N} g \, \ln \, M[g]
=  \frac{d}{dt} H(g) - {\rho \, N \over 2 \, \EE} \frac{d}{dt} \EE \\
&=&- D_{H,\alpha}(g) -  {\rho \, N \over 2 \, \EE} \, (1-\alpha) \, D_{\EE} (g).
\eean
Next from Lemma~\ref{DHeDH1} and the estimates 
(\ref{unifL1mg}), (\ref{defT*}), (\ref{bdHkg}) and (\ref{lowerbdg}) we have 
\begin{eqnarray*}
\frac{d}{dt} H(g|M[g]) 
&=& - D_{H,1}(g,g) +  \OO(1-\alpha) \quad\hbox{on}\quad (t_1,T_*).
\end{eqnarray*}
Then from  (\ref{EEPineg}), we are then led to the following differential inequation on the  relative entropy
$$
\frac{d}{dt} H(g|M[g]) \le - C_6 \, H(g|M[g])^2 + C_7 \, (1-\alpha) \quad \hbox{on}\quad (t_1,T_*). 
$$
By straightforward computations we deduce that {\em independently of the value 
of $H(g_{t_1}|M[g_{t_1}])$} (this ``loss of memory" effect is typical of differential equations 
with overlinear damping terms), we have 
$$
\forall \, t \in [t_1,T_*], \quad H(g_t|M[g_t]) \le C_8 \, (1-\alpha)^{1/2}\, 
\frac{1+ e^{-C_{9} \, (1-\alpha)^{1/2} \, (t-t_1)}}{1- e^{-C_{9} \, (1-\alpha)^{1/2} \, (t-t_1)}}
$$
for some explicit constants. As a conclusion, defining $t_2 := t_1 + C_9 ^{-1} \, (1-\alpha)^{-1/2}$ and choosing 
$\bar \alpha \in [\alpha_3,1)$ in such a way that $t_2 < T_*$ we have for $\alpha \in [\alpha',1)$ 
$$
\forall \, t \in [t_2,T_*]  \quad H(g(t)|M[g]) \le C_{10} \, (1-\alpha)^{1/2}.
$$
Finally, using Csisz\'ar-Kullback-Pinsker inequality (\ref{CKineg}), as well as H\"older inequality, we obtain under 
the same conditions on $\alpha$ and the time variable: 
\beqn\label{bdRH}
\| g - M[g] \|_{L^1_3} \le C \, \| g - M[g] \|_{L^1}^{1/2} \, \| g \|_{L^1_6}^{1/2} 
\le C \, H(g|M[g])^{1/4} \le C \, (1-\alpha)^{1/8}.
\eeqn

\smallskip\noindent
{\em Step~4}. Now let us go back to the energy equation (\ref{eqEg}). 
First, with the help of the moment bound (\ref{unifL1mg}), one may write
$$ %\label{eqEg2}
\EE'(t) = 2 \, (1-\alpha) \, [\rho \, \EE - b_1 \, D_\EE(g) + \OO(1-\alpha)].
$$
Thanks to (\ref{bdRH}) we deduce
$$ %\label{eqEg3}
\EE'(t) = 2 \, (1-\alpha) \, (\rho \, \EE - b_1 \, D_\EE(M[g]) + \OO((1-\alpha)^{1/8}))
\quad\hbox{on}\quad(t_2,T_*).
$$
Finally, thanks to (\ref{defPsi}), (\ref{propPsi}) and the relation $\EE(g) = \rho \, N \, \theta(g)$, we get on $(t_2,T_*)$
\beqn\label{eqEg4}
\EE'(t) := \Sigma(\EE(t),\alpha) =  (1-\alpha) \, [k_3 \, \EE  \, (\bar\EE_1^{1/2} - \EE^{1/2}) + \OO((1-\alpha)^{1/8})],
\eeqn
where $\bar\EE_1 = \rho \, N \, \bar\theta_1$ with $\bar\theta_1$ is the quasi-elastic self-similar temperature 
defined in~(\ref{tempSS}).  We may then choose $\alpha'' \in [\alpha',1)$ such that 
$\Sigma(c_1,\alpha) > 0$ for any $\alpha \in  [\alpha'',1)$. 
We conclude by maximum principle that $T_* = \infty$ for  $\alpha \in [\alpha'',1)$. 
In particular, all the previous estimates on $g$ are uniform on $(t_2,\infty)$. 

\smallskip
\noindent
{\em Step~5}. Thanks to (\ref{eqEg4}) we easily get 
$$ %\label{eqEg5}
{d \over dt} (\EE-\bar\EE_1)^2 \le - (1-\alpha) \, [k_5 \, (\EE-\bar\EE_1)^2 + \OO((1-\alpha)^{1/8})],
$$
so that (for some constants $a,b>0$)
$$ %\label{eqEg6}
\forall \, t \ge t_2, \quad |\EE(t) -\bar\EE_1| \le  |\EE(t_2) -\bar\EE_1| \, e^{-a\, (1-\alpha) \, (t-t_2)} + b \, (1-\alpha)^{1/8}.
$$
Setting $T(\alpha) = \max \{ t_2,c \, (1-\alpha)^{-1}\}$ for some suitable constant $c>0$, we then obtain
\beqn\label{eqEg6}
|\EE-\bar\EE_1|= \OO((1-\alpha)^{1/8}) \quad \hbox{on}\quad [T(\alpha),\infty).
\eeqn
In order to conclude that (\ref{estimateonGe3}) holds, we write 
$$
 g(t) - \bar G_\alpha = (g(t) - M[g(t)])  + (M[g(t)]- \bar G_1 ) + (\bar G_1 - \bar G_\alpha),
$$
and we estimate the first term thanks to (\ref{bdRH}), the second term thanks to (\ref{eqEg6}) 
and the third term by (\ref{estimateonGe2}).  \qed

\subsection{Decomposition and Liapunov functional for smooth initial datum}

The proof of the gobal convergence (point (v) of Theorem~\ref{theo:uniq}) for smooth initial data 
only amounts to connect the two previous results of Propositions~\ref{domaineunif} 
and~\ref{Attractiveg} by choosing $\alpha$ such that $\varphi(\alpha) \le \eps$ 
where $\eps$ is the size of the attraction domain in Proposition~\ref{domaineunif} and 
$\varphi(\alpha)$ is defined in  Propositions~\ref{Attractiveg}. More precisely, we state without 
proof the straightforward combination of Propositions~\ref{Attractiveg} and Proposition~\ref{domaineunif}.

\begin{cor}\label{coro:cvgHk}
Let us fix an exponential weight function $m$ as in (\ref{defdem}), 
with exponent $s \in (0,1)$.  
Then for any $\rho, \, \EE_0, \, M_0$  there exists $C$ and $\alpha_5 \in [\alpha_4,1)$ 
(depending on $\rho, \, \EE_0, \, M_0,m$) 
such that for any $\alpha \in [\alpha_5,1)$ and any initial datum 
$0 \le g_{\mbox{{\em \scriptsize{in}}}} \in L^1(m^{-q_2}) \cap H^{k_2} $ satisfying 
$$
g_{\mbox{{\em \scriptsize{in}}}}\in \Cc_{\rho,0,\EE_0}, \qquad 
\| g_{\mbox{{\em \scriptsize{in}}}}\|_{ L^1(m^{-q_2}) \cap H^{k_2} } \le M_0,
$$
the solution $g$ associated to the rescaled equation~(\ref{eqresca}) satisfies 
$$ %\label{energyGa2}
\forall \,  t \ge 0, \quad \| \Pi_\alpha \, (g_t - \bar G_\alpha ) \|_{L^1(m^{-1})} 
\le C \, e^{\mu_\alpha  \, t},
$$
$$ %\label{lingtoGa2}
\forall \,  t \ge 0, \quad \| (\mbox{{\em Id}} - \Pi_\alpha) \, (g_t - \bar G_\alpha )  \|_{L^1(m^{-1})} 
\le C \, e^{(3/2) \, \mu_\alpha  \, t}.
$$
\end{cor}
\begin{rem}
Note that the constant $C$ in the rate of decay does not depend on $\alpha$. This 
comes from the fact the size of the linearized stability domain is uniform as 
$\alpha$ goes to $1$ in Proposition~\ref{domaineunif}, which allows in Propositions~\ref{Attractiveg} 
to pick a fixed $\alpha'$ such that in the estimate (\ref{estimateonGe3}) 
$\varphi(\alpha')$ is less than this size, and therefore 
that the time $T(\alpha')$ required to enter this neighborhood does not blow-up as $\alpha$ 
goes to $1$.
\end{rem} 

As a by-product of the previous propositions, we state and prove a result which provides a 
partial answer to the question (important from the physical 
viewpoint) of finding Liapunov functionals for this particles system. 
Let us define the required objects. We consider a fixed mass $\rho$ and 
some restitution coefficient $\alpha$ whose range will be specified below. 
At initial times, non-linear effects dominate and therefore we define 
$$
{\mathcal H}_1 (g) := H(g|M[g]) + \big(\EE- \bar\EE_\alpha \big)^2
$$
where $\bar \EE_\alpha = \EE(\bar G_\alpha)$ is the energy of the self-similar profile corresponding to $\alpha$ 
and the mass $\rho$.
At eventual times, linearized effects dominate. Therefore we define a quite natural candidate 
from the spectral study: 
$$
{\mathcal H}_2(g) := \| h^1 \|_{L^1(m^{-1})}^2 + (1-\alpha) \int_0 ^{+\infty} \left\| R_\alpha(s) \, h^2 \right\|^2 _{L^2} \, ds,
$$
with $h^1 = \Pi_\alpha h$, $h^2 = \Pi_\alpha^\perp h$ and $h = g - \bar G_\alpha$.

\begin{prop} \label{prop:lyap} 
There is $k_4 \in \N$ big enough (this value is specified in the proof) such that for any exponential weight 
function $m$ as defined in (\ref{defdem}), any time $t_0 \in (0,\infty)$ and any $\rho,\EE_0, M_0 \in (0,\infty)$, 
there exists $\kappa_* \in (0,\infty)$ and $\alpha_6 \in [\alpha_5,1)$ such that for any $\alpha \in [\alpha_6,1]$ 
and initial datum $g_{\mbox{\scriptsize{{\em in}}}} \in H^{k_4} \cap L^1(m^{-1})$ satisfying
$$ %\label{hypginGabis}
g_{\mbox{\scriptsize{{\em in}}}} \in \Cc_{\rho,0,\EE_0}, 
\qquad \| g_{\mbox{\scriptsize{{\em in}}}} \|_{H^{k_4} \cap L^1(m^{-1})} \le M_0, \qquad
g_{\mbox{\scriptsize{{\em in}}}}(v) \ge M_0 ^{-1} \, e^{-M_0 \, |v|^8}, 
$$
the solution $g$ to the rescaled equation (\ref{eqresca}) with initial datum 
$g_{\mbox{\scriptsize{{\em in}}}}$ is such that the functional 
$$
{\mathcal H}(g_t) = {\mathcal H}_1(g_t) \, {\bf 1}_{\big\{{\mathcal H}_1(g_t) \ge \kappa_*\big\}} 
+ {\mathcal H}_2(g_t) \, {\bf 1}_{\big\{{\mathcal H}_1(g_t) \le \kappa_*\big\}} 
$$
is decreasing for all times $t \in [0,+\infty)$. Moreover, $\HH(g(t, \cdot))$ is strictly decreasing as 
long as $g(t, \cdot)$ has not reached the self-similar state $\bar G_\alpha$. 
\end{prop}

\medskip
\noindent{\sl Proof of Proposition \ref{prop:lyap}.} 
We split the proof into three steps. 

\smallskip\noindent{\sl Step 1: Initial times.} 
Taking $k_4 \ge k_2$ and $\alpha \in [\alpha_4,1)$, we know from the proof of 
Proposition~\ref{Attractiveg} that the solution $g$ satisfies that 
$$
\forall \, t \in [t_0,\infty),  \quad  \| g(t, \cdot) \|_{H^{k_4} \cap L^1(m^{-1})} \le M_1,
\quad
g(t,v) \ge M_1^{-1} \, e^{- M_1 \, |v|^8},
$$
for some constant $M_1 \in (0,\infty)$ (recall that 
$\alpha_4$ was adjusted in terms of $\rho,\EE_0, M_0$). 
Coming back then to Steps 3 and 4 in the proof of Proposition~\ref{Attractiveg}, we obtain
the two following differential equation on $(t_0,\infty)$
$$
\frac{d}{dt} H(g|M[g]) \le - K_1 \, H(g|M[g])^2 +  \OO(1-\alpha) 
$$
and 
$$
\frac{d}{dt}\EE =  2 \, \rho \, (1-\alpha) \, \left[  K_2 \, \EE \, (\bar \EE_\alpha ^{1/2} - \EE^{1/2} ) \, (\EE-\bar \EE_\alpha) 
+  \OO( (1-\alpha)^{1/8}) \right],
$$
for some constants $K_i \in (0,\infty)$. We easily deduce that for any 
$\kappa \in (0,\infty)$ there exists $\alpha_\kappa \in [\alpha_5,\infty)$ such that 
\beqn\label{HH1'<0}
\frac{d}{dt} {\mathcal H}_1(g_t) < 0 \ \mbox{ for any } \ 
t \in (0,\infty)\,\, \hbox{such that}\,\,\, {\mathcal H}_1(g_t) \ge \kappa.
\eeqn

\smallskip\noindent{\sl Step 2: Eventual times.} Let us first remark that from 
point (iii) in Theorem~\ref{thLalpha} (iii) and the interpolation inequality (\ref{hWk1qm2}), 
for any $q \in \N^*$ there exists $k,k' \in \N$ and $C_i \in (0,\infty)$ such that 
\bean
\left\| R_\alpha \, h^2 \right\|_{L^2}
&\le&  C_1 \, \left\| R_\alpha \, h^2  \right\| _{W^{k,1}(m^{-q/2})}\\
&\le& C_2 \, e^{\bar\mu \, s} \, \left\| h^2 \right\| _{W^{k+2,1}_2(m^{-q/2})}
\le  C_3 \, e^{\bar\mu\, s} \, \left\| h \right\| _{H^{k'} \cap L^1(m^{-q})},
\eean
so that, taking $k_4$ big enough, the functional $\HH_2(g(t,.))$ is well-defined 
for any times $t \in (0,\infty)$. 
First observe that from (\ref{estimh1}) there holds
\beqn\label{Liapestimh1}
\frac{d}{dt} \| h^1 \|^2_{L^1(m^{-1})}  \le (1-\alpha) \, \Big[ K_1 \,   \| h \|^{5/2}_{L^1(m^{-1})} 
- K_2 \, \| h^1 \|^2_{L^1(m^{-1})} \, \Big].
\eeqn
Second, we compute (with the notation of Subsection~\ref{subsec:decaysemi})
\bean
{d \over dt} \int_0 ^{+\infty} \left\| R_\alpha(s) \, h^2_t \right\|^2 _{L^2} \, ds = 
 2 \, \int_0 ^{+\infty} \int_{\R^N} (e^{s \, \bar \LL_\alpha} \, h^2 )\, [e^{s \, \bar  \LL_\alpha} \, (\bar \LL_\alpha h^2 + \Pi_\alpha^\perp Q_\alpha(h,h))] \, ds \, dv.
\eean
On the one hand, 
\bean
I_1 &=&  2 \, \int_0 ^{+\infty} \int_{\R^N} (e^{s \, \bar \LL_\alpha} \, h^2 )\, [e^{s \, \bar \LL_\alpha} \, \bar \LL_\alpha h^2 ] \, dsdv \\
&=&   \int_0 ^{+\infty} {d \over ds} \| e^{s \, \bar \LL_\alpha} \, h^2 \|_{L^2}^2 \, ds =  - \| h^2 \|_{L^2}^2.
\eean
On the other hand, 
\bean
I_2 
&=&  2 \, \int_0 ^{+\infty} \int_{\R^N} (R_\alpha (s) \, h^2 )\, [R_\alpha(s) \, \Pi_\alpha^\perp Q_\alpha(h,h)) ] \, dsdv \\
&\le&  2 \, C_1^2 \int_0 ^{+\infty} \| R_\alpha (s) \, h^2\|_{W^{k_1,1}(m^{-q/2})} \, 
\| R_\alpha(s) \, \Pi_\alpha^\perp Q_\alpha(h,h) \|_{W^{k_1,1}(m^{-q/2})}  \, ds \\
&\le&  C'_2 \left(  \int_0 ^{+\infty} e^{2\bar\mu \, s}  \, ds\right)   \| h^2\|_{W^{k_1+1,1}_2(m^{-q/2})} \, 
\| Q_\alpha(h,h) \|_{W^{k_1+1,1}_2(m^{-q/2})} \\
&\le&  C'_3 \,  \| h^2\|_{L^2}^{3/4} \,  \| h^2\|_{H^{k_3} \cap L^1(m^{-1})}^{1/4} \,   \| h\|_{L^2}^{3/2} \,  \| h\|_{H^{k_3} \cap L^1(m^{-1})}^{1/2} ,
\eean
for some $k_3 \in \N$ given by Proposition~\ref{interpolineg}. Taking $k_4 \ge k_3$, we then obtain
\beqn\label{Liapestimh2}
{d \over dt} \int_0 ^{+\infty} \left\| R_\alpha(s) \, h^2_t \right\|^2 _{L^2} \, ds \le
K_3 \, \| h \|_{L^2}^{9/4} - \| h^2 \|_{L^2}^2.
\eeqn
Gathering (\ref{Liapestimh1}) and (\ref{Liapestimh2}) and using some interpolation again, we deduce that there exists $\kappa' \in (0,\infty)$ such that 
\beqn\label{HH2'<0}
\frac{d}{dt} {\mathcal H}_2(g_t) < 0 \ \mbox{ for any } \  t \in (0,\infty)\,\, \hbox{such that}\,\,\, \| h_t \|_{L^1} \le \kappa'.
\eeqn

\smallskip\noindent{\sl Step 3.} 
We conclude putting together (\ref{HH1'<0}) and (\ref{HH2'<0}), and 
using (\ref{CKineg}), (\ref{estimateonGe2}) in order to prove that 
$$
\HH_1(g) \le \kappa \quad \hbox{implies}\quad \| h_t \|_{L^1} \le \kappa',
$$
for $\alpha \in [\alpha_6,1]$ for some $\alpha_6 \in [\alpha_5,1)$. 
\qed

%%%%%%%%%%%%%%%%%%%%%%%%%%%%%%%%%%%%%%%%%%%%
\subsection{Global stability for general initial datum}  

\smallskip

We first state and prove a regularity result on the iterated gain term which is the inelastic collision operator version of the same result proved for the elastic collision operator in \cite{MW99,Abraham}. 

\smallskip

 \begin{lem} \label{itarteQ+} There exists a constant $C$ such that for 
 any $f,g,h \in L^1_2(\R^N)$ and any $\alpha \in (0,1]$ there holds
\beqn\label{L3estim}
\| Q^+_\alpha(f,Q_\alpha^+(g,h)) \|_{L^3} \le C \, \| f \|_{L^1_2} \, \| g \|_{L^1_2} \,  \| h \|_{L^1_2}.
\eeqn
\end{lem}

\smallskip\noindent
{\sl Proof of Lemma~\ref{itarteQ+}. } We follow \cite[lemma 2.1]{MW99} and  \cite[lemma 2.1]{Abraham} 
and we make use of the Carleman representation introduced in~\cite[Proposition 1.5]{MMII}. 
Let us consider $f,g,h \in L^1_2(\R^N)$ and $\phi \in L^\infty(\R^N)$. 
We apply twice the weak formulation of the gain term
\bean
&& \int_{\R^N} Q^+(f,Q^+(g,h)) (v) \, \phi(v) \, dv \\
&&= \int_{\R^N} Q^+(g,h) (v) \, \left[ \int_{\R^N} f(v_2) \, |v-v_2| \int_{S^2} \, \phi(w'_2) \, d\sigma_2 dv_2 \right] dv \\
&&= \int_{\R^N} \!\! \int_{\R^N} \!\!  \int_{\R^N} g(v) \, h(v_1) \, f(v_2) \, \left[ 
|v-v_1|  \int_{\Sp^2}  |v'_1-v_2| \int_{\Sp^2} \, \phi(v''_2) \, d\sigma_2 \, d\sigma_1 \right] \, dv \, dv_1  \, dv_2
\eean 
with $w'_2 = V'(v_2,v,\sigma_2)$, $v' _1= V' (v,v_1,\sigma_1)$ and therefore
$v''_2 = V'(v_2,v'_1,\sigma_2)$. Recall that for any given $v,v_*,\sigma \in \R^N$, we define 
$$
w = {v+v_* \over 2}, \quad u = v-v_*, \quad \gamma = {1 +e \over 2}, \quad
u' = (1-\gamma) \, u + \gamma \, |u| \, \sigma
$$
and then 
$$
V' = V'(v,v_*,\sigma) = {w \over 2} + {u' \over 2} = v + {\gamma \over 2} \, (|u| \, \sigma - u)
$$ 
$$
V'_* = V'_* (v,v_*,\sigma) = {w \over 2} - {u' \over 2} = v_* -  {\gamma \over 2} \, (|u| \, \sigma - u).
$$
We denote by $\Phi = \Phi(v,v_1,v_2)$ the term between brackets in the last integral. 
Introducing the point $w_1$ and the set $S_{v,v_1,\eps}$ defined by 
$$
w_1 := (1-\gamma/2) \, v + (\gamma/2) \, v_1, \quad\,\, 
S_{v,v_1,\eps} := \left\{ z \in \R^N; \,\, \Big|  | z - w_1|  - (\gamma/2) \,  |v-v_1| \Big| \le  \eps /2\right\},
$$
we get 
\beqn\label{defPhi}
\Phi = {(2/\gamma)^2 \over |v-v_1|} \lim_{\eps \to 0} { \Psi_\eps \over \eps}, 
\quad  \Psi_\eps = \int_{\R^N} \! \int_{\Sp^2}  
{\bf 1}_{S_{v,v_1,\eps}}(v'_1) \, |v'_1 - v_2| \, \phi(v''_2) \, d\sigma_2 \, dv'_1.
\eeqn
Remarking that $v''_2 = v_2 + (\gamma/2) \, (|u_2| \, \sigma_2 - u_2)$ with $u_2 = v'_1 - v_2$, we observe that the integral term $ \Psi_\eps$ is very similar to the collision term $Q^+$ (here $v_2$ (resp. $v_1$, $\sigma_2$, $\gamma$, $v''_2$) plays the role of $v$ (resp. $v_1$, $\sigma$, $\beta$, $'v$) in the gain term) and therefore we may give a Carleman representation of $\Psi_\eps$. The same computations as performed in~\cite[Proposition 1.5]{MMII} yield
$$
\Psi_\eps = {4 \over \gamma^2}  \int_{\R^N} \! \int_{E_{v_2,v''_2}}  
{\bf 1}_{S_{v,v_1,\eps}}(v'_1) \, |v''_2 - v_2|^{-1} \, \phi(v''_2) \, dE(v''_3) \, dv''_2
$$
where $E_{v_2,v''_2}$ is the hyperplan orthogonal to the vector $v_2 - v''_2$ and passing through the 
point $\Omega_{v_2,v''_2} = v_2 + (1-\gamma^{-1}) \, (v_2 - v''_2)$. Here $v''_3$ stands for the post collision velocity issued from $v'_1$, that is $v''_3 = V'_* (v_2,v'_1,\sigma_2)$, and then, thanks to the momentum conservation, $v'_1 := v''_2 + v''_3 - v_2$.  We finally define $\Pi_{v_2,v''_2}$ the hyperplan orthogonal to the vector $v_2 - v''_2$ and passing through the 
point $\Omega'_{v_2,v''_2} = v''_2 + (1-\gamma^{-1}) \, (v_2 - v''_2)$ and we get 
\beqn\label{Psieps3}
\Psi_\eps = {4 \over \gamma^2}  \int_{\R^N} \! \int_{\Pi_{v_2,v''_2}}  
{\bf 1}_{S_{v,v_1,\eps}}(v'_1) \, |v''_2 - v_2|^{-1} \, \phi(v''_2) \, dE(v'_1) dv''_2.
\eeqn
Now, arguing as in \cite[lemma 2.1]{Abraham}, we see that  the measure of the intersection $\Gamma_\eps$ 
between the plane $\Pi_{v_2,v''_2}$  and the thickened sphere $S_{v,v_1,\eps}$ is bounded by $\pi \, \eps \, \gamma \, |v-v_1|$ and that  $v''_1 \in \Gamma_\eps$ implies that $v''_2 \in B^\eps$ with 
$$
B^\eps := \Big\{ z \in \R^N; |z|^2 \le |v|^2 + |v_1|^2 + 2 \, \eps \, (|v|+|v_1|) +  \eps^2 \, |v_2|^2 \Big\}.
$$
Gathering these estimates with (\ref{defPhi}) and (\ref{Psieps3}) we get 
\bean
\Phi &=& {(2/\gamma)^4 \over |v-v_1|}  \lim_{\eps \to 0}  {1 \over \eps} 
 \int_{\R^N}  { \phi(v''_2) \over |v''_2 - v_2| } \, \hbox{mes}(\Gamma_\eps) \,  dv''_2 \\
&\le&  {2^4 \, \pi \over \gamma^3} \lim_{\eps \to 0}   \int_{\R^N}  
{ \phi(v''_2) \over |v''_2 - v_2|} \, {\bf 1}_{B^\eps }(v''_2)\,  dv''_2
=  {2^4 \, \pi \over \gamma^3}  \int_{\R^N}  
{ \phi(v''_2) \over |v''_2 - v_2|} \, {\bf 1}_{B^0 } (v''_2)\,  dv''_2
\eean
where we have defined  $B^0 :=  \{ z \in \R^N; |z|^2 \le |v|^2 + |v_1|^2 \}$. Using \cite[lemma 2.2]{Abraham} we may conclude as in the end of \cite[lemma 2.1]{Abraham}  and therefore (\ref{L3estim}) follows. 
\qed

\smallskip

We second establish that the solution $g$ of the rescaled equation (\ref{eqresca}) decomposes 
between a regular part and a small remaining part as it has been proved for the elastic 
Boltzmann equation in \cite{MV**}, and then partially extended to the inelastic Boltzmann 
equation in \cite{MMII}. As compared to this last paper, this result relaxes the assumption on the 
initial datum to $g_{\mbox{\scriptsize{in}}} \in L^1_3$, but at the price of the hypothesis of a lower bound 
on the energy.

\smallskip

 \begin{lem} \label{decompos2} Consider $g_{\mbox{\scriptsize{{\em in}}}} \in L^1_3$ and the 
 associated solution $g \in C([0,\infty);L^1_3)$ to the rescaled equation (\ref{eqresca}). 
 Assume that for some constant $\rho, c_1,M_1,T \in (0,\infty)$ there holds
 \beqn\label{assumgdecomp1}
 g_{\mbox{\scriptsize{{\em in}}}} \in \Cc_{\rho,0}, \qquad 
 \| g_{\mbox{\scriptsize{{\em in}}}} \|_{L^1_3} \le M_1, \qquad \forall \, t \in [0,T], \quad \EE(g(t,\cdot)) \ge c_1.
 \eeqn
Then, there are $\alpha_7 \in [\alpha_6,1)$ and  $\lambda \in (-\infty,0)$, and for any exponential weight function $m$ (as defined in (\ref{defdem}) and any $k \in \N$, there exists a constant $K$ (which depends on $\rho, c_1,M_1,k,m$) 
such that for any $\alpha \in [\alpha_7,1]$, we may split $g = g^S + g^R$ with 
\beqn\label{concludegdecomp1bis}
\forall \, t \in [0,T], \quad  \| g^S(t, \cdot) \|_{H^k \cap L^1(m^{-1})} \le K, \qquad 
\| g^R(t, \cdot) \|_{ L^1_3} \le K \, e^{\lambda \, t}.
\eeqn
 \end{lem}

\smallskip\noindent
{\sl Proof of Lemma~\ref{decompos2}} 
The starting point is to write the rescaled equation (\ref{eqresca}) in the following way
$$ %\label{decomproof1}
{\partial g \over \partial t}  + \tau_\alpha \, v \cdot \nabla_v g + \ell\, g = Q^+_\alpha(g,g),
$$
with $\ell(t,v) := \tau_\alpha \, N + L(g(t,\cdot))(v)$. Introducing the linear semigroup 
$$
(U_t \, h) (v) =  h(v \, e^{-\tau_\alpha \, t}) \, \exp \left( - \int_0^t \ell(s,v) \, ds\right)
$$
and using the Duhamel formula, we have
\[ 
g_t = U_t g_{\mbox{\scriptsize{in}}} + \int_0 ^t U_{t-s} Q_\alpha^+(g_s,g_s) \, ds.
 \]
We iterate that last identity and we obtain $g= g^R_1 + g^S_1$ with 
\[
 g^R_1= U_t g_{\mbox{\scriptsize{in}}} + \int_0 ^t U_{t-s} Q_\alpha^+(g_s,U_s g_{\mbox{\scriptsize{in}}} ) \, ds, \quad
g^S_1= \int_0 ^t \!\! \int_0 ^s U_{t-s} Q_\alpha^+(g_s,  U_{s-u} Q_\alpha^+(g_u,g_u) ) \,du \, ds.
\]
On the one hand, the energy lower bound (\ref{assumgdecomp1}) and Lemma~\ref{estimatesonL} imply that there exists a constant  $c_2 \in (0,\infty)$ such that 
$$
(U_t \, h) (v) \le e^{-c_2 \, t} \, (V_{\xi_t} h)(v)
\quad\hbox{with}\quad (V_\xi h)(v) = h(\xi \, v) \quad\hbox{and}\quad \xi_t = e^{-\tau_\alpha \, t}
$$
On the other hand, straightforward homogeneity arguments leads to 
$$
Q^+_\alpha(g,V_\xi h) = \xi^{-N-1} \, V_{\xi^{-1}} \, Q^+_\alpha(V_{\xi^{-1}} g,h)
$$ 
and $\| h_\xi \, |.|^q\|_{L^p} = \xi^{-q-N/p} \, \| h \, |.|^q\|_{L^p}$ for any functions $g,h$ and 
positive real $\xi$. From these considerations we deduce that 
$$
\| g_1^R(t) \|_{L^1} 
\le e^{(N \, \tau_\alpha-c_2) \, t} \, \| g_{\mbox{\scriptsize{in}}} \|_{L^1} + 
e^{((N+1) \, \tau_\alpha-c_2) \, t} \, \| g_{\mbox{\scriptsize{in}}} \|_{L^1_1} \, \sup_{s\ge0} \, \| g_s \|_{L^1_1} 
\le C \, e^{-(c_2/2) \, t},
$$
for some constant $C$ and for any $(1-\alpha)$ small enough. In the same way, we have 
\bean
\| g_1^S(t) \|_{L^3} 
%&\le& 
%\int_0 ^t \!\! \int_0 ^s e^{-c_2 \, (t-\sigma)} \, V_{t-s} Q_\alpha^+(g_s,  V_{s-\sigma} Q_\alpha^+(g_\sigma,g_\sigma) )  \|_{L^3}  \,d\sigma ds \\
&\le& 
\int_0 ^t \!\! \int_0 ^s e^{[(2N/3+1)\, \tau_\alpha-c_2] \, (t-\sigma)} \, \| Q_\alpha^+(V_{\xi^{-1}_{s-\sigma}}\, g_s, Q_\alpha^+(g_\sigma,g_\sigma) )  \|_{L^3}  \,d\sigma ds.
\eean
Taking $(1-\alpha)$ smaller if necessary and using Lemma~\ref{itarteQ+}, we obtain 
\bean
\| g_1^S(t) \|_{L^3} 
&\le& 
\left( \int_0 ^t \!\! \int_0 ^s e^{ -(c_2/2)\, (t-\sigma)} \, d\sigma \, ds \right)
\sup_{s \ge 0} \| g_s \|_{L^1 _2}^{3} ,
\eean
which ends the proof of (\ref{concludegdecomp1bis}) in the case $k= 0$, with the help of point (i) of 
Lemma~\ref{estimMSg}. The general case $k \in \N^*$ is then treated by following the strategy 
introduced in \cite{MV**} and using the result of appearance of regularity proved in \cite{MMII} 
(and recalled in point (iii) of Lemma~\ref{estimMSg}).
\qed
\smallskip

We third recall a classical $L^1$ stability result for the elastic Boltzmann equation which has been established in \cite[Proposition 3.2]{MMII} for the rescaled equation (\ref{eqresca}). 

\begin{lem} \label{stablocal} 
Consider $0 \le g_{\mbox{\scriptsize{{\em in}}}}^1,  g_{\mbox{\scriptsize{{\em in}}}}^2 \in L^1_3 \cap \Cc_{\rho,0}$ 
and the two associated solution $g_i  \in C([0,\infty);L^1_3) \cap L^\infty(0,\infty; L^1_3)$ to the rescaled 
equation (\ref{eqresca}). There exists $C_{\mbox{\scriptsize{{\em stab}}}} \in (0,\infty)$ (only depending 
on $b$ and $\sup_{t\ge 0} \|g^1 + g^2\|_{L^1 _3}$) such that
\[ \forall \, t \ge 0, \ \ \
\|g^2_t -g^1_t \|_{L^1 _2} \le  \|g^2_{\mbox{\scriptsize{{\em in}}}} - g^1_{\mbox{\scriptsize
{{\em in}}}} \|_{L^1 _2} \,
e^{C_{\mbox{\scriptsize{{\em stab}}}} t}. \]
\end{lem}

\medskip\noindent{\sl Proof of point (iv) of Theorem~\ref{theo:uniq}. } 
Let us consider $g_{\mbox{\scriptsize{in}}} \in L^1_3 \cap \Cc_{\rho,0,\EE_{\mbox{\scriptsize{in}}}}$ 
with $\| g_{\mbox{\scriptsize{in}}} \|_{L^1_3} \le M_0$ for some fixed $\EE_{\mbox{\scriptsize{in}}}, \, M_0 \in (0,\infty)$ 
and $g$ the associated solution to the rescaled equation (\ref{eqresca}) which has been built in \cite{MMII}. 
We know that there exists $M_1 \in (0,\infty)$ such that 
\beqn\label{bdL3}
\sup_{(0,\infty)} \| g(t,\cdot) \|_{L^1_3} \le M_1.
\eeqn

\medskip\noindent{\sl Step 1.  } We define 
$$ 
T_* := \sup\big\{ T \in (0,\infty), \,\, \EE(g(t, \cdot)) \ge c_1 \,\, \forall \, t \in [0,T] \big\}, \quad 
c_1 := \min \{ \EE_{\mbox{\scriptsize{in}}}, \bar \EE_1 \} /2. 
$$
We shall prove that  $T_* =+ \infty$. We argue by contradiction, assuming that $T_* < \infty$. From the equation on the energy (\ref{eqEg}) and the uniform estimate (\ref{bdL3}) and from the definition of $T_*$ we have 
\beqn\label{T*&EE'}
T_* \ge C_1 \, (1-\alpha)^{-1} \quad\hbox{and}\quad \EE'(T_*) \le 0.
\eeqn
Thanks to Lemma~\ref{decompos2}, we may decompose 
$$
g =  g^S + g^R \quad\hbox{on}\quad (0,t_1),
$$
with $t_1 \in (0, T_*)$ to be fixed. At time $t_1$ we initiate a new flow 
starting from the smooth part of $g$. More precisely,  we decompose
$$
g = \tilde g^S + \tilde g^R  \quad\hbox{on}\quad (t_1,T_*),
$$
with $\tilde g^S(t_1) = [\rho / \rho(g^S(t_1))] \, g^S(t_1)$, $\tilde g^S$ solution (with mass $\rho$!) 
to the equation (\ref{eqresca}) on $(t_1,T_*)$ and $\tilde g^R := 
g - \tilde g^S$. On the one hand, from (\ref{concludegdecomp1bis}) and Lemma \ref{stablocal} we have
$$
\| \tilde g^R (t) \|_{L^1_3}  \le C \, e^{C_{\mbox{\scriptsize{stab}}} \, (T_*-t_1) + \lambda \, t_1} \quad
\hbox{on}\quad (t_1,T_*).
$$
We choose $t_1 = \eta \, T_*$ with $\eta \in (0,1)$ in such a way that 
$C_{\mbox{\scriptsize{stab}}} \, (1-\eta) + \lambda \, \eta = \lambda/2$. We have then proved  
\beqn\label{bdgrL3}
\| \tilde g^R (t) ]  \|_{L^1_3} \le C \, e^{(\lambda/2) \, C_1 \,  (1-\alpha)^{-1}} \quad
\hbox{on}\quad (t_1,T_*).
\eeqn
On the other hand, following Step 3 in the proof of Proposition \ref{Attractiveg}, 
we deduce a similar estimate as (\ref{bdRH}), namely
\beqn\label{bdRHgS}
\| \tilde g^S(T_*, \cdot) - M[\tilde g^S(T_*, \cdot)] \|_{L^1_3} = \OO((1-\alpha)^{1/8}) 
\eeqn
for any $(1-\alpha)$ small enough chosen in such a way that the intermediate time 
$t_2$ defined in Step 3 of the proof of Proposition~\ref{Attractiveg} satisfies $t_1 +t_2 \le T_*$. 
Gathering (\ref{bdgrL3}) and (\ref{bdRHgS}) we obtain
$$ %\label{bdRHg}
\| \tilde g(T_*, \cdot) - M[g(T_*, \cdot)] \|_{L^1_3} = \OO((1-\alpha)^{1/8} ).
$$
Coming back to the equation (\ref{eqEg}) on the energy and proceeding like in Step 4 in the proof of 
Proposition~\ref{Attractiveg}, we get 
$$
\EE'(T_*) \ge (1-\alpha) \, \Big[ k_3 \, c_1 \, (\bar \EE_1^{1/2} - c_1^{1/2}) - C \, (1-\alpha)^{1/8} \Big] > 0
$$
for any $(1-\alpha)$ small enough. That is in contradiction with (\ref{T*&EE'}) and we conclude 
that $T_*=+\infty$.

\medskip\noindent{\sl Step 2.  } 
Thanks to the previous step, we have a uniform in time lower bound 
on the energy, and therefore we can run the decomposition theorem 
for all times. 

By applying the decomposition theorem as in Step~1 for a given 
time $t \in (0,\infty)$, starting a new flow at $t_1 = \eta \, t$ taking 
$[\rho / \rho(g^S(t_1))] \, g^S(t_1,\cdot)$ 
as initial datum, and then using Corollary~\ref{coro:cvgHk} on the smooth part 
$\tilde g^S(s,\cdot)$, $s \in [t_1,t]$, we find that at time $t$, the solution $g_t$ decomposes 
as $\tilde g^S _t + \tilde g^R _t$, where $\tilde g^S _t$ approaches the self-similar profile 
with rate $C \, e^{\mu_\alpha \, (t-t_1)}$, that is $C \, e^{(1-\eta) \, \mu_\alpha \, t}$, 
and the remaining part $\tilde g^R$ goes to $0$ with rate $C \, e^{(\lambda/2) \, t}$.  
Since $|\lambda/2|$ is larger than $(1-\eta) \, |\mu_\alpha|$ for $(1-\alpha)$ small 
enough, it concludes the proof of (\ref{cvgirreg}). 
\qed

\appendix

\section{Appendix: Moments of Gaussians}
\setcounter{equation}{0}
\setcounter{theo}{0}

We state here some moments of tensor product of Gaussians.

\begin{lem}\label{MomentsGausians} The following identities hold
\bear
\label{Mv2}
&&\int_{\R^N} M_{1,0,1} \, |v|^2 \, dv = N,  \\
\label{Mv4}
&& \int_{\R^N} M_{1,0,1} \, |v|^4 \, dv = N \, (N+2), \\
\label{MMu3}
&&\int_{\R^N \times \R^N} M_{1,0,1} \, (M_{1,0,1})_* \, |u|^3 \, dv \, dv_* = 
     2^{3/2} \, \int_{\R^N} M_{1,0,1} \, |v|^3 \, dv, \\ \label{MMv2u3}
&&\int_{\R^N \times \R^N} M_{1,0,1} \, (M_{1,0,1})_* \,|v|^2 \,  |u|^3 \, dv \, dv_* = 
 \sqrt{2} \, (2N+3) \, \int_{\R^N} M_{1,0,1}(v) \, |v|^3 \, dv.
\eear    
\end{lem}

\smallskip\noindent{\sl Proof of Lemma~\ref{MomentsGausians}. }
The proof of (\ref{Mv2}) and (\ref{Mv4}) being straightforward and
the proof of (\ref{MMu3}) being very similar to the proof of (\ref{MMv2u3}) we
only prove  (\ref{MMv2u3}). We first notice that 
  \bean 
  && \int_{\R^N \times \R^N} M_{1,0,1} \, (M_{1,0,1})_* \, |v|^2 \, |u|^3 \, dv \, dv_* \\
  && \qquad = \frac{1}{2} \, \int_{\R^N \times \R^N} M_{1,0,1} \, (M_{1,0,1})_* \, 
     (|v|^2 +|v_*|^2) \, |u|^3 \, dv \, dv_* \\
  && \qquad = \frac{1}{4} \, \int_{\R^N \times \R^N} M_{1,0,1} \, (M_{1,0,1})_* \, 
     (|v+v_*|^2 +|v-v_*|^2) \, |u|^3 \, dv \, dv_*. 
  \eean
Making use of the change of variable $(v,v_*) \to (x=(v+v_*)/\sqrt{2},y=(v-v_*)/\sqrt{2})$, we then   get 
  \bean 
  && \int_{\R^N \times \R^N} M_{1,0,1} \, (M_{1,0,1})_* \, |v|^2 \, |u|^3 \, dv \, dv_* \\
  && \qquad = \sqrt{2} \, \int_{\R^N \times \R^N} M_{1,0,1}(x) \, M_{1,0,1}(y) \, 
     (|x|^2 +|y|^2) \, |y|^3 \, dx \, dy \\ 
  && \qquad = \sqrt{2} \, N \, \int_{\R^N} M_{1,0,1}(v) \, |v|^3 \, dv 
              + \sqrt{2} \, \int_{\R^N} M_{1,0,1}(v) \, |v|^5 \, dv \\
  && \qquad = \sqrt{2} \, (2N+3) \, \int_{\R^N} M_{1,0,1}(v) \, |v|^3 \, dv.
  \eean
\qed

\section{Appendix: Interpolation inequalities}
\setcounter{equation}{0}
\setcounter{theo}{0}

\begin{lem} \label{interpolineg} 
\begin{itemize}
\item[(i)] 
For any $k,k^*,q,q^* \in \ZZ$ with $k\ge k^*$, $q\ge q^*$ and any $\theta \in (0,1)$ there is  $C \in (0,\infty)$ such that for 
$h \in W^{k^{**},1}_{q^{**}}(m^{-1})$
\beqn\label{hWk1qm}
\| h \|_{W^{k,1}_q(m^{-1})} \le C \,\| h \|_{W^{k^*,1}_{q^*}(m^{-1})}^{1-\theta} \,  \| h \|_{W^{k^{**},1}_{q^{**}}(m^{-1})}^\theta.
\eeqn
with $k^{**},q^{**} \in \ZZ$ such that $k = (1-\theta) \, k^{*} + \theta \, k^{**}$, $q = (1-\theta) \, q^{*} + \theta \, q^{**}$.
\item[(ii)] 
For any $k, q \in N^*$ and any exponential weight function $m$ as defined in (\ref{defdem}), 
there exists $C \in (0,\infty)$ such that for any $h \in H^{k^\ddagger} \cap L^1(m^{-12})$ 
with $k^\ddagger := 8 k + 7 (1+N/2)$ 
\beqn\label{hWk1qm2}
\| h \|_{W^{k,1}_q(m^{-1})} \le C \, \| h \|_{H^{k^\ddagger}}^{1/4} \, \,\| h \|_{L^1( m^{-12})}^{1/4} \,  \| h \|_{L^1(m^{-1})}^{3/4}.
\eeqn
\end{itemize}
\end{lem}

\medskip\noindent{\sl Proof of Lemma~\ref{interpolineg}.}
The inequality (\ref{hWk1qm}) in point (i) is a classical result from interpolation theory. Let us focus on 
point (ii). We prove the inequality (\ref{hWk1qm2}) for $h \in \Ss(\R^N)$ and then argue by density. 
On the one hand, we observe that for any $\ell$ there exists $C$ such that 
$$
\| h \|_{H^\ell}^2 \le  C \, \| h \|_{L^1} \, \| h \|_{H^{\ell^\dagger} }, \quad \ell^\dagger := 2\ell + 1+ N/2.
$$
Iterating twice that inequality, we get (for some related exponents $k^\dagger, k^\ddagger$)
\beqn\label{HkL1}
\| h \|_{H^{k^\dagger}} \le C \, \| h \|_{L^1}^{3/4} \, \| h \|_{H^{k^\ddagger} }^{1/4}.
\eeqn

On the other hand, using first Cauchy-Schwartz inequality, plus the same argument as 
above and H\"older's inequality, we obtain
\bear\nonumber
\| h \|_{W^{k,1}_q(m^{-1})} &\le& C \| h \|_{H^k(m^{-3/2})} \le C \, \| h \|_{L^1(m^{-3})}^{1/2} \, \| h \|_{H^{k^\dagger} }^{1/2} \\
\label{WkL1}
 &\le&  C \, \| h \|_{L^1(m^{-12})}^{1/8} \, \| h \|_{L^1}^{3/8} \, \| h \|_{H^{k^\dagger} }^{1/2}.
\eear
We conclude gathering (\ref{WkL1}) and (\ref{HkL1}).
\qed

\bigskip
\noindent
{\bf{Acknowledgments.}} We thank Alexander Bobylev, Jos\'e Antonio Carrillo and  C\'edric Villani
for fruitful discussions on the inelastic Boltzmann equation.

\end{document}